\documentclass[preprint,12pt]{article}
\usepackage{amssymb}
\usepackage{mathrsfs}
\usepackage{amsfonts}
\usepackage{graphicx}
\usepackage{colortbl,dcolumn}
\usepackage{amsmath}
\usepackage{psfrag}
\usepackage{booktabs}
\usepackage{cite}

\usepackage{amsthm}
\numberwithin{equation}{section}
\usepackage[top=0.8in, bottom=0.8in, left=0.8in, right=0.8in]{geometry}

\newcommand{\R}{\mathbb{R}}
\newcommand{\N}{\mathbb{N}}
\newcommand{\E}{\mathbb{E}}
\renewcommand{\P}{\mathbb{P}}
\newcommand{\tr}{\operatorname{trace}}

\newcommand{\dd}{\text{d}}

\newtheorem{thm}{Theorem}[section]

\newtheorem{lem}[thm]{Lemma}
\newtheorem{prop}[thm]{Proposition}
\newtheorem{cor}[thm]{Corollary}
\newtheorem{rem}[thm]{Remark}

\newtheorem{assumption}[thm]{Assumption}

\begin{document}
\title{Mean-square convergence rates of implicit Milstein type methods for SDEs with non-Lipschitz coefficients
%: 
%applications to financial models
%
%New implicit Milstein type methods for SDEs with non-Lipschitz coefficients: applications to financial models
%Implicit Milstein type methods for SDEs with non-globally Lipschitz coefficients
%and their applications of approximating financial models
%Mean-square convergence rates of implicit Milstein methods for SDEs with super-linearly growing coefficients
%and their applications to stochastic volatility processes
 \footnotemark[1]}

\author{
Xiaojie Wang\footnotemark[2]
%,  
%\, Mengchao Wang \footnotemark[3], \, Yuying Zhao \footnotemark[4]
%\quad
%Ruisheng Qi$\,^\text{b}$
%\quad
%Fengze Jiang$\,^\text{c}$
\\
\footnotesize School of Mathematics and Statistics, HNP-LAMA, Central South University, Changsha, China\\
%\footnotesize x.j.wang7@csu.edu.cn\; and \;x.j.wang7@gmail.com\\
%\footnotesize $\,^\text{b}$ School of Mathematics and Statistics, Northeastern University at Qinhuangdao, Qinhuangdao, China\\
%\footnotesize qirsh@neuq.edu.cn\; and \;qiruisheng123@sohu.com
%\footnotesize $\,^\text{c}$ School of Mathematics and Statistics, Huazhong University of Science and Technology, Wuhan, China\\
}

%      \date{}
       \maketitle

%       \footnotetext{\footnotemark[1] This work was supported by Natural Science Foundation of China (12071488, 11971488)
%                and Natural Science Foundation of Hunan Province for Distinguished Young Scholars (2020JJ2040). 
%                The author also wants to thank Dr. Xu Yang for his useful comments based on carefully reading the draft.
%                }
        \footnotetext{\footnotemark[2] x.j.wang7@csu.edu.cn, x.j.wang7@gmail.com}
%        \footnotetext{\footnotemark[3] ...@...}
%        \footnotetext{\footnotemark[4] ...@...}
        % Corresponding author:
       \begin{abstract}
          {\rm\small
            A class of implicit Milstein type methods is introduced and analyzed in the present article for
           stochastic differential equations (SDEs) with non-globally Lipschitz drift and diffusion coefficients.
           By incorporating a pair of method parameters $\theta, \eta \in [0, 1]$ into 
           both the drift and diffusion parts, the new schemes are indeed a kind of drift-diffusion double implicit methods.
           Within a general framework, we offer upper mean-square error bounds for the proposed schemes,
           based on certain error terms only getting involved with the exact solution processes.
           Such error bounds help us to easily analyze mean-square convergence rates of
           the schemes, without relying on a priori high-order moment estimates of numerical approximations.
           Putting further globally polynomial growth condition, we successfully recover the expected mean-square convergence rate of order one
           for the considered schemes with $\theta \in [\tfrac12, 1], \eta \in [0, 1]$.
           Also, some of the proposed schemes are applied to 
           solve three SDE models evolving in the positive domain $(0, \infty)$.
           More specifically,
           the particular drift-diffusion implicit Milstein method ($ \theta = \eta = 1 $) is utilized to
           approximate the Heston $\tfrac32$-volatility model and the stochastic Lotka-Volterra competition model. 
           The semi-implicit Milstein method ($\theta =1, \eta = 0$)
           is used to solve the Ait-Sahalia interest rate model.
           Thanks to the previously obtained error bounds, we reveal the optimal mean-square convergence rate of
           the positivity preserving schemes under more relaxed conditions, compared with existing relevant results
           in the literature. Numerical examples are also reported to confirm the previous findings.
           } \\

%\textbf{PACS: 02.60.Lj}

\end{abstract}
%%%%%%%
%%%%%%

%%%%%
%%%%%
%
\section{Introduction}	
%
%
%As is known to all,
%It is well-known that 
Stochastic differential equations (SDEs) find applications in  a wide range of scientific areas such as
finance, chemistry, biology, engineering and many other branches of science.
In general, analytical solutions to nonlinear SDEs are usually not available and
%exact solutions are not explicitly available and 
development and analysis of numerical methods for simulation of SDEs are of significant interest in practice.
%many applications.
%one usually relies on their numerical approximations in applications. 
To analyze the numerical approximations,
a global Lipschitz condition is often imposed on the coefficient functions of SDEs \cite{milstein2013stochastic,kloeden1992numerical}.
% the standard assumption
Nevertheless,  SDEs arising from applications rarely obey such a traditional but restrictive condition.
Notable examples of SDEs with non-globally Lipschitz continuous coefficients 
include numerous models such as the $\tfrac32$-volatility model \cite{heston1997simple,lewis2000option},
\begin{equation}
\label{eq:intro-32model}
d X_t = X_t ( \mu - \alpha X_t ) dt + \beta X_t ^ {3/2} \, \dd W_t,
 \quad 
 X_0 = x_0 > 0,
 \quad
 \mu , \alpha , \beta > 0,
\end{equation}
and the Ait-Sahalia interest rate model \cite{ait-sahalia1996testing},
\begin{equation}
\label{eq:intro-Ait-Sahalia-model-SDE}
\dd X_t = ( \alpha_{-1} X_t^{-1} - \alpha_0 + \alpha_1 X_t - \alpha_2 X_t^\kappa  ) \, \dd t + \sigma X_t^{\rho} \, \dd W_t,
\quad
X_0 = x_0 > 0,
\end{equation}
from mathematical finance,
where $\alpha_{-1}, \alpha_0, \alpha_1, \alpha_2, \sigma > 0$ are positive constants and $\kappa > 1, \rho >1$.
Evidently, coefficients of these two models violate the global Lipschitz condition.
%
%This results in essential ...
As already shown in \cite{hutzenthaler2011strong}, the popularly used Euler-Maruyama method
produces divergent numerical approximations when used to solve a large class of SDEs with super-linearly growing coefficients,
such as \eqref{eq:intro-32model} and \eqref{eq:intro-Ait-Sahalia-model-SDE}.
%including the above two SDE examples.
Therefore special care must be taken to design and analyze convergent numerical schemes
in the absence of the Lipschitz regularity of coefficients.
%a non-globally Lipschitz setting.
Recent years have witnessed a prosper growth of relevant works devoted to the numerical analysis of 
SDEs under non-globally Lipschitz conditions,  
with an emphasis on analyzing implicit schemes \cite{higham2002strong,zong2018convergence,higham2013convergence,
beyn2016stochastic,beyn2017stochastic,andersson2017mean,neuenkirch2014first,alfonsi2013strong,mao2013strong,
mao2013strongJCAM,wang2018mean,yao2018stability,cui2019convergence,hu96semi-implicit},
and devising explicit methods based on modifications of traditionally explicit schemes \cite{guo2018truncated,Hutzenthaler12,Hutzenthaler15numerical,mao2015truncated,mao2016convergence,
zhang2017order-preserving,Szpruch18V-integrability,Tretyakov2013fundamental,sabanis2013note,sabanis2016euler,wang2013tamed,
hutzenthaler2018exponential,
chassagneux2016explicit,kumar2019milstein,gan2020tamed,hutzenthaler2020perturbation,
kelly2019strong,kelly2017adaptive,fang2020adaptive},
%and investigating adaptive time-stepping schemes \cite{kelly2019strong}.
%
to just mention a few. 
%and one can get a comp... list of relevant references therein... 
%an extended literature review on this direction is available in \cite{Hutzenthaler15numerical}.
Although explicit methods such as tamed methods \cite{Hutzenthaler12,sabanis2016euler} 
and truncated schemes \cite{mao2015truncated,guo2018truncated},  computationally more efficient than
implicit ones for one time step, are able to well tackle non-stiff SDEs with super-linearly growing coefficients, they 
usually face a severe stepsize restriction due to stability issues when used to solve stiff SDE systems \cite{milstein2013stochastic}. 
Moreover, explicit time stepping schemes like tamed methods, similarly to the classical explicit Euler/Milstein schemes,  
are usually not positivity preserving when applied to approximate financial models 
whose solutions naturally remain positive (see, e.g., \cite{kahl2008structure,higham2013convergence,szpruch2011numerical}).
%
%%%%%%%%%%%% small noise %%%%%%%%
%\cite{buckwar2006multistep}
%\cite{romisch2006stepsize}
%\cite{buckwar2010stochastic}
%\cite{milstein1997numerical}
%\cite{anderson2016multilevel}
%\cite{milstein1997mean}

%%%%%%%%%%%Monograph%%%%%%
%\cite{milstein2013stochastic}
%\cite{kloeden1992numerical}
%
%
%On the basis of ...
%
%

%Both convergence and stability of the stochastic theta methods (STMs) with $\theta \in [0, 1]$ has been extensively
%studied in the literature \cite{...}.

In this article we are concerned with implicit Milstein schemes for mean-square approximations of It\^o  SDEs 
with non-globally Lipschitz continuous coefficients, in the form of
\begin{equation} \label{eq:intro-SODE}
%\begin{split}
%& 
\dd X_t = f ( X_t ) \, \dd t + g ( X_t ) \, \dd W _ t, 
\quad
t \in (0, T],
%\\
%& 
\quad X_0 = x_0,
%\end{split}
\end{equation}
where  $W \colon [0,T] \times \Omega \rightarrow \mathbb{R}^m $ stands for the $\mathbb{R}^m$-valued 
standard Brownian motion, $f \colon \mathbb{R}^d \rightarrow \mathbb{R}^d$ the drift coefficient function, and 
$g \colon \mathbb{R}^d \rightarrow \mathbb{R}^{d \times m}$ the diffusion coefficient function.
Mean-square approximations are of particular importance for the computation of 
statistical quantities of the solution process of \eqref{eq:intro-SODE}
through computationally efficient multilevel Monte Carlo (MLMC) methods  \cite{giles08multilevel}.
%This work attempts to seek new implicit Milstein-type schemes for \eqref{eq:intro-SODE},
Recall that  Milstein-type schemes achieve a higher mean-square convergence rate than the Euler-type schemes and can be combined with 
the MLMC approach to reduce computational costs further \cite{giles08multilevel,giles19analysis,giles2008improved}.
In the literature,  various Milstein type methods \cite{beyn2017stochastic, Buckwar2011comparative, guo2018truncated,higham2013convergence,kelly2019strong,kumar2019milstein,
wang2013tamed,WGW12,zhang2017order-preserving,li2020explicit,kruse2019randomized,bossy2015strong} have been studied
and the present work proposes a class of implicit Milstein-type schemes and establish a mean-square convergence theory for the new schemes.
%justifies an efficient Multilevel Monte Carlo method  for SDEs with non-Lipschitz coefficients.
On a uniform mesh constructed over $[0, T]$ with a uniform time stepsize $h = \tfrac{T}{N}, N \in \N$, 
we develop a family of double implicit Milstein-type methods with a pair of method parameters $(\theta, \eta)$
for \eqref{eq:intro-SODE} as follows:
\begin{equation}
\label{eq:introduction-scheme}
 \begin{split}
 Y_{ n + 1} = & Y_{ n } + \theta f ( Y_{n+1} )h
         + ( 1-\theta ) f( Y_n )h
            + g(Y_{n}) \Delta W_n
            + \sum^m_{j_1,j_2=1}
            \mathcal{L}^{j_1} g_{j_2}(Y_{n}) I_{j_1,j_2}^{t_n,t_{n+1}}
            \\
            & \quad
            +
             \tfrac {\eta} {2} \sum_{ j = 1}^{ m} \mathcal{L}^{j } g_{j }(Y_{n}) h
            -
            \tfrac {\eta} {2} \sum_{ j = 1}^{ m} \mathcal{L}^{j } g_{j }(Y_{n+1}) h,
            \qquad
            Y_{0}= X_0,
\end{split}
\end{equation}
where  $\theta, \eta \in [0, 1]$, 
$\Delta W_{ n } : = W_{t_{ n+1 } } - W_{t_ {n} } $, $ n \in \{0, 1, 2, \ldots, N-1\}$ 
and $\mathcal{L}^{j_1} g_{j_2}, I_{j_1,j_2}^{t_n,t_{n+1}}$ are precisely defined by \eqref{eq:scheme-notation}.
%
%aims to examine implicit Milstein type schemes for
%
When $d = m = 1$, the schemes \eqref{eq:introduction-scheme} coincide with the proposed ones in \cite{higham2013convergence},  
where the authors used the positivity preserving schemes to 
solve the $\tfrac32$-volatility model \eqref{eq:intro-32model} and proved its strong convergence 
with no convergence rate revealed.  
%@@Related results on Milstein type schemes...
After assigning $\eta = 0$, the proposed scheme reduces to the classical $\theta$-Milstein method,
which has been studied in \cite{kloeden1992numerical,zong2018convergence,Buckwar2011comparative}.
But in the regime of possibly super-linearly growing diffusion coefficients $g$, 
the strong convergence rate of the $\theta$-Milstein method is, up to the best of our knowledge, still an open problem.
This paper shall fill the gap.

Also, we mention that 
%the underlying schemes work for multi-dimensional SDEs with non-commutative noise,  but 
an order reduction would be caused due to additional costs of approximating multiple stochastic integrals
%the computational efficiency will be reduced due to additional efforts spent for approximating multiple stochastic integrals 
$I_{j_1,j_2}^{t_n,t_{n+1}}$ when the multi-dimensional SDEs are driven by non-commutative noise.
As clarified in \cite[section 7]{rossler2010runge},  the effective order of Milstein methods
should be $\tfrac23$ in the case of non-commutative noise when the multiple stochastic integrals are efficiently approximated.
We refer to the simulation method proposed by Wiktorsson \cite{wiktorsson2001joint} 
and see also \cite{gilsing2007sdelab} for implementation issues.  
Compared to the order 0.5 strong Euler-type schemes, which attains
the effective order $0.5$, there is still a significantly improved convergence 
for the Milstein methods in the non-commutative noise setting.
Furthermore, we mention that the application of fully implicit (stochastically implicit) methods are unavoidable 
for stiff systems where the stochastic part plays an essential role (see \cite[pp.33]{milstein2013stochastic}
and \cite{milstein1998balanced} for detailed comments and an illustrative example).
By incorporating a pair of method parameters $\theta, \eta \in [0, 1]$ into the drift and diffusion parts, here we construct
a kind of fully implicit methods for general multi-dimension SDE systems with non-Lipschitz coefficients. 
Finally, we point out that proving the expected convergence rate of the proposed schemes for SDEs 
in non-Lipschitz settings, especially for the above two financial models, is highly non-trivial and  remains an unsolved problem.
The present work aims to fill these gaps by successfully establishing a first order of mean-square convergence 
for the scheme \eqref{eq:introduction-scheme} in different settings, covering the two aforementioned financial models.

%The present work successfully establish 
%a first order of mean-square convergence for the scheme \eqref{eq:introduction-scheme}
%within a general framework, which turns out to fill the above two gaps. 
%and also reveals a first order of mean-square convergence for 
%But their strong convergence rates are still open, which partly motivates
%the present work.
%@@ Many previous studies 
%          
%\cite{WGW12}                             %  Wang-Gan-Wang-implicit-Milstein 
%\cite{higham2013convergence} % Milstein-3/2-model
%\cite{Buckwar2011comparative} % Stability-Euler-Milstein

By formulating certain generalized monotonicity conditions in a domain $D \subset \R^d$
(Assumption \ref{ass:monotonicity-condition}), 
we develop an easy and novel approach to derive upper mean-square error bounds for the proposed schemes, 
which only get involved with the exact solution processes (Theorem \ref{thm:upper-error-bound}). 
The framework is broad and covers the two aforementioned SDE financial models.
Such error bounds are powerful as they help us to easily analyze mean-square convergence rates of
the schemes, without relying on a priori high-order moment estimates of numerical approximations.
Putting further globally polynomial growth and coercivity conditions in $\R^d$ (Assumption \ref{ass:f-polynomial-growth}), 
we utilize the derived upper error bound  to successfully identify a mean-square convergence rate 
of order one for the schemes \eqref{eq:introduction-scheme} solving general SDEs \eqref{eq:intro-SODE} 
(see Theorem \ref{thm: convergence-rate-of-STM}, Corollaries \ref{cor:MS-rate-Milstein1}, \ref{cor:MS-rate-Milstein2}).
%We recall that ...
% this paper is the first one to..., and also the first one to ...
%
%With the derived upper error bound at hand, it becomes very easy to recover mean-square convergence rates
%of the underlying scheme ... 
%
%

Later in section \ref{sect:applications},  we turn our attention to two scalar SDE models \eqref{eq:intro-32model} and 
\eqref{eq:intro-Ait-Sahalia-model-SDE} arising in mathematical finance and 
a stochastic Lotka-Volterra (LV) competitive model \eqref{eq:LVmodel} from ecology.
Since the considered models evolve in the positive domain $D = (0, \infty)$, instead of the whole space $\R$,
the convergence theory developed in section \ref{sec:conv-rates-global-polynomial}
cannot be applied in this situation. In order to address such issues, we apply two particular schemes 
covered by \eqref{eq:introduction-scheme} to approximate these specific models, 
which are capable of preserving positivity of the continuous models.
%i.e., 
%the $\tfrac32$-volatility model \eqref{eq:intro-32model} and 
%the Ait-Sahalia interest rate model \eqref{eq:intro-Ait-Sahalia-model-SDE} 
%mentioned above. 
More precisely,  the drift-diffusion double implicit Milstein method with parameters
$ \theta = \eta = 1 $ is utilized to approximate the Heston $\tfrac32$-volatility model \eqref{eq:intro-32model} 
and the stochastic LV competitive model \eqref{eq:LVmodel},
resulting in a recurrence of a quadratic equation with an explicit solution. 
%
%\begin{equation}
%\label{eq:intro-32model}
%d X_t = X_t ( \mu - \alpha X_t ) dt + \beta X_t ^ {3/2} \, \dd W(t),
% \quad 
% X_0 = x_0 > 0,
% \quad
% \mu , \alpha , \beta > 0,
%\end{equation}
%
And the semi-implicit Milstein method with a pair of parameters $\theta =1, \eta = 0$
is used to solve the Ait-Sahalia interest rate model \eqref{eq:intro-Ait-Sahalia-model-SDE}
in both a standard and a critical regime.
Both schemes are able to preserve positivity of the underlying models 
and their mean-square convergence rates are carefully analyzed.
%
%described by
%\begin{equation}
%\label{eq:intro-Ait-Sahalia-model-SDE}
%\dd X_t = ( \alpha_{-1} X_t^{-1} - \alpha_0 + \alpha_1 X_t - \alpha_2 X_t^\kappa  ) \, \dd t + \sigma X_t^{\rho} \, \dd W_t,
%\quad
%X_0 = x_0 > 0,
%\end{equation}
%where $\alpha_{-1}, \alpha_0, \alpha_1, \alpha_2, \sigma > 0$ are positive constants and $\kappa > 1, \rho >1$.
With the aid of the previously obtained error bounds, 
we prove a first order of mean-square convergence for both schemes  under mild assumptions
for the first time, which fills the gap left by  \cite{higham2013convergence,szpruch2011numerical}.
Compared with existing relevant results for first order schemes, more relaxed conditions are put here.
Specifically,  the drift-diffusion double implicit Milstein scheme is shown to achieve a mean-square convergence rate of order one
when used to solve the Heston $\tfrac32$-volatility model \eqref{eq:intro-32model} with model parameters 
obeying $ \tfrac{\alpha}{\beta^2} \geq \tfrac52  $ (Theorem \ref{thm:Milstein-convergence-rate-32model}).
Also, the semi-implicit Milstein method is proved to retain a mean-square convergence rate of order one, 
when solving the Ait-Sahalia interest rate model \eqref{eq:intro-Ait-Sahalia-model-SDE}, for full model parameters 
in the standard regime $ \kappa + 1 > 2 \rho $ (Theorem \ref{thm:Mistein-Ait-Sahalia-convergence-non-critcal}) 
and for model parameters obeying $\tfrac{ \alpha_2 } { \sigma^2 } \geq 2 \kappa - \tfrac32 $ 
and $\tfrac{ \alpha_2 }{ \sigma^2 } > \tfrac { \kappa + 1 } { 2 \sqrt{2} }$
in the general critical case $\kappa + 1 = 2 \rho$ (Theorem \ref{thm:Milstein-Ait-Sahalia-convergence-critcal}).

Recall that  a kind of Lamperti-backward Euler method was proposed and analyzed in \cite{neuenkirch2014first}  for 
a class of scalar SDEs defined in a domain, covering the above two financial models.
There a mean-square convergence rate of order one was proved for the scheme applied to 
the $\tfrac32$-volatility model  with parameters satisfying $ \tfrac{\alpha}{\beta^2} \geq 5 $ (see \cite[Proposition 3.2]{neuenkirch2014first}).
Also, the scheme used to approximate the Ait-Sahalia interest rate model owns a first mean-square convergence order 
for full model parameters in the case $ \kappa + 1 > 2 \rho $
and for parameters obeying $\tfrac{ \alpha_2 } { \sigma^2 } > 5 $ in a special critical case $\kappa = 2, \rho = 1.5$ 
(see Propositions 3.5, 3.6 from \cite{neuenkirch2014first}). 
Unlike the Lamperti transformed scheme introduced in \cite{neuenkirch2014first}, 
we propose and analyze the implicit Milstein-type schemes  applied to SDEs directly.
%Comparing our findings with relevant results in \cite{neuenkirch2014first}, 
From the above discussions, one can easily detect that our convergence results improve relevant ones in \cite{neuenkirch2014first}.
On the one hand, we prove the expected convergence rate for the $\tfrac32$-volatility model
on the condition $ \tfrac{\alpha}{\beta^2} \geq \tfrac52 $, also improving the restriction $ \tfrac{\alpha}{\beta^2} \geq 5 $ 
required in \cite{neuenkirch2014first}.
On the other hand, our approach is able to treat the Ait-Sahalia model in the general critical case $\kappa + 1 = 2 \rho$, 
with a first mean-square convergence order identified under conditions
$\tfrac{ \alpha_2 } { \sigma^2 } > 2 \kappa - \tfrac32 $ 
and $\tfrac{ \alpha_2 }{ \sigma^2 } > \tfrac { \kappa + 1 } { 2 \sqrt{2} }$, 
which is, as far as we know, missing in the literature.
For the special critical case $\kappa = 2, \rho = 1.5$ studied in \cite{neuenkirch2014first}, 
the restriction $ \tfrac{ \alpha_2 } { \sigma^2 } > 5 $ there
is moderately relaxed to $\tfrac{ \alpha_2 } { \sigma^2 } > \tfrac52 $ here.

%The backward Milstein scheme\cite{MS13,AK17,BIK16,BIK17}.
%The implicit schemes have advantages for stiff SDEs

%Before closing the introduction section, we would like to mention that,
%
%As opposed to the implicit methods,  explicit methods based on modifications of the usual It\^o-Taylor schemes,
%such as the balanced/tamed schemes \cite{Hut12con, HutJen14,HutJent2015, LiuMao13, SabaEul13, 
%Saba16, ZWH14, Tretyakov13, Zhang17}, \cite{WG13,Zhang17}, and the truncated methods
%\cite{Mao15,Mao16,Guo18truncated,li2020explicit},
%%\cite{HJ11, CHJ13, HutJen14, HutJent2015, LiuMao13, Mao15, Mao16, Guo18truncated, SabaEul13, Saba16, WG13,
%%ZWH14, Tretyakov13, Zhang17, BIK16, BIK17} (see also the references therein).
%are also able to be strongly convergent under certain non-globally Lipschitz conditions.
%%%%
%%%%

To conclude, the main contributions of the article are summarized as follows: (i) a family of double implicit Milstein-type schemes 
is introduced for multi-dimension SDE systems with non-Lipschitz coefficients;  
(ii) a novel approach of the error analysis is developed to recover the mean-square convergence rate of order one for the schemes, 
which fills several gaps in the literature;
(iii) the optimal mean-square convergence rate of the positivity preserving schemes applied to two financial models is obtained for the first time 
and more relaxed conditions are required, compared with existing relevant results for first order schemes in the literature.
%improved convergence results are obtained when the schemes are used to solve financial models. 
Therefore, this work can justify an efficient Multilevel Monte Carlo method \cite{giles08multilevel} for SDEs 
with non-globally Lipschitz coefficients including the above financial models.
% parameter range

The remainder of this article is structured as follows. In the forthcoming section, a setting is formulated and 
a family of new Milstein-type schemes are introduced. Upper mean-square error bounds of the proposed schemes
are then elaborated in section \ref{sect:upper-error-bound}. Equipped with the obtained error bounds, mean-square convergence 
rates of the schemes are analyzed in section \ref{sec:conv-rates-global-polynomial} for a general class of SDEs, 
under further globally polynomial growth conditions.
Additionally,  applications of the error bounds to two schemes for several SDE models in practice 
are examined in section \ref{sect:applications}, with an optimal convergence rate revealed. 
Further,  some numerical tests are provided to confirm the theoretical findings and a brief conclusion is made at the end of the article.
%some preliminaries are collected...

\section{SDEs and the proposed schemes}
%\section{Preliminaries}
\label{sect:setting}
%

%Let $ d,m, N \in \mathbb{N} $, $ T \in (0, \infty) $ and let $ \left( \Omega, \mathcal{ F }, \P,  \{ \mathcal{ F }_t \}_{ t \in [0,T] }
%\right) $ be a filtered probability space. 
%
Throughout this paper, we use $\N$ to denote the set of all positive integers
%and denote $\N_0 = \{ 0 \} \cup \N$.
and let $ d,m \in \N$, $ T \in (0, \infty) $ be given. Let $\| \cdot \|$ and $ \langle \cdot, \cdot \rangle $ denote the Euclidean
norm and the inner product of vectors in $\R^d$, respectively. 
Adopting the same notation as the vector norm, we denote $\|A\| : =\sqrt{\tr(A^{T}A)}$ as the trace norm of a matrix $A \in \R^{d \times m}$.
Given a filtered probability space $ \left( \Omega, \mathcal{ F }, \{ \mathcal{ F }_t \}_{ t \in [0,T] }, \P 
\right) $, we use $\E$ to mean the expectation and $L^{r} (\Omega; \R^d ), r \geq 1 $, to denote 
the family of $\R^d$-valued random variables $\xi$ satisfying $\E[ \|\xi \|^{r}]<\infty$.
Let us consider the following SDEs of It\^o type:
\begin{equation} \label{eq:SODE}
\begin{split}
\left\{
    \begin{array}{ll}
    \dd X_t  = f ( X_t ) \, \dd t + g ( X_t ) \,\dd W_t,
    \quad
    t \in (0, T],
    \\
   X_0 = x_0,
 \end{array}\right.
 \end{split}
\end{equation}
where $f \colon \mathbb{R}^d \rightarrow \mathbb{R}^d$ is the drift coefficient function, and
$g \colon \mathbb{R}^d \rightarrow \mathbb{R}^{d \times m}$ is the diffusion coefficient function, 
frequently written as $g = (g_{i,j})_{d \times m} =  (g_1, g_2,..., g_m)$ for $g_{i,j} \colon \mathbb{R}^d \rightarrow \mathbb{R}$
and $g_j \colon \mathbb{R}^d \rightarrow \mathbb{R}^{d }, i \in \{ 1, 2,..., d \}, j \in \{ 1, 2,..., m \} $.
Moreover, $W_{\cdot} \colon [0,T] \times \Omega \rightarrow \mathbb{R}^m $ stands for the $\mathbb{R}^m$-valued
standard Brownian motions with respect to $ \{ \mathcal{ F }_t \}_{ t \in [0,T] } $ and the initial data $X_0 \colon \Omega \rightarrow \mathbb{R}^d $ is assumed to be $\mathcal{ F }_0$-measurable.

In general, the  system of SDEs \eqref{eq:SODE} does not have a closed-form solution. 
In order to approximate \eqref{eq:SODE}, we construct a uniform mesh on $[0, T]$
with $ h = \tfrac{T}{N}$ being the stepsize, for any $ N \in \N$.
%
%%%%
%
On the uniform mesh, we propose a family of double implicit Milstein methods
with a pair of method parameters $(\theta,\eta)$, given by
\begin{equation}
\label{eq:general-scheme}
 \begin{split}
 Y_{ n + 1} = & Y_{ n } + \theta f ( Y_{n+1} )h
         + ( 1-\theta ) f( Y_n )h
            + g(Y_{n}) \Delta W_n
            + \sum^m_{j_1,j_2=1}
            \mathcal{L}^{j_1} g_{j_2}(Y_{n}) I_{j_1,j_2}^{t_n,t_{n+1}}
            \\
            & \quad
            +
             \tfrac {\eta} {2} \sum_{ j = 1}^{ m} \mathcal{L}^{j } g_{j }(Y_{n}) h
            -
            \tfrac {\eta} {2} \sum_{ j = 1}^{ m} \mathcal{L}^{j } g_{j }(Y_{n+1}) h,
            \qquad
            Y_{0}= X_0,
\end{split}
\end{equation}
where $\theta, \eta \in [0, 1]$,
$\Delta W_{ n } : = W_{t_{ n + 1 } } - W_{t_ {n} } $, $ n \in \{0, 1, 2, \ldots, N - 1\}$, and
\begin{equation}\label{eq:scheme-notation}
\mathcal{L}^{j_1} : = \sum^d_{k=1} g_{k,j_1}\frac{\partial}{\partial x^k},
\quad
I_{j_1,j_2}^{t_n,t_{n+1}} : = \int^{t_{n+1}}_{t_n} \int^{s_2}_{t_n} d W_{s_1}^{j_1}
d W_{s_2}^{j_2},
\quad
j_1, j_2 \in \{1, 2,..., m\}.
\end{equation}
%%%%
%Evidently, 
In the following we use $ \tfrac{ \partial \phi } { \partial x} $ to denote  the Jacobian matrix of the vector function
$\phi \colon \mathbb{R}^d \rightarrow \mathbb{R}^{d }$ and one can observe that, 
for $g_{j_2} \colon \mathbb{R}^d \rightarrow \mathbb{R}^{d }, j_2 \in \{ 1, 2,..., m \}$,
\begin{equation} \label{eq:derivative-operator}
\mathcal{L}^{j_1} g_{j_2} ( x ) 
=
 \sum^d_{k=1} g_{k,j_1}\frac{\partial  g_{j_2} ( x )}{\partial x^k}
 = 
\frac{ \partial g_{j_2} } { \partial x} ( x ) g_{j_1} ( x ), \quad x \in \R^d.
\end{equation}
%where  we denote by $ \tfrac{ \partial g_{j_2} } { \partial x} $ the Jacobian matrix of the vector function
%$g_{j_2} \colon \mathbb{R}^d \rightarrow \mathbb{R}^{d }, j_2 \in \{ 1, 2,..., m \}$.
%
By incorporating a pair of method parameters $\theta, \eta \in [0, 1]$ 
into the drift and diffusion coefficients, the newly proposed schemes are implicitly defined 
when $\theta + \eta \neq 0$ and their well-posedness will be discussed later.
Taking $\eta = 0$ in \eqref{eq:general-scheme}, the above double implicit Milstein methods \eqref{eq:general-scheme}
reduce to the classic $\theta$ Milstein methods \cite{kloeden1992numerical}, which are drift implicit and given by
\begin{equation}
\label{eq:semi-implicit-scheme}
 \begin{split}
 Y_{ n + 1} = & Y_{ n } + \theta f ( Y_{n+1} )h
         + ( 1-\theta ) f( Y_n )h
            + g(Y_{n}) \Delta W_n
            + \sum^m_{j_1,j_2=1}
            \mathcal{L}^{j_1} g_{j_2}(Y_{n}) I_{j_1,j_2}^{t_n,t_{n+1}},
            \quad
            Y_{0}= X_0.
\end{split}
\end{equation}

In general, a straightforward introduction of implicitness into approximations of the diffusion term containing random variables
suffers from unbounded numerical approximations with positive probability,  see \cite[Chapter 1.3.4]{milstein2013stochastic} for clarifications.
%
%In \cite{WGW12}, a family of fully implicit Milstein methods were constructed, 
%but only valid for SDEs with a special class of linear diffusion coefficients. 
%
When the diffusion coefficient $ g $ fulfills the so-called commutativity condition, namely,
\begin{equation} \label{eq:cc}
\mathcal{L}^{j_1} g_{ j_2} = \mathcal{L}^{j_2} g_{ j_1}, \quad j_1,j_2 \in \{ 1,...,m \},
\end{equation}
by recalling 
%$
\begin{equation}
I_{j_1,j_2}^{t_n,t_{n+1}} + I_{j_2, j_1}^{t_n,t_{n+1}} =  \Delta W_n^{j_1}\Delta W_n^{j_2}, 
\quad
j_1,j_2 \in \{ 1,...,m\}, j_1 \neq j_2
\end{equation}
%$
and 
%$
\begin{equation}
I_{j,j}^{t_n,t_{n+1}}  = \tfrac12 ( | \Delta W_n^j |^2 - h ), 
\quad
j \in \{ 1,...,m\},
\end{equation}
%$ 
%%one can recast the drift implicit schemes \eqref{eq:semi-implicit-scheme} as
%%\begin{equation}
%%\label{eq:commutative-drift-implicit-scheme}
%% \begin{split}
%% Y_{ n + 1} = & Y_{ n } + \theta f ( Y_{n+1} )h
%%         + ( 1-\theta ) f( Y_n )h
%%            + g(Y_{n}) \Delta W_n
%%            \\
%%            & \quad
%%            + \tfrac12
%%            \sum^m_{j_1,j_2=1}
%%            \mathcal{L}^{j_1} g_{j_2}(Y_{n}) \Delta W_n^{j_1} \Delta W_n^{ j_2 }
%%             -
%%             \tfrac { 1 } {2} \sum_{ j = 1}^{ m}
%%             \mathcal{L}^{j } g_{j }(Y_{n}) h
%%             ,
%%            \qquad
%%            Y_{0}= X_0.
%%\end{split}
%%\end{equation}
%take advantage of ...
one can recast the proposed double implicit Milstein method \eqref{eq:general-scheme} as
\begin{equation}
\label{eq:commutative-scheme}
 \begin{split}
 Y_{ n + 1} = & Y_{ n } + \theta f ( Y_{n+1} )h
         + ( 1-\theta ) f( Y_n )h
            + g(Y_{n}) \Delta W_n
            + \tfrac12
            \sum^m_{j_1,j_2=1}
            \mathcal{L}^{j_1} g_{j_2}(Y_{n}) \Delta W_n^{j_1} \Delta W_n^{ j_2 }
            \\
            & \quad
            -
             \tfrac { ( 1 - \eta ) } {2} \sum_{ j = 1}^{ m}
             \mathcal{L}^{j } g_{j }(Y_{n}) h
            -
            \tfrac {\eta} {2} 
            \sum_{j = 1}^{m} \mathcal{L}^{j } g_{j}(Y_{n+1}) h,
            \qquad
            Y_{0}= X_0.
\end{split}
\end{equation}
Here an implicit approximation is introduced with an additional method parameter $\eta \in [0, 1]$ only in the last term
that does not contain any random variable. In \cite{Buckwar2011comparative}, such schemes were applied to scalar linear
SDEs with several multiplicative noise terms ($m > 1$) and their mean-square stability properties were studied.
In particular, the commutativity condition \eqref{eq:cc} is fulfilled when $m= 1$
%
%@@MORE COMMENTS!!!@@
%the scalar SDEs ($ d = m = 1 $) 
and the newly proposed schemes 
\eqref{eq:general-scheme} (or \eqref{eq:commutative-scheme} equivalently) applied to the scalar SDEs ($ d = m = 1 $) reduce to
\begin{equation}
\label{eq:scheme-dm1}
 \begin{split}
 Y_{ n + 1} = & Y_{ n } + \theta f ( Y_{ n + 1} )h
         + ( 1-\theta ) f( Y_n )h
            + g(Y_{n}) \Delta W_n
            +{\tfrac{1}{2}} g'g( Y_{n} )\Delta W_n^2
            \\
             &- { \tfrac{(1-\eta)}{2} } g'g( Y_{n} ) h
             - { \tfrac{\eta}{2} } g'g( Y_{ n + 1} ) h ,
             \qquad
             Y_{0}= X_0.
  \end{split}
\end{equation}
Such schemes have been examined in \cite{higham2013convergence}, where the authors recovered the strong
convergence rate only under globally Lipschitz conditions. 
Moreover, the authors used \eqref{eq:scheme-dm1} to solve the $\tfrac32$-volatility model \eqref{eq:intro-32model} 
and proved its strong convergence  with no convergence rate revealed. 
Roughly speaking, the main difficulty of recovering the convergence rate
is caused by the super-linearly growing diffusion coefficients of SDEs.
In the literature, a lot of researchers \cite{andersson2017mean,beyn2016stochastic,beyn2017stochastic,gan2020tamed,guo2018truncated, 
Hutzenthaler15numerical,hutzenthaler2020perturbation,kumar2019milstein,li2020explicit,mao2016convergence,mao2013strongJCAM, 
mao2013strong,sabanis2016euler,Szpruch18V-integrability,Tretyakov2013fundamental,wang2018mean,zhang2017order-preserving} 
attempt to analyze strong approximations of SDEs with super-linearly growing diffusion coefficients.
%Actually,  
However, the strong convergence rate of the classical $\theta$-Milstein method 
in the regime of possibly super-linearly growing diffusion coefficients is, up to the best of our knowledge, still an open problem.
The present article aims to establish a mean-square convergence theory for the generalized $\theta$-Milstein schemes 
\eqref{eq:introduction-scheme} within a general framework, which fills several gaps in the literature 
and provides improved convergence results for computational finance.
Finally, it is worthwhile to emphasize that the newly proposed double implicit Milstein methods \eqref{eq:general-scheme} 
do not require the commutativity condition \eqref{eq:cc} and thus work for non-commutative noise driven SDEs.

%When $ \theta = \eta = 1 $ ,
% we have
% \begin{equation*}
% \begin{split}
% e_{k+1} - \Delta f^{ X, Y }_{k+1} h
% + { \tfrac{ 1 }{ 2 } } \Delta g'g ^{ X, Y }_{ k+1 } h
%               = &
%               e_k + \Delta g^{ X, Y }_k \Delta W_k
%               + { \tfrac {1}{2} } \Delta g'g^{ X, Y}_k \Delta W_k^2  - R_{k+1}
% \end{split}
%\end{equation*}
%and
% \begin{equation}
% \label{ eq : * }
% \begin{split}
%   \|  e_{k+1} - \Delta f^{ X, Y }_{k+1} h
%     + { \tfrac{ 1 }{ 2 } } \Delta g'g ^{ X, Y }_{ k+1 } h \| ^2
%      = &
%     \| e_k \| ^2 + \| \Delta g^{ X, Y }_k \Delta W_k  \| ^2
%     + { \tfrac {1}{4} } \| \Delta g'g^{ X, Y }_k \Delta W_k^2 \| ^2 + \| R_{ k+1 } \| ^2
%     \\
%     & + 2 \langle e_k , \Delta g^{ X, Y }_k \Delta W_k \rangle
%       + \langle e_k , \Delta g'g^{ X, Y }_k \Delta W_k^2 \rangle
%       - 2 \langle e_k , R_{ k+1 } \rangle
%      \\
%     & +  \langle \Delta g^{ X, Y }_k \Delta W_k , \Delta g'g^{ X, Y }_k \Delta W_k^2 \rangle
%       - 2 \langle \Delta g^{ X, Y }_k \Delta W_k , R_{ k+1 } \rangle
%        \\
%     & - \langle \Delta g'g^{ X, Y }_k \Delta W_k^2 , R_{ k+1 } \rangle .
% \end{split}
%\end{equation}
%%%%
%%
%
\section{Upper mean-square error bounds for the schemes}
\label{sect:upper-error-bound}
The aim of the present section is to derive upper mean-square error bounds of the implicit Milstein type methods
for SDEs taking values in  a domain $D \subset \R^d$,
which will help us to easily analyze the mean-square convergence rate of the schemes later.
%
%To guarantee the well-posedness of SDEs and schemes \eqref{eq:general-scheme} in the domain $D$, we make the following assumptions.
%%
%\begin{assumption}[Well-posedness of SDEs and STMs]
%\label{ass:well-possedness}
%%Under Assumption \ref{ass:monotonicity-condition},
%Assume SDE \eqref{eq:SODE} possesses a unique $ \{ \mathcal{ F }_t \}_{ t \in [0,T] } $-adapted
%$D$-valued solution with continuous sample paths,
%$X \colon [0,T] \times \Omega \rightarrow  D \subset \mathbb{R}^d $,
%satisfying $\sup_{ s \in [ 0, T]} \E [ \|  X_{ s } \|^2 ] < \infty$ and $\sup_{ s \in [ 0, T]} \E [ \| f ( X_{ s } ) \|^2 ] < \infty$.
%%obeying, for any $p \leq p^*$,
%%\begin{equation}
%%\E [ \|  X_t \|^p ] \leq ...
%%\end{equation}
%Moreover, suppose the proposed scheme \eqref{eq:general-scheme}, as a drift-implicit time-stepping scheme ($\theta >0$),
%admits a unique $ \{ \mathcal{ F }_{t_n} \}_{n = 0}^N $-adapted solution $\{ Y_n \}_{n = 0}^N$, taking values in the domain $D$.
%\end{assumption}
%
To this end, we set up a general framework  by making two key assumptions as follows.
\begin{assumption}[Generalized monotonicity conditions in a domain]
\label{ass:monotonicity-condition}
Assume that the diffusion coefficients $ g_j \colon D \rightarrow \mathbb{R}^{d }, j \in \{ 1, 2,..., m \}  $ 
are differentiable in a domain $D \subset \R^d$ and  that the drift coefficient $f \colon D \rightarrow \mathbb{R}^d$ 
and the diffusion coefficient $g =  (g_1, g_2,..., g_m) \colon D \rightarrow \mathbb{R}^{d \times m}$ of SDEs \eqref{eq:SODE}
satisfy certain monotonicity conditions in $D \subset \R^d$. More accurately,
for method parameters $\theta, \eta \in [0, 1]$
there exist constants $q \in (2, \infty) $, $ \varrho \in (1, \infty) $, $L_1, L_2 \in [0, \infty)$ and $ h_0 \in (0, T]$ such that,
$\forall x, y \in D$, $ h = \tfrac{T}{N} \in (0, h_0)$,
 \begin{align}
 \label{eq:mono-condition1}
% \begin{split}
 &
 2 \langle x - y ,f ( x ) - f ( y ) \rangle
  +  ( q - 1 ) \|g ( x ) - g ( y )\| ^2
  + \tfrac{ \varrho }{ 2 } h
  \sum_{j_1,j_2=1}^m \big\| 
   \mathcal{L}^{j_1}g_{j_2}( x ) - \mathcal{L}^{j_1}g_{j_2}( y ) \big\| ^2
 \\ \nonumber
 & \quad
  + \eta h  \Big \langle  \sum_{j=1}^m 
  \big[ \mathcal{L}^{j}g_{j}( x )
   - \mathcal{L}^{j}g_{j}( y ) \big] , f ( x ) - f ( y ) \Big \rangle
  +  ( 1 - 2 \theta ) h \| f ( x ) - f ( y ) \|^2
   \leq L _ 1 \| x - y \|^2 ,
%\:
%\forall x, y \in D,
% \end{split}
% \end{equation}
 %
 %%
% \begin{equation}\label{ineq:mono2}
% \begin{split}
% &
% 2\langle x ,f ( x ) \rangle
%  +  ( p ^ * - 1 ) \|g ( x ) \| ^2
%  + \tfrac{ h }{ 2 } \| g'g ( x )\| ^2
%  + \eta h \langle  g'g ( x ) , f ( x ) \rangle
% \\
% &
%  +  ( 1 - 2 \theta ) h \| f ( x ) \|^2
%   \leq L _ 2 ( 1 + \| x \|^2 ) ,
% \end{split}
% \end{equation}
%
%%
\\
% \begin{equation}
 \label{eq:mono-condition2}
% \begin{split}
 &
 \Big \langle x - y ,\theta [ f ( x ) - f ( y ) ]
 - \tfrac{\eta}{2} 
 \sum_{j=1}^m [ \mathcal{L}^{j}g_{j}( x )
  - \mathcal{L}^{j} g_{j}( y )]  \Big \rangle
   \leq L_2 \| x - y \|^2
%\:
%\forall x, y \in D
.
% \end{split}
 \end{align}
 %
 %%
% \begin{equation}\label{ineq:mono4}
% \begin{split}
% &
% \langle x ,\theta f ( x )- \tfrac{\eta}{2} g'g ( x ) \rangle
%   \leq L _ 4 ( 1 + \| x \|^2 ) ,
% \end{split}
% \end{equation}
 %
 %%
\end{assumption}
%
%%%%%%%%%%%%%%%%%The assumption seems unnecessary on 01.09.2019%%%%%%%%%%%
%
%\begin{assumption}
%\label{ass:fg-greater-c0-c1}
%%
%For method parameters $\theta, \eta \in [0, 1]$, there exist constants $ \rho_1 \in [0,1), \rho_2 \in [0,1) , c_0, c_1 \in \R$ such that, for any $ x, y \in D $,
% \begin{equation}\label{eq:mono-condition3}
% \begin{split}
% &
% \rho_1 \theta ^ 2 \| f(x) - f ( y ) \| ^ 2 +
% \tfrac{\rho_2 \eta ^2}{4} \|  \sum_{j=1}^m [ \mathcal{L}^{j}g_{j}( x )
%  - \mathcal{L}^{j}g_{j}( y )] \|^2
% \\
% & \quad - \eta \theta \langle f ( x ) - f ( y ) ,   \sum_{j=1}^m [ \mathcal{L}^{j}g_{j}( x )
%  - \mathcal{L}^{j}g_{j}( y )] \rangle
%  \geq c_0 + c_1 \| x - y \|^2.
% \end{split}
% \end{equation}
%%
%\end{assumption}
%%%%%
%%%%%%%%%%%%%%%%%%%%%%%%%%%%%%%%%%%%%%%%%%%%%%%%%%%%%%%%

%\begin{thm}[Well-posedness of SDEs and STMs]\label{thm: well-possedness}
%Under Assumption \ref{ass:monotonicity-condition}, SDE \eqref{eq:SODE} possesses a unique solution
%$X \colon [0,T] \times \Omega \rightarrow \mathbb{R}^d $, obeying, for any $p \leq p^*$,
%\begin{equation}
%\E [ \|  X_t \|^p ] \leq ...
%\end{equation}
%Furthermore, the STMs are well-defined, as a drift-implicit time-stepping scheme.
%\end{thm}
%{\it Proof of Theorem \ref{thm: well-possedness}.}
%%%%%%
%%%%%%
%
%
{\color{black}Conditions in Assumption \ref{ass:monotonicity-condition} are crucial to the  error analysis 
for  the proposed schemes and are called generalized monotonicity conditions in a domain $D$.
When $\eta = 0$, the implicit methods \eqref{eq:general-scheme} reduce to 
the classic $\theta$ Milstein methods \eqref{eq:semi-implicit-scheme} and the above two conditions are satisfied as
$\theta \in [\tfrac12, 1]$, $g$, $\mathcal{L}^{j_1}g_{j_2}, j_1, j_2 \in \{ 1, 2,..., m \} $ satisfy the globally Lipschitz condition
\begin{equation}
\|g ( x ) - g ( y )\| ^2
  + 
  \sum_{j_1,j_2=1}^m \big\| 
   \mathcal{L}^{j_1}g_{j_2}( x ) - \mathcal{L}^{j_1}g_{j_2}( y ) \big\| ^2
\leq L \| x - y \|^2,
\quad
\forall x, y \in D,
\end{equation} 
and $f$ obeys the monotonicity condition
\begin{equation}
\label{eq:f-monocity-condiion}
\langle x - y ,  f ( x ) - f ( y ) \rangle
   \leq L \| x - y \|^2,
\quad
\forall x, y \in D.
\end{equation}
Such a global monotonicity condition \eqref{eq:f-monocity-condiion} is frequently used in the literature, 
to ensure the well-posedness of drift-implicit  methods and to derive their strong convergence rates. 
When the diffusion $g$ is not  globally Lipschitz, which is the case for the aforementioned models \eqref{eq:intro-32model},
\eqref{eq:intro-Ait-Sahalia-model-SDE},  things become much more involved. 
As one can see later, Assumption \ref{ass:cor-monotonicity-condition-theta-Milstein}
and Assumption \ref{ass:cor-monotonicity-condition-backward-Milstein} below provide sufficient conditions that
imply Assumption \ref{ass:monotonicity-condition} and allow for non-globally Lipschitz diffusion coefficient.}
Since Assumption \ref{ass:monotonicity-condition} alone does not suffice to guarantee
the well-posedness of SDEs and the considered schemes in the domain $D$, we additionally require the following assumptions.
\begin{assumption}[Well-posedness of SDEs and schemes]
\label{ass:well-possedness}
%Under Assumption \ref{ass:monotonicity-condition},
Assume $\sum_{ j_1, j_2 = 1}^m  \E [ \| \mathcal{L}^{j_1} g_{j_2} ( X_{ 0 } ) \|]^2 < \infty$ 
and SDE \eqref{eq:SODE} possesses a unique $ \{ \mathcal{ F }_t \}_{ t \in [0,T] } $-adapted
$D$-valued global solution with continuous sample paths,
$X \colon [0,T] \times \Omega \rightarrow  D \subset \mathbb{R}^d $,
satisfying $\sup_{ s \in [ 0, T]} \E [ \|  X_{ s } \|^2 ] + \sup_{ s \in [ 0, T]} \E [ \| f ( X_{ s } ) \|^2 ] 
%+ \sup_{ s \in [ 0, T]}  \sum_{ j_1, j_2 = 1}^m  \E [ \| \mathcal{L}^{j_1} g_{j_2} ( X_{ s } ) \|]^2 
< \infty$.
%obeying, for any $p \leq p^*$,
%\begin{equation}
%\E [ \|  X_t \|^p ] \leq ...
%\end{equation}
Moreover, for $\theta, \eta \in [0, 1]$ specified in Assumption \ref{ass:monotonicity-condition} 
suppose the proposed scheme \eqref{eq:general-scheme} 
%as an implicit time-stepping scheme ($\theta >0$),
admits a unique $ \{ \mathcal{ F }_{t_n} \}_{n = 0}^N $-adapted solution $\{ Y_n \}_{n = 0}^N$, $N \in \N$, taking values in the domain $D$.
\end{assumption}

{\color{black}We mention that Assumption \ref{ass:well-possedness} is necessary but not strict.
For example, by taking $D$ to be the whole space $\R^d$, i.e., $D = \R^d$, Assumption \ref{ass:monotonicity-condition}
and Assumption \ref{ass:f-polynomial-growth} below 
together suffice to imply Assumption \ref{ass:well-possedness}.
In addition, some models in practice taking values in $D = (0, \infty)$ are also given in 
Section \ref{sect:applications} to satisfy the above assumptions.}
Under the above two assumptions, we are able to formulate the following main result of this section
that offers upper mean-square error bounds  for the underlying schemes.
%\subsection{Upper error bounds of STMs}
%\subsection{Upper bound of the global approximation error}
%%%%%%%%
\begin{thm}[Upper mean-square error bounds]
\label{thm:upper-error-bound}
Let Assumptions \ref{ass:monotonicity-condition}, \ref{ass:well-possedness} hold
with $\theta \in [\tfrac12, 1]$ and $2 L_2 h \leq \nu $ for some $ \nu \in (0, 1) $.
Let $\{ X_{ t} \}_{ t \in [0, T]}$ and $ \{Y_n\}_{0\leq n\leq N} $ be solutions to 
\eqref{eq:SODE} and \eqref{eq:general-scheme}, respectively.
Then there exists a uniform constant $C$ such that, for any $ n \in \{ 1,2,..., N \}$, $N \in \N$,
\begin{equation}
\begin{split}
\label{eq:Thm-Error-Bound}
\E \big [ \big \lVert X_{ t_n } - Y_n \big \rVert ^2 \big ]
\leq
C
%\big[ ( \tfrac{q -  1 }{q - 2 } + \tfrac{ 1 }{ \varepsilon_1 }
%    + \tfrac{ \eta ^2 }{ \theta ^2 } \tfrac{ 1 }{ 2 \varepsilon_2 } )
%       e ^ T \sum_{ i = 1 }^{ k } \E [ \| R_i \| ^2 ] \big]
%    + \tfrac{ \theta ^ 2 + 1 }{ \theta ^ 2 h } e ^ T
%    \sum_{ i = 1 }^{ k }
%    \E [ \| \E( R_{ i } \, | \mathcal{F}_{ i-1} ) \| ^2 ]
\Big(  
\sum_{ i = 1 }^{ n } \E [ \| R_i \| ^2 ]
    + \tfrac{ 1 }{ h }
    \sum_{ i = 1 }^{ n }
    \E [ \| \E( R_{ i } \, | \mathcal{F}_{ i-1} ) \| ^2 ]
 \Big)
,
\end{split}
\end{equation}
where 
%for $ i \in \{ 1, 2, ..., n \} $ 
we denote
\begin{equation}
 \begin{split}
 \label{eq:Error-Remainder-Defn}
R_{ i }  := 
            & 
            \theta \! \int_{ t_{i-1} } ^ { t_{ i } } f( X_{s} ) -  f(X_{ t_{i} } ) \, \dd s
             +
             (1 - \theta ) \! \int_{ t_{i-1} } ^ { t_{ i } } f( X_{s} ) -  f(X_{ t_{i-1} }) \, \dd s
             + \int_{ t_{i-1} } ^ { t_{i} } g( X_{s} ) - g(X_{ t_{i-1} }) \, \dd W_s
             \\
             &
             - \! \sum^m_{j_1,j_2=1}
            \mathcal{L}^{j_1} g_{j_2}( X_{ t_{i-1} } ) I_{j_1,j_2}^{t_n,t_{n+1}}
            +
             \tfrac {\eta} {2} h \sum_{ j = 1}^{ m} \big[ \mathcal{L}^{j } g_{j }( X_{ t_{i} } )  - \mathcal{L}^{j } g_{j }( X_{ t_{i-1} } ) \big]
%             
%             \\
%           = & \theta \int_{ t_{i-1} }^{ t_{i} }
%            [ f( X_s )- f (X_{t_{i} }) ] \,\dd s
%            + (1-\theta) \int_{ t_{i-1} }^{ t_{i} }
%             [ f( X_s )- f (X_{t_{i-1 }}) ] \,\dd s
%             \\
%            & + \int_{ t_{i-1} }^{ t_{i} } [g( X_s )- g ( X_{t_{i-1} } ) ] \,\dd W_s
%            - { \tfrac{1}{2} } \sum_{j_1,j_2=1}^m
%              \mathcal{L}^{j_1}g_{j_2}( X_{ t_{i-1} } )\Delta W_{i-1}^{j_1}
%              \Delta W_{i-1}^{j_1} \\
%            &
%             + \sum_{j_1,j_2=1}^m \big[ \tfrac{(1-\eta)}{2}
%              \mathcal{L}^{j_1}g_{j_2}( X_{ t_{i-1} } ) h
%             +  { \tfrac{\eta}{2} }
%              \mathcal{L}^{j_1}g_{j_2}(X_{ t_{i} }) h \big] 
              ,
              \:\,  i \in \{ 1, 2, ..., n \}.
  \end{split}
\end{equation}
\end{thm}
%%%%%%%%%%%%
Throughout this paper, by $C$ we denote a generic deterministic positive constant, 
which might vary for each appearance but is independent of the time stepsize $h  =  \tfrac{T}{N} >0$, $N \in \N$.
%Throughout this paper, the notation C (with or without a subscript) will be used to denote a
%generic deterministic positive constant which might not be the same in each appearance.
%
It is interesting to observe that the term $R_i,  i \in \{ 1, 2, ..., N \}, N \in \N$ defined by \eqref{eq:Error-Remainder-Defn} 
only gets involved with the exact solutions to SDEs. Such error bounds can be used to analyze mean-square 
convergence rates of the schemes without relying on a priori high-order moment estimates of numerical approximations.
The proof of Theorem \ref{thm:upper-error-bound} is postponed, which requires the following two lemmas.
%For the purpose of rigorously validating Theorem \ref{thm:upper-error-bound}, we start from the following lemmas.
%
\begin{lem} \label{lem:R-k-bounded}
%
%Assume $\theta \in [0, 1]$ and $0 < h < \frac { 1 } { 2 L_4 }$ for $\theta \neq 0$ and $ 0 < h < \infty $ for $\theta = 0$.
%Further we assume $\E [ \| X_0 \|^2 ] < ...$
%It holds for all $n \in \{ 1,2,..., N \}$ that
%
Let Assumptions \ref{ass:monotonicity-condition}, \ref{ass:well-possedness} hold
and let $R_i,  i \in \{ 1, 2, ..., N \}$,  $ N \in \N $  be defined by \eqref{eq:Error-Remainder-Defn}.
%In addition, we assume $\sup_{ s \in [ 0, T]} \E [ \| f ( X_{ s } ) \|^2 ] < \infty$.
Then it holds 
%for all $k \in \{ 1,2,..., N \}$, $ N \in \N $ that
\begin{equation} \label{eq:Rk-bound}
\E [ \| R_i \|^2 ] < \infty,
\quad
\forall \, i \in \{ 1,2,..., N \}, N \in \N.
\end{equation}
\end{lem}
{\it Proof of Lemma \ref{lem:R-k-bounded}.}
%First we verify $ \E [ \| R_k \|^2 ] < \infty $ for all $ k \in \{  0,1,..., N\}$.
In light of \eqref{eq:mono-condition1}, one can show, $\forall x, y \in D$,
\begin{equation} \label{eq:g-g'g-difference0}
\begin{split}
& ( q - 1  ) \|g ( x ) - g ( y )\| ^2 
+ 
\tfrac{ \varrho }{ 2 } h
  \sum_{j_1,j_2=1}^m \big\| 
   \mathcal{L}^{j_1}g_{j_2}( x ) - \mathcal{L}^{j_1}g_{j_2}( y ) \big\| ^2
\\
&  \quad
\leq L _ 1 \| x - y \|^2 - 2 \langle x - y ,f ( x ) - f ( y ) \rangle
   -
   ( 1 - 2 \theta ) h  \| f ( x ) - f ( y ) \|^2
\\
& \qquad   
%   - \varrho \tfrac{ h }{ 2 } \sum_{j_1,j_2=1}^m \big\| 
%   \mathcal{L}^{j_1}g_{j_2}( x ) - \mathcal{L}^{j_1}g_{j_2}( y ) \big\| ^2
  - \eta h  \Big \langle  \sum_{j=1}^m 
  \big[ \mathcal{L}^{j}g_{j}( x )
   - \mathcal{L}^{j}g_{j}( y ) \big] , f ( x ) - f ( y ) \Big \rangle
\\
&   \quad
   \leq
   ( L _ 1 + 1 )  \| x - y \|^2
  + 
  \tfrac{ \varrho h } { 4 }   \sum_{j=1}^m 
  \big\| \mathcal{L}^{j}g_{j}( x )
   - \mathcal{L}^{j}g_{j}( y ) \big\|^2 
   + 
   \tfrac{ \varrho + (m + \varrho) h } { \varrho }  
   \| f ( x ) - f ( y ) \|^2,
\end{split}
\end{equation}
and thus, $\forall x, y \in D$,
\begin{equation} \label{eq:g-g'g-difference}
( q - 1  ) \|g ( x ) - g ( y )\| ^2 
+ 
\tfrac{ \varrho }{ 4 } h \!
  \sum_{j_1,j_2=1}^m \big\| 
   \mathcal{L}^{j_1}g_{j_2}( x ) - \mathcal{L}^{j_1}g_{j_2}( y ) \big\| ^2
   \leq
   ( L _ 1 + 1 )  \| x - y \|^2
   +
   \tfrac{ \varrho + (m + \varrho) h } { \varrho }    \| f ( x ) - f ( y ) \|^2.
\end{equation}
Combining this with Assumption \ref{ass:well-possedness} guarantees, for any $i \in \{ 1, 2,..., N\}, N \in \N$ and $s \in [t_{i-1}, t_i]$,
\begin{align}\label{eq:g-g-bound}
%\begin{split}
\E [ \| g ( X_s ) - g ( X_{t_{i-1}} )  \|^2] & < \infty,
\\
 \nonumber
h
\sum_{j_1,j_2=1}^m
       \E
       [
       \| \mathcal{L}^{j_1}  g_{j_2} ( X_{ t_{i-1} } ) \|^2
       ]
 & \leq
 2 h
 \sum_{j_1,j_2=1}^m
       \E
       [
       \| \mathcal{L}^{j_1}  g_{j_2} ( X_{ t_{i-1} } ) -  \mathcal{L}^{j_1}  g_{j_2} ( X_{ 0 } ) \|^2
       ]
 \\ 
 &
 \quad
 +
 2 h
  \sum_{j_1,j_2=1}^m
       \E
       [
       \| \mathcal{L}^{j_1}  g_{j_2} ( X_{ 0 } ) \|^2
       ]
       <
       \infty.
       \label{eq:g'g-bound}
%\end{split}
\end{align}
%The assumption $ \sum_{j_1, j_2 = 1}^{m} \sup_{ s \in [ 0, T]} \E [ \| \mathcal{L}^{j_1} g_{j_2} ( X_{ 0 } ) \|^2 < \infty $ 
This in turn implies, for any $i \in \{ 1, 2,..., N\}, N \in \N$ and $s \in [t_{i-1}, t_i]$,
\begin{equation}\label{eq:g-g'g-bound}
\begin{split}
\E \Big[ \Big\|
\int_{ t_{i-1} } ^ { t_{i} } g( X_{s} ) - g(X_{ t_{i-1} }) \, \dd W_s
\Big\|^2 \Big]
        =
\int_{ t_{i-1} } ^ { t_{i} }
\E [ \| g ( X_s ) - g ( X_{t_{i-1}} )  \|^2] \, \dd s
& < 
\infty,
\\
\E \Big[ \Big\| \sum_{j_1,j_2=1}^m  \mathcal{L}^{j_1}  g_{j_2} ( X_{ t_{i-1} } ) I^{t_k, t_{k+1}}_ {  j_1 , j_2  }  
       \Big\|^2 \Big]
        =
       \tfrac{h^2}{2} 
       \sum_{j_1,j_2=1}^m
       \E
       [
       \| \mathcal{L}^{j_1}  g_{j_2} ( X_{ t_{i-1} } ) \|^2
       ]
       & <
       \infty,
       \\
      \E
      \Big[ \Big\|
      \tfrac {\eta} {2} \sum_{ j = 1}^{ m} \mathcal{L}^{j } g_{j }( X_{ t_{i} } ) h
                   -
             \tfrac {\eta} {2} \sum_{ j = 1}^{ m} \mathcal{L}^{j } g_{j }( X_{ t_{i-1} } ) h
       \Big\|^2 \Big]
         & < \infty.
\end{split}       
\end{equation}  
The desired assertion follows, by taking \eqref{eq:g-g'g-bound} and 
the assumption $\sup_{ s \in [ 0, T]} \E [ \| f ( X_{ s } ) \|^2 ]  < \infty$ into account.
\qed
%
%To conclude , we derive the following lemma:
%

Based on the boundedness of $\E [ \| R_i \|^2 ]$, $i \in \{ 1, 2,..., N\}$, 
one can arrive at the subsequent moment bounds.
\begin{lem} \label{lem:MS-bounded}
%
%Assume $\theta \in [0, 1]$ and $0 < h < \frac { 1 } { 2 L_4 }$ for $\theta \neq 0$ and $ 0 < h < \infty $ for $\theta = 0$.
%Further we assume $\E [ \| X_0 \|^2 ] < ...$
%It holds for all $n \in \{ 1,2,..., N \}$ that
%
Let $\theta \in (0, 1], \eta \in [0, 1]$, $ 2  L_2  h \leq \nu $ for some $\nu \in (0, 1)$ and let
Assumptions \ref{ass:monotonicity-condition}, \ref{ass:well-possedness} hold.
%In addition, we assume $\sup_{ s \in [ 0, T]} \E [ \| f ( X_{ s } ) \|^2 ] < \infty$.
Then it holds for all $k \in \{ 0,1, 2,..., N \}$, $ N \in \N $ that
%%%
%
%\begin{equation}\label{eq:lem-ms-bounded}
%\begin{split}
%& \E [ \| Y_k - \theta f^{ Y }_k h
%      + { \tfrac{ \eta }{2} } g'g^{ Y }_k h \| ^ 2 ] < \infty , \quad
% \E [ \| Y_k \| ^ 2 ] < \infty ,
%\\
% & \E [ \| f^{ Y }_k \| ^ 2 ] < \infty , \quad
% \E [ \| g^{ Y }_k \| ^ 2 ] < \infty , \quad
% \E [ \| g'g^{ Y }_k \| ^ 2 ] < \infty,
%\end{split}
%\end{equation}
%
%%%
\begin{equation}\label{eq:lemma-numerical-bound}
\begin{split}
&  
%@@LESS!!!
%\E [ \| R_k \|^2 ] < \infty, \quad
\E \Big[ \Big\| e_k - \theta \Delta f^{ X, Y }_k h
      + \tfrac {\eta} {2} \sum_{ j = 1}^{ m} \Delta ( \mathcal{L}^{j} g_{j})_{k+1}^{X, Y} h \Big\| ^ 2 \Big] < \infty , \quad
 \E [ \| e_k \| ^ 2 ] < \infty ,
\\
 & \E [ \| \Delta f^{ X, Y }_k \| ^ 2 ] < \infty , \quad
 \E [ \| \Delta g^{ X, Y }_k \| ^ 2 ] < \infty , \quad
%\\
% &
%{\color{red}{\eta}} 
\sum_{j_1,j_2=1}^m 
   \E \big[
    \big\| 
   \Delta \big(\mathcal{L}^{j_1}g_{j_2} \big)_k^{ X, Y} \big\| ^2
   \big] < \infty,
%\:
%\E [ \| \tfrac {\eta} {2} \sum_{ j = 1}^{ m} \Delta ( \mathcal{L}^{j} g_{j})_{k+1}^{X, Y} h \| ^ 2 ] < \infty
%,
\end{split}
\end{equation}
where for any $k \in \{ 0,1,..., N \}, j_1, j_2 \in \{ 1, 2,..., m \}$ we denote  
%$ 
\begin{equation} \label{eq:short-notation}
\begin{split}
e_k : = X_{t_k} - Y_k ,
\quad
\Delta f^{ X, Y }_k  &: = f( X_{t_k} ) - f( Y_k ),
\quad
\Delta g^{ X, Y }_k : = g( X_{t_k} ) - g( Y_k ),
\\
\Delta (\mathcal{L}^{j_1} g_{j_2})^{ X, Y }_k & : =
     \mathcal{L}^{j_1}g_{j_2}( X_{t_k} ) - \mathcal{L}^{j_1}g_{j_2}( Y_k ).
\end{split}
\end{equation}
%$
\end{lem}
{\it Proof of Lemma \ref{lem:MS-bounded}.}
%{\color{red}{Estimate  $ \E [ \| e_k - \theta \Delta f^{ X, Y }_k h
%      + { \tfrac{ \eta }{2} }  \tilde{L} \Delta g^{ X, Y }_k
%       h \|^2 ] < \infty $... directly!!!} }
%
%
We first note that, for any $k \in \{ 0,1,..., N-1 \}, j_1, j_2 \in \{ 1, 2,..., m \}$,
\begin{equation}
\label{eq:solution}
 \begin{split}
X_{ t_{ k+1 } } = & 
X_{ t_{ k } } 
+
\int_{ t_k }^{ t_{k+1} } f ( X_s ) \, \dd s 
+
\int_{ t_k }^{ t_{k+1} } g ( X_s ) \, \dd W_s 
\\
= &
X_{ t _{ k } } + \theta f ( X _{ t _{k+1} } )h
         + ( 1-\theta ) f( X _{ t _ {k} } )h
            + g( X _ { t _{k} } ) \Delta W_k
               + \sum^m_{j_1,j_2=1}
            \mathcal{L}^{j_1} g_{j_2}( X _ { t _{k} } ) I_{j_1,j_2}^{t_n,t_{n+1}}
            \\
            & \quad
            +
             \tfrac {\eta} {2} \sum_{ j = 1}^{ m} \mathcal{L}^{j } g_{j }(X _ { t _{k} }) h
            -
            \tfrac {\eta} {2} \sum_{ j = 1}^{ m} \mathcal{L}^{j } g_{j }(X _ { t _{k+1} }) h  
              + R_{ k+1 } ,
  \end{split}
\end{equation}
where $ R_{k+1} $ is defined by \eqref{eq:Error-Remainder-Defn}.
Using the short-hand notation \eqref{eq:short-notation},  
we subtract  \eqref{eq:general-scheme}  from \eqref{eq:solution} to get
\begin{equation}
\label{eq:error-identity}
 \begin{split}
   e_{k+1} = & e_k + \theta \Delta f^{ X, Y}_{k+1} h
             + (1-\theta) \Delta f^{ X, Y }_k h + \Delta g^{ X, Y }_k \Delta W_k
             + \sum^m_{j_1,j_2=1}
             \Delta ( \mathcal{L}^{j_1} g_{j_2})_k^{X, Y}  I_{j_1,j_2}^{t_n,t_{n+1}}
             \\
             & 
             + \tfrac {\eta} {2} \sum_{ j = 1}^{ m} \Delta ( \mathcal{L}^{j} g_{j})_k^{X, Y}  h
            -
            \tfrac {\eta} {2} \sum_{ j = 1}^{ m} \Delta ( \mathcal{L}^{j} g_{j})_{k+1}^{X, Y} h 
             + R_{ k+1 } ,
             \quad
             k \in \{ 0,1,..., N-1 \}.
  \end{split}
\end{equation}
Denoting further
\begin{equation} \label{eq:J-short}
\begin{split}
 &
 \mathcal{J}^{ X, Y }_ k : =
 e_k - \theta \Delta f^{ X, Y }_k h
 +
      \tfrac {\eta} {2} \sum_{ j = 1}^{ m} \Delta ( \mathcal{L}^{j} g_{j})_{k}^{X, Y} h,
      \quad
      k \in \{ 0,1,..., N \},
% \\
% &
%\quad
% I^k_ {  j_1 , j_2  } : = { \tfrac{ 1 }{ 2 } }
% ( \Delta W_k^{j_1}\Delta W_k^{j_2} - h ) ,
 \end{split}
\end{equation}
one can recast \eqref{eq:error-identity} as
\begin{equation}
\mathcal{J}^{ X, Y }_ {k+1} 
=  
\mathcal{J}^{ X, Y }_ k  
+
 \Delta f^{ X, Y }_k h + \Delta g^{ X, Y }_k \Delta W_k
             + 
             \sum^m_{j_1,j_2=1}
             \Delta ( \mathcal{L}^{j_1} g_{j_2})_k^{X, Y}  I_{j_1,j_2}^{t_k,t_{k+1}}
             +
             R_{k+1}.
\end{equation}
%
%As a consequence,
Squaring both sides of the above equality yields
\begin{equation}
 \label{eq:error-identity-0}
 \begin{split}
 \| \mathcal {J}^ { X, Y }_ { k+1 }  \| ^2
  &=
       \Big\|
       \mathcal{J}^{ X, Y }_ k +\Delta f^{ X, Y }_k h + \Delta g^{ X, Y }_k \Delta W_k
       +
          \sum^m_{j_1,j_2=1}
             \Delta ( \mathcal{L}^{j_1} g_{j_2})_k^{X, Y}  I_{j_1,j_2}^{t_k,t_{k+1}}
             +
             R_{k+1} 
             \Big \| ^2
 \\
  &=
      \| \mathcal {J}^ { X, Y }_ { k }  \| ^2 + h ^2 \| \Delta f^{ X, Y }_k \| ^2
    + \| \Delta g^{ X, Y }_k \Delta W_k  \| ^2
    + \Big \| \sum^m_{j_1,j_2=1}
             \Delta ( \mathcal{L}^{j_1} g_{j_2})_k^{X, Y}  I_{j_1,j_2}^{t_k,t_{k+1}}  \Big \| ^2
    + \| R_{ k+1 }  \| ^2
 \\
  & \quad +
     2 h \langle \mathcal{J}^{ X, Y }_ k , \Delta f^{ X, Y }_k \rangle
     + 2 \langle \mathcal{J}^{ X, Y }_ k , \Delta g^{ X, Y }_k \Delta W_k \rangle
     + 2 \Big \langle \mathcal{J}^{ X, Y }_ k ,
              \sum^m_{j_1,j_2=1}
             \Delta ( \mathcal{L}^{j_1} g_{j_2})_k^{X, Y}  I_{j_1,j_2}^{t_k,t_{k+1}} \Big \rangle 
              \\
     & \quad
     + 2 \langle \mathcal{J}^{ X, Y }_ k , R_{ k+1 } \rangle
     +
     2 h \langle \Delta f^{ X, Y }_k , \Delta g^{ X, Y }_k \Delta W_k \rangle 
     + 2 h \Big \langle \Delta f^{ X, Y }_k ,
              \sum^m_{j_1,j_2=1}
             \Delta ( \mathcal{L}^{j_1} g_{j_2})_k^{X, Y}  I_{j_1,j_2}^{t_k,t_{k+1}} \Big \rangle 
     \\
     & \quad
    + 2 h \langle \Delta f^{ X, Y }_k , R_{ k+1 }  \rangle 
    +
    2  \Big \langle \Delta g^{ X, Y }_k \Delta W_k ,
             \sum^m_{j_1,j_2=1}
             \Delta ( \mathcal{L}^{j_1} g_{j_2})_k^{X, Y}  I_{j_1,j_2}^{t_k,t_{k+1}} \Big \rangle
    \\
  & \quad 
  + 2 \langle \Delta g^{ X, Y }_k \Delta W_k , R_{ k+1 }  \rangle
  +
     2  \Big \langle
              \sum^m_{j_1,j_2=1}
             \Delta ( \mathcal{L}^{j_1} g_{j_2})_k^{X, Y}  I_{j_1,j_2}^{t_k,t_{k+1}} ,
    R_{ k+1 } \Big \rangle.
\end{split}
\end{equation}
With this at hand, we first prove $\E [ \| \mathcal {J}^ { X, Y }_ { k } \|^2] < \infty$ for all $k \in \{ 0, 1,..., N\}$ 
based on an induction argument.
Noting that $ X_0 = Y_0 $ we thus have $\E [ \| \mathcal {J}^ { X, Y }_ { 0 } \|^2] = 0$. We assume 
$\E [ \| \mathcal {J}^ { X, Y }_ { k } \|^2] < \infty$ for some $k \in \{ 0, 1,..., N-1\}$, 
which together with \eqref{eq:mono-condition2} implies
\begin{equation}\label{ineq:bound}
\begin{split}
  \infty  > \E [ \| \mathcal {J}^ { X, Y }_ { k } \|^2] & = \E \Big[ \Big\| e_k - \theta \Delta f^{ X, Y }_k h
 +
      \tfrac {\eta} {2} \sum_{ j = 1}^{ m} \Delta ( \mathcal{L}^{j} g_{j})_{k}^{X, Y} h \Big\|^2 \Big] \\
      & = 
       \E [ \| e_k \|^2 ] 
%       +  
%       h^2
%       \E [ \|  \theta \Delta f^{ X, Y }_k \|^2 ]
       + h^2
       \E \Big[ \Big\| \tfrac{\eta}{2}  \sum_{ j = 1}^{ m} \Delta ( \mathcal{L}^{j} g_{j})_{k}^{X, Y}
       - \theta \Delta f^{ X, Y }_k \Big \|^2 \Big]  \\
      & \quad
      - 
      2 \theta h \E [ \langle e_k, \Delta f^{ X, Y }_k \rangle ]
      + 
      \eta h \E \Big[ \Big \langle e_k , 
                  \sum_{ j = 1}^{ m} \Delta ( \mathcal{L}^{j} g_{j})_{k}^{X, Y} \Big \rangle \Big]
%      - 
%      \theta \eta h^2 \E \Big[ \Big \langle \Delta f^{ X, Y }_k   
%      ,\sum_{ j = 1}^{ m} \Delta ( \mathcal{L}^{j} g_{j})_{k}^{X, Y} \Big \rangle \Big] 
            \\
  & \geq 
     (1 - 2 h L_2 ) \E [ \| e_k \|^2 ] 
      + 
      h^2 
      \E \Big[ \Big\| \tfrac{\eta}{2}  \sum_{ j = 1}^{ m} \Delta ( \mathcal{L}^{j} g_{j})_{k}^{X, Y}
       - \theta \Delta f^{ X, Y }_k \Big \|^2 \Big].
%      \\
%  & \geq 
%      \E [ \| e_k \|^2 ] + \theta^2 h^2
%       \E [ \| \Delta f^{ X, Y }_k \|^2 ]
%       + \tfrac{\eta^2 h^2}{4}
%       \E [ \| \sum_{ j = 1}^{ m} \Delta ( \mathcal{L}^{j} g_{j})_{k}^{X, Y} \|^2 ]  \\
%      & \quad
%      - 2 h L_2 \E [ \| e_k \|^2 ] 
%      - \theta \eta h^2 \E [\langle \Delta f^{ X, Y }_k
%      ,\sum_{ j = 1}^{ m} \Delta ( \mathcal{L}^{j} g_{j})_{k}^{X, Y} \rangle] \\
%   & \geq
%     (1 - 2 h L_2 + c_1 h^2)\E [ \| e_k \|^2 ] + (1 - \rho_1)\theta^2 h^2
%       \E [ \| \Delta f^{ X, Y }_k \|^2 ]
%    \\
%    & \quad 
%    + (1 - \rho_2) \tfrac{\eta^2 h^2}{4}
%       \E \Big [ \Big \| \sum_{ j = 1}^{ m} \Delta ( \mathcal{L}^{j} g_{j})_{k}^{X, Y} \Big \|^2 \Big ]  
%    + c_0 h^2,
\end{split}
\end{equation}
%Here conditions \eqref{eq:mono-condition1}, \eqref{eq:mono-condition2} were used. 
Therefore, for 
$ 2 h L_2 \leq \nu < 1 $, $ \theta \in (0, 1] $ and for some $k \in \{ 0, 1, 2, ..., N-1\}$ it holds
\begin{equation}\label{ineq:bound2}
  \E [ \| e_k \|^2 ] < \infty , \quad 
  \E \Big[ \Big\| \tfrac{\eta}{2}  \sum_{ j = 1}^{ m} \Delta ( \mathcal{L}^{j} g_{j})_{k}^{X, Y}
       - \theta \Delta f^{ X, Y }_k \Big \|^2 \Big] < \infty.
%  \E [ \| \Delta f^{ X, Y }_k \|^2 ] < \infty,  \quad
%   \E [ \| \sum_{ j = 1}^{ m} \Delta ( \mathcal{L}^{j} g_{j})_{k}^{X, Y} \|^2 ]   < \infty .
\end{equation}
This along with the generalized monotonicity condition \eqref{eq:mono-condition1} shows,
 for some $k \in \{ 0, 1, 2, ..., N-1\}$,
\begin{equation}\label{eq:g-bound}
  \begin{split}
&
    (q-1) \E [\| \Delta g_k^{X, Y} \|^2 ]
    +
     \tfrac{ \varrho }{ 2 } h \sum_{j_1,j_2=1}^m 
   \E \big[
    \big\| 
   \Delta \big(\mathcal{L}^{j_1}g_{j_2} \big)_k^{ X, Y} \big\| ^2
   \big]
% \\     
%      & \quad 
%      \leq -2 \E [ \langle e_k, \Delta f^{ X, Y }_k \rangle ]
%      - \eta h \E [\langle \sum_{ j = 1}^{ m} \Delta ( \mathcal{L}^{j} g_{j})_{k}^{X, Y} ,
%       \Delta f^{ X, Y }_k \rangle] \\
%      & \qquad
%      - (1-2\theta)h \E [ \| \Delta f^{ X, Y }_k \|^2 ] 
%      + L_1 \E [ \| e_k \|^2 ]  
 \\     
      & \quad 
      \leq 
      L_1 \E [ \| e_k \|^2 ] 
      -
      2 \E [ \langle \mathcal {J}^ { X, Y }_ { k }, \Delta f^{ X, Y }_k \rangle ]
      - 
      h \E [ \| \Delta f^{ X, Y }_k \|^2 ]
 \\     
      & \quad 
      \leq 
L_1 \E [ \| e_k \|^2 ]  
+ 
\tfrac{1}{h} 
\E [ \| \mathcal {J}^ { X, Y }_ { k } \|^2 ]           
       < \infty.
  \end{split}
\end{equation}
In view of \eqref{eq:J-short}, \eqref{ineq:bound2}, \eqref{eq:g-bound} and the assumption $\theta >0$, one can easily see
\begin{equation} \label{eq:g-f-g'g-bound}
\E [\| \Delta g_k^{X, Y} \|^2 ] < \infty, 
\quad
\sum_{j_1,j_2=1}^m 
   \E \big[
    \big\| 
   \Delta \big(\mathcal{L}^{j_1}g_{j_2} \big)_k^{ X, Y} \big\| ^2
   \big] < \infty,
\quad
\E [ \| \Delta f^{ X, Y }_k \|^2 ]  < \infty
\end{equation}
 for some $k \in \{ 0, 1, 2, ..., N-1\}$.
%
%Thanks to \eqref{eq:mono-condition1}, \eqref{ineq:bound2}, \eqref{eq:g-bound}
%and the definition of $ R_i $, we can derive that
%$ \E [ \| R_{k+1} \|^2 ] < \infty $.
%
These bounded moments suffice to ensure,  for some $k \in \{ 0, 1, 2, ..., N-1\}$,
\begin{equation}\label{eq:stoc-moment-bounds}
  \begin{split}
    &  
       \E \Big[ \Big\| \sum_{j_1,j_2=1}^m  \Delta ( \mathcal{L}^{j_1}  g_{j_2} )^{ X, Y}_{k} I^{t_k, t_{k+1}}_ {  j_1 , j_2  }  
       \Big\|^2 \Big]
       =
       \tfrac{h^2}{2} 
       \sum_{j_1,j_2=1}^m
       \E
       [
       \| \Delta ( \mathcal{L}^{j_1}  g_{j_2} )^{ X, Y}_{k} \|^2
       ]
       <
       \infty,
       \\ &
       \E [ \| \Delta g_k^{X, Y} \Delta W_k \|^2 ]
       = h \E [ \| \Delta g_k^{X, Y} \|^2 ] < \infty, 
 \end{split}
 \end{equation}
 and
 \begin{equation} \label{eq:stoc-zero}
 \begin{split}
& \E \Big [ \Big \langle \mathcal{J}^{ X, Y }_ k ,
    \sum_{j_1,j_2=1}^m  \Delta ( \mathcal{L}^{j_1}  g_{j_2} )^{ X, Y}_{k} I^{t_k, t_{k+1}}_ {  j_1 , j_2  } \Big \rangle \Big]  = 0,
              \\
    &
      \E \Big[ \Big \langle \Delta f^{ X, Y }_k ,
             \sum_{j_1,j_2=1}^m \Delta ( \mathcal{L}^{j_1}  g_{j_2} )^{ X, Y}_{k} I^{t_k, t_{k+1}}_ {  j_1 , j_2  } \Big \rangle \Big]
      = 0,
%      \\
%      &
\quad
      \E[\langle \Delta f^{ X, Y }_k , \Delta g^{ X, Y }_k \Delta W_k \rangle]  = 0 \\
    & \E \Big[ \Big \langle \Delta g^{ X, Y }_k \Delta W_k ,
             \sum_{j_1,j_2=1}^m \Delta ( \mathcal{L}^{j_1}  g_{j_2} )^{ X, Y}_{k} I^{t_k, t_{k+1}}_ {  j_1 , j_2  } \Big \rangle \Big]
      = 0,
%      \\
%      &
\quad
      \E \big[ \langle \mathcal{J}^{ X, Y }_ k , \Delta g^{ X, Y }_k \Delta W_k \rangle \big]
      = 0.
  \end{split}
\end{equation}
%
%\begin{equation}
%\begin{split}
%   \| \mathcal {J}^ { X, Y }_ { k+1 }  \| ^2
%    & \leq
%      2 \| \mathcal {J}^ { X, Y }_ { k }  \| ^2 + 2 h ^2 \| \Delta f^{ X, Y }_k \| ^2
%    + 2 \| \Delta g^{ X, Y }_k \Delta W_k  \| ^2
%    + 2 \| \sum_{j_1,j_2=1}^m
%              \mathcal{L}^{j_1} \Delta g^{ X, Y}_{j_2,k} I^k_ {  j_1 , j_2  }  \| ^2
%    + 5 \| R_{ k+1 }  \| ^2   \\
%  & \quad +
%     2 h \langle \mathcal{J}^{ X, Y }_ k , \Delta f^{ X, Y }_k \rangle
%     + 2 \langle \mathcal{J}^{ X, Y }_ k , \Delta g^{ X, Y }_k \Delta W_k \rangle
%     + 2 \sum_{j_1,j_2=1}^m \langle \mathcal{J}^{ X, Y }_ k ,
%              \mathcal{L}^{j_1} \Delta g^{ X, Y}_{j_2,k} I^k_ {  j_1 , j_2  } \rangle \\
%     & \quad +
%     2 h \langle \Delta f^{ X, Y }_k , \Delta g^{ X, Y }_k \Delta W_k \rangle  
%    + 2 h \sum_{j_1,j_2=1}^m \langle \Delta f^{ X, Y }_k ,
%              \mathcal{L}^{j_1} \Delta g^{ X, Y}_{j_2,k} I^k_ {  j_1 , j_2  } \rangle \\
%   & \quad +
%     2 \sum_{j_1,j_2=1}^m \langle \Delta g^{ X, Y }_k \Delta W_k ,
%              \mathcal{L}^{j_1} \Delta g^{ X, Y}_{j_2,k} I^k_ {  j_1 , j_2  } \rangle
%\end{split}
%\end{equation}
Equipped with these estimates and taking expectations on both sides of \eqref{eq:error-identity-0},
one can derive
\begin{equation}\label{eq:error-expect-identity}
\begin{split}
  \E [\| \mathcal {J}^ { X, Y }_ { k+1 }  \| ^2] 
  &
  =
    \E [ \| \mathcal {J}^ { X, Y }_ { k }  \| ^2 ] + h ^2 \E [ \| \Delta f^{ X, Y }_k \| ^2 ]
    + h \E [ \| \Delta g^{ X, Y }_k  \| ^2 ]
 \\
  & \quad 
    + \tfrac{h^2}{2} \sum^m_{j_1,j_2=1} \E \big[ \big\| 
             \Delta ( \mathcal{L}^{j_1} g_{j_2})_k^{X, Y}   \big \| ^2 \big]
    + \E [  \| R_{ k+1 }  \| ^2 ]
    +
     2 h  \E [ \langle \mathcal{J}^{ X, Y }_ k , \Delta f^{ X, Y }_k \rangle ]
 \\
  & \quad 
     + 2 \E [  \langle \mathcal{J}^{ X, Y }_ k , R_{ k+1 } \rangle ]
    + 2 h \E [ \langle \Delta f^{ X, Y }_k , R_{ k+1 }  \rangle ]
    %      +
%    2  \E [ \langle \Delta g^{ X, Y }_k \Delta W_k ,
%             \sum^m_{j_1,j_2=1}
%             \Delta ( \mathcal{L}^{j_1} g_{j_2})_k^{X, Y}  I_{j_1,j_2}^{t_k,t_{k+1}} \rangle ]
  + 2 \E [ \langle \Delta g^{ X, Y }_k \Delta W_k , R_{ k+1 }  \rangle ]
    \\
  & \quad 
  +
     2  \E \Big[ \Big \langle
              \sum^m_{j_1,j_2=1}
             \Delta ( \mathcal{L}^{j_1} g_{j_2})_k^{X, Y}  I_{j_1,j_2}^{t_k,t_{k+1}} ,
    R_{ k+1 }  \Big \rangle \Big]
.
\end{split}
\end{equation}
Owing to the assumption that 
$ \E [\| \mathcal {J}^ { X, Y }_ { k }  \| ^2] < \infty $  for some $k \in \{ 0, 1, 2, ..., N-1\}$ 
and its consequence \eqref{eq:g-f-g'g-bound} 
as well as \eqref{eq:Rk-bound}, one can use the Cauchy-Schwarz inequality to infer
\begin{equation}
\begin{split}
\E [\| \mathcal {J}^ { X, Y }_ { k+1 }  \| ^2] 
  & \leq
          3 \E [ \| \mathcal {J}^ { X, Y }_ { k }  \| ^2 ] + 3 h ^2 \E [ \| \Delta f^{ X, Y }_k \| ^2 ]
    + 2 h \E [ \| \Delta g^{ X, Y }_k  \| ^2 ]
 \\
  & \quad 
  + h^2 \sum^m_{j_1,j_2=1} \E \big[ \big\| 
             \Delta ( \mathcal{L}^{j_1} g_{j_2})_k^{X, Y} \big \| ^2 \big]
   + 5 \E [  \| R_{ k+1 }  \| ^2 ]
   < \infty
  .
  \end{split}
\end{equation}
Based on the induction argument, the assertion 
$ \E [\| \mathcal {J}^ { X, Y }_ { k }  \| ^2] < \infty $
holds for all $ k \in \{ 0,1,2,\cdots, N \} $.
Following the same lines as used in \eqref{ineq:bound}-\eqref{eq:g-f-g'g-bound},
the boundedness of $ \E [\| \mathcal {J}^ { X, Y }_ { n }  \| ^2] $ for all $n \in \{ 0, 1, 2, ..., N\}$
ensures that
\begin{equation}
\E [ \| e_{k} \| ^ 2 ] < \infty , 
\:\:
\E [\| \Delta g_k^{X, Y} \|^2 ] < \infty, 
\:\:
\sum_{j_1,j_2=1}^m 
   \E \big[
    \big\| 
   \Delta \big(\mathcal{L}^{j_1}g_{j_2} \big)_k^{ X, Y} \big\| ^2
   \big] < \infty,
\:\:
\E [ \| \Delta f^{ X, Y }_k \|^2 ]  < \infty
\end{equation}
hold for all $ k \in \{ 0,1,2,\cdots, N \} $.
The desired assertion are thus justified.
\qed
%

%\begin{cor}\label{cor:MS-error-bounded}
% For $ \theta \in [ 0, 1 ] $ and $ h < \frac{ 1 }{ 2 L_4 } $ for $ \theta \neq 0 $ and $ 0 < h < \infty $ for $ \theta = 0 $ ,
%it holds for all $ k \in \{ 1,2,...,N \}  $
%that
%\begin{equation}
%\begin{split}
%& \E [ \| e_k - \theta \Delta f^{ X, Y }_k h
%      + { \tfrac{ \eta }{2} } \Delta g'g^{ X, Y }_k h \| ^ 2 ] < \infty , \quad
% \E [ \| e_k \| ^ 2 ] < \infty ,
%\\
% & \E [ \| \Delta f^{ X, Y }_k \| ^ 2 ] < \infty , \quad
% \E [ \| \Delta g^{ X, Y }_k \| ^ 2 ] < \infty , \quad
%%\\
%% &
% \E [ \| \Delta g'g^{ X, Y }_k \| ^ 2 ] < \infty .
%\end{split}
%\end{equation}
%\end{cor}
%
%{\it Proof of Corollary \ref{cor:MS-error-bounded}.}
%Further, taking into account the following facts
%\begin{equation*}
%\begin{split}
% \E [ \| X_{ t_{ k+1 } } - \theta f ( X _ { t_{ k+1 }} )h
% + { \tfrac{ \eta }{2} } g'g ( X _ { t_{ k+1 } } ) h \| ^2 ] < \infty , \quad
% \E [ \| X_{ t_{ k+1 } } \| ^2 ] < \infty ,
%\\
% \E [ \| f ( X _ { t_{ k+1 }} ) \| ^2 ] < \infty , \quad
% \E [ \| \| g ( X _ { t_{ k+1 }} ) \| ^2 ] < \infty ,\quad
% \E [ \| g'g ( X _ { t_{ k+1 }} ) \| ^2 ] < \infty
%\end{split}
%\end{equation*}
%for all $ k \in \{ 1,2,...,N \} $  , one can readily arrive at the desired assertions.
%\\
%\qed
%

Before proceeding further, we point out that the moment bounds in \eqref{eq:lemma-numerical-bound}, depending on $N$,
are not proved to be uniformly bounded with respect to $N$. 
This means that the moment bounds might depend on the step number $N$.
However, such moment bounds are enough 
for the subsequent error analysis, which does not rely on the precise uniform moment bounds of the numerical approximations.
Now we are well prepared to prove Theorem \ref{thm:upper-error-bound}.

{\it Proof of Theorem \ref{thm:upper-error-bound}.}
%
%Noticing that
%\begin{equation}
%\begin{split}
%  \E [ \| I _ { ( 1 , 1 ) } \| ^2 ]
%  =
%  \E [ \| { \tfrac{ 1 }{ 2 } } ( \Delta W_k^2 - h ) \| ^ 2 ]
%  =
%  { \tfrac{ 3 }{ 4 } } h^2
%   - { \tfrac{ 1 }{ 2 } } h^2
%   + { \tfrac{ 1 }{ 4 } } h^2
% =   { \tfrac{ 1 }{ 2 } } h^2 ,
%\end{split}
%\end{equation}
%Bearing in mind that 
Recalling
$
 \mathcal{J}^{ X, Y }_ k : =
 e_k - \theta \Delta f^{ X, Y }_k h
 +
      \tfrac {\eta} {2} \sum_{ j = 1}^{ m} \Delta ( \mathcal{L}^{j} g_{j})_{k}^{X, Y} h
$
and using its consequence
$
 \Delta f^{ X, Y }_k h   =
{ \tfrac{ 1 }{\theta} } ( e_k - \mathcal{J}^{ X, Y }_ k
      +  \tfrac {\eta} {2} \sum_{ j = 1}^{ m} \Delta ( \mathcal{L}^{j} g_{j})_{k}^{X, Y} h )
$
and \eqref{eq:lemma-numerical-bound},
we derive from \eqref{eq:error-expect-identity} that, for any $ k \in \{ 0,1,2,\cdots, N-1 \} $,
\begin{equation}
 \begin{split}
 \E [ \| \mathcal{J}^{ X, Y }_{ k+1 } \| ^2 ]
  & =
      \E [ \| \mathcal {J}^ { X, Y }_ { k }  \| ^2 ] + (1 - 2 \theta ) h ^2 \E [ \| \Delta f^{ X, Y }_k \| ^2 ]
    + h \E [ \| \Delta g^{ X, Y }_k  \| ^2 ]
 \\
  & \quad 
  + \tfrac{h^2}{2} \sum^m_{j_1,j_2=1} \E \big[ \| 
             \Delta ( \mathcal{L}^{j_1} g_{j_2})_k^{X, Y}  \| ^2 \big]
   + \E [  \| R_{ k+1 }  \| ^2 ]
  +
     2 h  \E [ \langle e_ k , \Delta f^{ X, Y }_k \rangle ]
     \\
     & \quad 
     + \eta h^2 \E \Big[ \Big\langle \sum^m_{j=1}
             \Delta ( \mathcal{L}^{j} g_{j})_k^{X, Y} , \Delta f^{ X, Y }_k \Big \rangle \Big]
     + 
     2 \E [  \langle \mathcal{J}^{ X, Y }_ k , R_{ k+1 } \rangle ]
    + \tfrac{2}{\theta}  \E [ \langle  e_k , R_{ k+1 }  \rangle ]
    \\
  & \quad 
%      +
%    2  \E [ \langle \Delta g^{ X, Y }_k \Delta W_k ,
%             \sum^m_{j_1,j_2=1}
%             \Delta ( \mathcal{L}^{j_1} g_{j_2})_k^{X, Y}  I_{j_1,j_2}^{t_k,t_{k+1}} \rangle ]
- \tfrac{2}{\theta} \E [ \langle \mathcal {J}^ { X, Y }_ { k } , R_{ k+1 }  \rangle ]
    + \tfrac{\eta h }{\theta} \E \Big[ \Big \langle \sum^m_{j=1}
             \Delta ( \mathcal{L}^{j} g_{j})_k^{X, Y}  , R_{ k+1 }  \Big \rangle \Big ]
  \\
  & \quad
  + 2 \E [ \langle \Delta g^{ X, Y }_k \Delta W_k , R_{ k+1 }  \rangle ]
%\\ & \quad
  +
     2  \E \Big [ \Big \langle
              \sum^m_{j_1,j_2=1}
             \Delta ( \mathcal{L}^{j_1} g_{j_2})_k^{X, Y}  I_{j_1,j_2}^{t_k,t_{k+1}} ,
    R_{ k+1 } \Big \rangle \Big]
%    \\
%    & @@=
%      \E [ \| \mathcal {J}^ { X, Y }_ { k }  \| ^2 ] + (1 - 2 \theta ) h ^2 \E [ \| \Delta f^{ X, Y }_k \| ^2 ]
%    + h \E [ \| \Delta g^{ X, Y }_k  \| ^2 ]
%    + \tfrac{h^2}{2} \sum^m_{j_1,j_2=1} \E \Big[ \| 
%             \Delta ( \mathcal{L}^{j_1} g_{j_2})_k^{X, Y}  \| ^2 \Big]
% \\
%  & \quad 
%   + \E [  \| R_{ k+1 }  \| ^2 ]
%  +
%     2 h  \E [ \langle e_ k , \Delta f^{ X, Y }_k \rangle ]
%     +
%     \tfrac{ \eta h}{\theta} \E[ \langle \sum^m_{j=1}
%             \Delta ( \mathcal{L}^{j} g_{j})_k^{X, Y} , e_k \rangle ]
%     \\
%     & \quad
%     -
%     \tfrac{ \eta h}{\theta} \E[ \langle \sum^m_{j=1}
%             \Delta ( \mathcal{L}^{j} g_{j})_k^{X, Y} ,  \mathcal {J}^ { X, Y }_ { k } \rangle ]
%     +
%      \tfrac{ \eta^2 h^2}{ 2 \theta} \E[ \| \sum^m_{j=1}
%             \Delta ( \mathcal{L}^{j} g_{j})_k^{X, Y} \|^2 ]
%     \\
%     & \quad 
%     + 2 \E [  \langle \mathcal{J}^{ X, Y }_ k , R_{ k+1 } \rangle ]
%    + \tfrac{2}{\theta}  \E [ \langle  e_k , R_{ k+1 }  \rangle ]
%    - \tfrac{2}{\theta} \E [ \langle \mathcal {J}^ { X, Y }_ { k } , R_{ k+1 }  \rangle ]
%    + \tfrac{\eta h }{\theta} \E [ \langle \sum^m_{j=1}
%             \Delta ( \mathcal{L}^{j} g_{j})_k^{X, Y}  , R_{ k+1 }  \rangle ]
%    \\
%  & \quad 
%  + 2 \E [ \langle \Delta g^{ X, Y }_k \Delta W_k , R_{ k+1 }  \rangle ]
%  +
%     2  \E [ \langle
%              \sum^m_{j_1,j_2=1}
%             \Delta ( \mathcal{L}^{j_1} g_{j_2})_k^{X, Y}  I_{j_1,j_2}^{t_k,t_{k+1}} ,
%    R_{ k+1 }  \rangle]
%    @@
    \\
    & =
      \E [ \| \mathcal {J}^ { X, Y }_ { k }  \| ^2 ] + (1 - 2 \theta ) h ^2 \E [ \| \Delta f^{ X, Y }_k \| ^2 ]
    + h \E [ \| \Delta g^{ X, Y }_k  \| ^2 ]
 \\
  & \quad 
  + \tfrac{ h^2 }{2} \sum^m_{j_1,j_2=1}  \E \Big[ \|
             \Delta ( \mathcal{L}^{j_1} g_{j_2})_k^{X, Y} \| ^2 \Big]
   + \E [  \| R_{ k+1 }  \| ^2 ]
  +
     2 h  \E [ \langle e_ k , \Delta f^{ X, Y }_k \rangle ]
     \\
     & \quad 
     +
     \eta h^2 \E \Big [ \Big \langle \sum^m_{j=1}
             \Delta ( \mathcal{L}^{j} g_{j})_k^{X, Y} , \Delta f^{ X, Y }_k \Big \rangle \Big ]
     + \tfrac{ 2 \theta - 2}{\theta} \E [  \langle \mathcal{J}^{ X, Y }_ k , \E ( R_{ k+1 } | \mathcal{F}_{t_k} ) \rangle ]
%    - \tfrac{2}{\theta} \E [ \langle \mathcal {J}^ { X, Y }_ { k } , \E ( R_{ k+1 } | \mathcal{F}_{t_k} )   \rangle ]
    \\
  & \quad 
%      +
%    2  \E [ \langle \Delta g^{ X, Y }_k \Delta W_k ,
%             \sum^m_{j_1,j_2=1}
%             \Delta ( \mathcal{L}^{j_1} g_{j_2})_k^{X, Y}  I_{j_1,j_2}^{t_k,t_{k+1}} \rangle ]
+ \tfrac{2}{\theta}  \E [ \langle  e_k , \E ( R_{ k+1 } | \mathcal{F}_{t_k} )   \rangle ]
    + \tfrac{\eta h }{\theta} \E \Big[ \Big \langle \sum^m_{j=1}
             \Delta ( \mathcal{L}^{j} g_{j})_k^{X, Y}  , R_{ k+1 }  \Big \rangle \Big ]
\\ & \quad
+ 2 \E [ \langle \Delta g^{ X, Y }_k \Delta W_k , R_{ k+1 }  \rangle ]
  +
     2  \E \Big[ \Big \langle
              \sum^m_{j_1,j_2=1}
             \Delta ( \mathcal{L}^{j_1} g_{j_2})_k^{X, Y}  I_{j_1,j_2}^{t_k,t_{k+1}} ,
    R_{ k+1 }  \Big \rangle \Big].
 \end{split}
\end{equation}
Using the Cauchy-Schwarz inequality and the Young inequality gives
\begin{equation}
\begin{split}
\tfrac{ 2 \theta - 2}{\theta} \E [  \langle \mathcal{J}^{ X, Y }_ k , \E ( R_{ k+1 } | \mathcal{F}_{t_k} ) \rangle ]
& \leq
h \E [ \| \mathcal{J}^{ X, Y }_ k \|^2 ] 
+
\tfrac{(\theta - 1)^2}{\theta^2 h} \E [ \| \E ( R_{ k+1 } | \mathcal{F}_{t_k} ) \|^2 ],
\\
\tfrac{2}{\theta}  \E [ \langle  e_k , \E ( R_{ k+1 } | \mathcal{F}_{t_k} )   \rangle ]
& \leq
h \E [ \| e_ k \|^2 ] 
+
\tfrac{1}{\theta^2 h} \E [ \| \E ( R_{ k+1 } | \mathcal{F}_{t_k} ) \|^2 ],
\\
\tfrac{\eta h }{\theta} \E \Big[ \Big \langle \sum^m_{j=1}
             \Delta ( \mathcal{L}^{j} g_{j})_k^{X, Y}  , R_{ k+1 }  \Big \rangle \Big ]
& \leq
\tfrac{\varrho-1}{4} h^2 
\sum^m_{j=1}
            \E [ \| \Delta ( \mathcal{L}^{j} g_{j})_k^{X, Y} \|^2 ]
    +
    \tfrac{\eta^2m}{ (\varrho-1) \theta^2 }
    \E [\| R_{k+1} \|^2],
    \nonumber
\end{split}
\end{equation}
%
%and
%\begin{equation}
%\begin{split}
\begin{align}
2 \E [ \langle \Delta g^{ X, Y }_k \Delta W_k , R_{ k+1 }  \rangle ] 
& \leq
(q - 2) h \E [ \| \Delta g^{ X, Y }_k \|^2 ] 
+
\tfrac{1}{q-2}\E [\| R_{k+1} \|^2],
\\
2  \E \Big[ \Big \langle
              \sum^m_{j_1,j_2=1}
             \Delta ( \mathcal{L}^{j_1} g_{j_2})_k^{X, Y}  I_{j_1,j_2}^{t_k,t_{k+1}} ,
    R_{ k+1 }  \Big \rangle \Big]
& \leq
\tfrac{\varrho -1}{4} h^2
\sum^m_{j_1,j_2=1} \E \big[ \big\| \Delta ( \mathcal{L}^{j_1} g_{j_2})_k^{X, Y}    \big \|^2 \big]
%\\ & \quad
+
\tfrac{2}{\varrho-1}\E [\| R_{k+1} \|^2].
\nonumber
%\\
%&=
%\tfrac{\varrho -1}{4} h^2
%\sum^m_{j_1,j_2=1} \E \Big[ \Big\| \Delta ( \mathcal{L}^{j_1} g_{j_2})_k^{X, Y}  \Big \|^2 \Big]
%+
%\tfrac{2}{\varrho-1}\E [\| R_{k+1} \|^2].
%\end{split}
%\end{equation}
\end{align}
Taking these estimates into consideration  and recalling \eqref{eq:mono-condition1} yield
\begin{equation}
 \begin{split}
 \E [ \| \mathcal{J}^{ X, Y }_{ k+1 } \| ^2 ]
& \leq
    ( 1 + h)  \E [ \| \mathcal {J}^ { X, Y }_ { k }  \| ^2 ] + (1 - 2 \theta ) h ^2 \E [ \| \Delta f^{ X, Y }_k \| ^2 ]
    + h (q - 1) \E [ \| \Delta g^{ X, Y }_k  \| ^2 ]
  \\ & \quad  
    + \tfrac{ \varrho h^2 }{2} \sum^m_{j_1,j_2=1}  \E \big[ \|
             \Delta ( \mathcal{L}^{j_1} g_{j_2})_k^{X, Y} \| ^2 \big]
    +
     \eta h^2 \E \Big [ \Big \langle \sum^m_{j=1}
             \Delta ( \mathcal{L}^{j} g_{j})_k^{X, Y} , \Delta f^{ X, Y }_k \Big \rangle \Big ]
 \\
  & \quad 
    + 
      ( \tfrac{q-1}{q-2} + \tfrac{\eta^2m}{(\varrho-1) \theta^2 } + \tfrac{2}{\varrho-1} ) \E [  \| R_{ k+1 }  \| ^2 ]
      +
     2 h  \E [ \langle e_ k , \Delta f^{ X, Y }_k \rangle ]
 \\
     & \quad
      +
    h \E [ \| e_k \|^2 ]
     +
     \tfrac{\theta^2 - 2 \theta + 2}{\theta^2 h} \E [ \| \E ( R_{ k+1 } | \mathcal{F}_{t_k} ) \|^2]
%     \\
%     & \quad 
%      +
%    2  \E [ \langle \Delta g^{ X, Y }_k \Delta W_k ,
%             \sum^m_{j_1,j_2=1}
%             \Delta ( \mathcal{L}^{j_1} g_{j_2})_k^{X, Y}  I_{j_1,j_2}^{t_k,t_{k+1}} \rangle ]
%    + \tfrac{\eta h }{\theta} \E \Big[ \Big \langle \sum^m_{j=1}
%             \Delta ( \mathcal{L}^{j} g_{j})_k^{X, Y}  , R_{ k+1 }  \Big \rangle \Big]
%\\ & \quad
%  +
%     2  \E \Big[ \Big \langle
%              \sum^m_{j_1,j_2=1}
%             \Delta ( \mathcal{L}^{j_1} g_{j_2})_k^{X, Y}  I_{j_1,j_2}^{t_k,t_{k+1}} ,
%    R_{ k+1 }  \Big \rangle \Big]
\\
& \leq
 ( 1 + h)  \E [ \| \mathcal {J}^ { X, Y }_ { k }  \| ^2 ]  + ( 1 + L_1) h \E [ \| e_k \|^2 ] 
\\ & \quad
 +
 \big ( \tfrac{q-1}{q-2} + \tfrac{\eta^2m}{(\varrho-1) \theta^2 } + \tfrac{2}{\varrho-1} \big) \E [  \| R_{ k+1 }  \| ^2 ] 
 +  
 \tfrac{\theta^2 - 2 \theta + 2}{\theta^2 h} \E [ \| \E ( R_{ k+1 } | \mathcal{F}_{t_k} ) \|^2] 
 .
 \end{split}
\end{equation}
%
%Then we take conditional expectation of both sides ,
%\begin{equation}
% \begin{split}
% \E [ \| \mathcal{J}^{ X, Y }_{ k+1 } \| ^2 ]
% \leq &
% \E [ \| \mathcal{J}^{ X, Y }_ k  \| ^2 ]
%    + L_1 h \E [ \| e_k \| ^2 ]
%    + ( \tfrac{q -  1 }{q - 2 } + \tfrac{ 1 }{ \varepsilon_1 } + \tfrac{ \eta ^2 }{ \theta ^2 }
%       \tfrac{ 1 }{ 2 \varepsilon_2 } ) \E [ \| R_{ k+1 } \| ^2 ]
% \\
% &
%    - \tfrac{ 2 }{ \theta } \E [ \langle e_k , \E( R_{ k+1 } | \mathcal{F}_k )  \rangle ]
%    +  2( \tfrac{ 1 }{ \theta } - 1 ) \E [ \langle \mathcal{J}^{ X, Y }_ k  , \E( R_{ k+1 } | \mathcal{F}_k )  \rangle ]
% \\
% \leq &
% \E [ \| \mathcal{J}^{ X, Y }_ k  \| ^2 ]
%    + L_1 h \E [ \| e_k \| ^2 ]
%    + ( \tfrac{q -  1 }{q - 2 } + \tfrac{ 1 }{ \varepsilon_1 } + \tfrac{ \eta ^2 }{ \theta ^2 }
%       \tfrac{ 1 }{ 2 \varepsilon_2 } ) \E [ \| R_{ k+1 } \| ^2 ]
% \\
% &
%    + h \E [ \| e_k \| ^2 ] + \tfrac{ 1 }{ \theta ^ 2 h } \E [ \| \E( R_{ k+1 } | \mathcal{F}_k ) \| ^2 ]
%    + h \E [ \| \mathcal{J}^{ X, Y }_ k \| ^2 ] + \tfrac{ (1 - \theta) ^ 2 }{ \theta ^ 2 h } \E [ \| \E( R_{ k+1 } | \mathcal{F}_k ) \| ^2 ]
% \\
% \leq &
%    ( 1 + h ) \E [ \| \mathcal{J}^{ X, Y }_ k  \| ^2 ]
%    + ( 1 + L_1 ) h \E [ \| e_k \| ^2 ]
%    + ( \tfrac{q -  1 }{q - 2 } + \tfrac{ 1 }{ \varepsilon_1 } + \tfrac{ \eta ^2 }{ \theta ^2 }
%       \tfrac{ 1 }{ 2 \varepsilon_2 } ) \E [ \| R_{ k+1 } \| ^2 ]
% \\
% &
%    + \tfrac{ \theta ^ 2 + 1 }{ \theta ^ 2 h } \E [ \| \E( R_{ k+1 } | \mathcal{F}_k ) \| ^2 ]
%    .
% \end{split}
%\end{equation}
%
By iteration and observing  $ \mathcal{J}^{ X, Y }_ 0 = e_0 - \theta \Delta f^{ X, Y }_0 h
                             + \tfrac {\eta} {2} h \sum_{ j = 1}^{ m} \Delta ( \mathcal{L}^{j} g_{j})_{0}^{X, Y}  = 0 $ we deduce
\begin{equation} \label{eq:estimate-Jk-iteration}
 \begin{split}
 \E [ \| \mathcal{J}^{ X, Y }_{ k+1 } \| ^2 ]
 \leq &
%    ( 1 + h )^{ k+1 } \E [ \| \mathcal{J}^{ X, Y }_ 0  \| ^2 ]
%    + 
    ( 1 + L_1 ) h \sum_{ i = 0 }^{ k }( 1 + h )^{ ( k - i ) } \E [ \| e_i \| ^2 ]
\\
&    
    + ( \tfrac{q-1}{q-2} + \tfrac{\eta^2m}{(\varrho-1) \theta^2 } + \tfrac{2}{\varrho-1} )
       \sum_{ i = 0 }^{ k } ( 1 + h )^{ ( k - i ) } \E [ \| R_{ i + 1} \| ^2 ]
 \\
 &
    + \tfrac{ \theta ^ 2 - 2 \theta + 2 }{ \theta ^ 2 h }
    \sum_{ i = 0 }^{ k }( 1 + h )^{ ( k - i ) }
    \E [ \| \E( R_{ i+1 } | \mathcal{F}_{t_i} ) \| ^2 ]
 \\
\leq &
 ( 1 + L_1 ) h e ^ T \sum_{ i = 0 }^{ k } \E [ \| e_i \| ^2 ]
    + ( \tfrac{q-1}{q-2} + \tfrac{\eta^2m}{ (\varrho-1) \theta^2 } + \tfrac{2}{\varrho-1} )
       e ^ T \sum_{ i = 0 }^{ k } \E [ \| R_{i + 1} \| ^2 ]
\\ &
    + \tfrac{ \theta ^ 2 - 2 \theta + 2 }{ \theta ^ 2 h } e ^ T 
    \sum_{ i = 0 }^{ k }
    \E [ \| \E( R_{ i+1 } | \mathcal{F}_{t_i} ) \| ^2 ]
    .
 \end{split}
\end{equation}
Additionally, the assumption \eqref{eq:mono-condition2} ensures
\begin{equation}
\begin{split}
\E [ \| \mathcal{J}^{ X, Y }_{ k+1 } \| ^2 ]
 &
 =
 \E \Big[ \Big\| e_{k+1} - \theta \Delta f^{ X, Y }_{k+1} h
       +
      \tfrac {\eta} {2} \sum_{ j = 1}^{ m} \Delta ( \mathcal{L}^{j} g_{j})_{k+1}^{X, Y} h \Big \| ^2 \Big]
\\
 & =
  \E [ \| e_{k+1} \| ^2 ] 
  +  
  \E \Big [ \Big \| \theta  h  \Delta f^{ X, Y }_{k+1} 
  - 
      \tfrac {\eta} {2} \sum_{ j = 1}^{ m} \Delta ( \mathcal{L}^{j} g_{j})_{k+1}^{X, Y} h \Big \| ^2 \Big]
\\
& \quad
  - 2 h \E \Big [ \Big \langle e_{k+1} , \theta \Delta f^{ X, Y }_{k+1} - 
      \tfrac {\eta} {2} \sum_{ j = 1}^{ m} \Delta ( \mathcal{L}^{j} g_{j})_{k+1}^{X, Y} h  \Big \rangle \Big ]
\\
 & \geq
  ( 1 - 2 L_2 h ) \E [ \| e_{k+1} \| ^2 ]
%  + c_0 h ^ 2
  .
 \end{split}
\end{equation}
Inserting this into \eqref{eq:estimate-Jk-iteration} yields
\begin{equation}
\begin{split}
( 1 - 2 L_2 h ) \E [ \| e_{k+1} \| ^2 ]
 \leq &
 ( 1 + L_1 ) h e ^ T \sum_{ i = 0 }^{ k  } \E [ \| e_i \| ^2 ]
    + ( \tfrac{q-1}{q-2} + \tfrac{\eta^2m}{ (\varrho-1) \theta^2 } + \tfrac{2}{\varrho-1} )
       e ^ T \sum_{ i = 0 }^{ k } \E [ \| R_{i+1} \| ^2 ]
 \\
 &
    + \tfrac{ \theta ^ 2 - 2 \theta +2 }{ \theta ^ 2 h } e ^ T \sum_{ i = 0 }^{ k }
    \E [ \| \E( R_{ i+1 } | \mathcal{F}_{t_i} ) \| ^2 ].
 \end{split}
\end{equation}
Owing to $ 2  L_2  h \leq \nu < 1$ by assumption and bearing the moment bounds \eqref{eq:lemma-numerical-bound} in mind,
%$ h \in (0,  \frac{ 1 }{ 2 L_2 } ) $, 
one can apply Gronwall's inequality to acquire the desired assertion.
%
%
%\begin{equation}
%\begin{split}
% \E [ \| e_{k+1} \| ^2 ]
% \leq &
%    C \big[ ( \tfrac{q -  1 }{q - 2 } + ? + ? )
%       e ^ T \sum_{ i = 1 }^{ k + 1 } \E [ \| R_i \| ^2 ] \big]
%    + \tfrac{ \theta ^ 2 + 1 }{ \theta ^ 2 h } e ^ T \sum_{ i = 1 }^{ k + 1 }
%    \E [ \| \E( R_{ k+1 } |\, \mathcal{F}_{t_k} ) \| ^2 ]
%%    - \tfrac{ c_0 h ^ 2 }{ 1 - 2 L_2 h }
%,
% \end{split}
%\end{equation}
%as required.
\qed
%%%%%%%
%%%%
%

%It is worthwhile to point out that, 
{\color{black}It is worthwhile to point out that,
conditions in Assumption \ref{ass:monotonicity-condition}
are not difficult to be fulfilled. For instance, the following assumption suffices to imply Assumption \ref{ass:monotonicity-condition}.}
\begin{assumption}%[Monotonicity type conditions in a domain]
\label{ass:cor-monotonicity-condition-theta-Milstein}
Assume that the diffusion coefficients $ g_j \colon \mathbb{R}^d \rightarrow \mathbb{R}^{d }, j \in \{ 1, 2,..., m \}  $ 
are differentiable in a domain $D \subset \R^d$. 
%for method parameters $\theta, \eta \in [0, 1]$
There exist constants $q \in (2, \infty) $, $ \zeta \in (0, \infty) $ and $L_3 \in [0, \infty)$ 
%and $ h_0 \in (0, T]$ 
such that,
$\forall x, y \in D$, $ h \in (0, 2 \zeta)$, the drift and diffusion coefficients of SDEs \eqref{eq:SODE} obey
 \begin{equation}
 \label{eq:cor-mono-condition1}
% \begin{split}
 2 \langle x - y ,f ( x ) - f ( y ) \rangle
  +  ( q - 1 ) \|g ( x ) - g ( y )\| ^2
  + \zeta \sum_{j_1,j_2=1}^m \big\| 
   \mathcal{L}^{j_1}g_{j_2}( x ) - \mathcal{L}^{j_1}g_{j_2}( y ) \big\| ^2
   \leq L _ 3 \| x - y \|^2
.
% \end{split}
 \end{equation}
 %
 %%
% \begin{equation}\label{ineq:mono4}
% \begin{split}
% &
% \langle x ,\theta f ( x )- \tfrac{\eta}{2} g'g ( x ) \rangle
%   \leq L _ 4 ( 1 + \| x \|^2 ) ,
% \end{split}
% \end{equation}
 %
 %%
\end{assumption}
{\color{black}
We mention that  such a condition was also used in \cite[Theorem 2.3]{beyn2017stochastic} 
for the backward Milstein method ($\theta = 1, \eta = 0$) and $D = \R^d$.
It is not difficult to check that, when the above condition \eqref{eq:cor-mono-condition1} holds, 
all conditions in Assumption \ref{ass:monotonicity-condition} are satisfied with $\theta \in [\tfrac12, 1], \eta = 0$, 
$L_1 = L_3$ and $L_2 = \tfrac{\theta L_3}{2}$.}
As a direct consequence of Theorem \ref{thm:upper-error-bound}, we get the following corollary.
\begin{cor}\label{cor:error-bounds-ass1}
%[An upper mean-square error bound]
\label{cor:theta-Milstein-upper-error-bound}
Let Assumptions \ref{ass:well-possedness}, \ref{ass:cor-monotonicity-condition-theta-Milstein} be fulfilled
with $ \theta L_3 h \leq \nu $ for some $\nu \in (0, 1)$ and $\theta \in [\tfrac12, 1], \eta = 0$.
Let $\{ X_{ t} \}_{ t \in [0, T]}$ and $ \{Y_n\}_{0\leq n\leq N} $ be solutions to 
SDEs \eqref{eq:SODE} and the semi-implicit Milstein method \eqref{eq:semi-implicit-scheme}, respectively.
Then the mean-square error upper bound \eqref{eq:Thm-Error-Bound} holds.
%%
%\begin{equation}
%\begin{split}
%\E \big [ \big \lVert X_{ t_n } - Y_n \big \rVert ^2 \big ]
%\leq
%C
%\Big(  
%\sum_{ i = 1 }^{ n } \E [ \| R_i \| ^2 ]
%    + \tfrac{ 1 }{ h }
%    \sum_{ i = 1 }^{ n }
%    \E [ \| \E( R_{ i } \, | \mathcal{F}_{ i-1} ) \| ^2 ]
% \Big)
%,
%\end{split}
%\end{equation}
%where $R_n, n\in \{ 1,2,..., N\}$ is defined as \eqref{eq:Error-Remainder-Defn}.
\end{cor}

Observe that the condition \eqref{eq:cor-mono-condition1} would impose a strict restriction 
on the polynomial growth of the diffusion coefficient, which excludes practical models
such as the $\tfrac32$-volatility model \eqref{eq:intro-32model} and the Ait Sahalia model 
\eqref{eq:intro-Ait-Sahalia-model-SDE}. 
This can be remedied by utilizing the following assumption.

\begin{assumption}%[Monotonicity type conditions in a domain]
\label{ass:cor-monotonicity-condition-backward-Milstein}
Assume that the diffusion coefficients $ g_j \colon D \rightarrow \mathbb{R}^{d }, j \in \{ 1, 2,..., m \}  $ 
are differentiable in the domain $D \subset \R^d$. For method parameters $\theta \in (\tfrac12, 1], \eta \in [0, 1]$,
there exist constants $q \in (2, \infty) $, $ \varrho \in (1, \infty) $ and $L_4, L_5, L_6 \in [0, \infty)$ 
%and $ h_0 \in (0, T]$ 
such that,
$\forall x, y \in D$,
the drift and diffusion coefficients of SDEs \eqref{eq:SODE} obey
 \begin{equation}
 \label{eq:cor2-mono-condition1}
\begin{split}
 2 \langle x - y ,f ( x ) - f ( y ) \rangle
  +  ( q - 1 ) \|g ( x ) - g ( y )\| ^2
   \leq L _ 4 \| x - y \|^2,
   \\
\tfrac{ \varrho }{ 2 } \sum_{j_1,j_2=1}^m \big\| 
   \mathcal{L}^{j_1}g_{j_2}( x ) - \mathcal{L}^{j_1}g_{j_2}( y ) \big\| ^2
  + \eta \Big \langle  \sum_{j=1}^m 
  \big[ \mathcal{L}^{j}g_{j}( x )
   - \mathcal{L}^{j}g_{j}( y ) \big] , f ( x ) - f ( y ) \Big \rangle
   \\
  +  ( 1 - 2 \theta ) \| f ( x ) - f ( y ) \|^2
   \leq L _ 5 \| x - y \|^2,
   \\
    \Big \langle x - y ,\theta [ f ( x ) - f ( y ) ]
 - \tfrac{\eta}{2} 
 \sum_{j=1}^m [ \mathcal{L}^{j}g_{j}( x )
  - \mathcal{L}^{j} g_{j}( y )]  \Big \rangle
   \leq L_6 \| x - y \|^2.
\end{split}
 \end{equation}
 %
 %%
% \begin{equation}\label{ineq:mono4}
% \begin{split}
% &
% \langle x ,\theta f ( x )- \tfrac{\eta}{2} g'g ( x ) \rangle
%   \leq L _ 4 ( 1 + \| x \|^2 ) ,
% \end{split}
% \end{equation}
 %
 %%
\end{assumption}
{\color{black}One can straightforwardly  verify that Assumption \ref{ass:cor-monotonicity-condition-backward-Milstein} implies 
Assumption \ref{ass:monotonicity-condition} and one gets the following corollary,
 as a direct consequence of Theorem \ref{thm:upper-error-bound}.}
\begin{cor}
%[An upper mean-square error bound]
\label{cor:backward-Milstein-upper-error-bound}
Let Assumptions \ref{ass:well-possedness}, \ref{ass:cor-monotonicity-condition-backward-Milstein} hold
with $\theta \in (\tfrac12, 1]$, $\eta \in [0,1]$ and $ 2 L_6 h \leq \nu$ for some $\nu \in (0, 1)$.
 Let $\{ X_{ t} \}_{ t \in [0, T]}$ and $ \{Y_n\}_{0\leq n\leq N} $ be solutions to 
SDEs \eqref{eq:SODE} and the double implicit Milstein method \eqref{eq:general-scheme}, respectively.
Then the mean-square error upper bound \eqref{eq:Thm-Error-Bound} holds.
%there exists a uniform constant $C$ such that, for any $ n \in \{ 1,2,..., N \}$, $N \in \N$,
%%
%\begin{equation}
%\begin{split}
%\E \big [ \big \lVert X_{ t_n } - Y_n \big \rVert ^2 \big ]
%\leq
%C
%%\big[ ( \tfrac{q -  1 }{q - 2 } + \tfrac{ 1 }{ \varepsilon_1 }
%%    + \tfrac{ \eta ^2 }{ \theta ^2 } \tfrac{ 1 }{ 2 \varepsilon_2 } )
%%       e ^ T \sum_{ i = 1 }^{ k } \E [ \| R_i \| ^2 ] \big]
%%    + \tfrac{ \theta ^ 2 + 1 }{ \theta ^ 2 h } e ^ T
%%    \sum_{ i = 1 }^{ k }
%%    \E [ \| \E( R_{ i } \, | \mathcal{F}_{ i-1} ) \| ^2 ]
%\Big(  
%\sum_{ i = 1 }^{ n } \E [ \| R_i \| ^2 ]
%    + \tfrac{ 1 }{ h }
%    \sum_{ i = 1 }^{ n }
%    \E [ \| \E( R_{ i } \, | \mathcal{F}_{ i-1} ) \| ^2 ]
% \Big)
%,
%\end{split}
%\end{equation}
%where $R_n, n\in \{ 1,2,..., N\}$ is defined as \eqref{eq:Error-Remainder-Defn}.
\end{cor}
%%%%%%%%%%%%%
%\begin{example}
%...
%\end{example}
%%%%%%%%%%%%%

In section \ref{sect:applications}, we will show that the above two financial models and their numerical schemes 
fulfill Assumption \ref{ass:cor-monotonicity-condition-backward-Milstein} and one can thus rely on 
Corollary \ref{cor:backward-Milstein-upper-error-bound} to obtain the desired convergence rate.
Before closing this section, we would like to mention that, the previously obtained mean-square error bound
\eqref{eq:Thm-Error-Bound} is powerful as it helps us to easily analyze mean-square convergence rates of the schemes, 
without relying on a priori high-order moment estimates of numerical approximations.
This will be seen in the forthcoming two sections, where we shall use the error bounds to recover the expected mean-square 
convergence rates of the proposed schemes in various circumstances.

\section{Mean-square convergence rates under globally polynomial growth conditions
%for SODEs with monotone coefficients
}
\label{sec:conv-rates-global-polynomial}
Equipped with the previously derived upper mean-square error bounds,
the present section aims to identify the expected mean-square convergence rate of 
the underlying schemes \eqref{eq:general-scheme} for SDEs in the whole space $\R^d$
under further globally polynomial assumptions.
To this end, we make the following globally polynomial growth and coercivity conditions on the drift and diffusion coefficients.
%%%%%%%%%%%%%
%%%%%%%%%%%%%
%\subsection{Convergence rates for SDEs with multiplicative noise}
%
%In order to analyze the convergence rates in the multiplicative noise case, we make the following
%polynomial growth assumption on the drift coefficient $f$.
%

%
\begin{assumption}
[Globally polynomial growth and coercivity conditions in $\R^d$]
\label{ass:f-polynomial-growth}
%Suppose that both SDEs \eqref{eq:SODE} and the schemes \eqref{eq:general-scheme} obey
%Assumption \ref{ass:monotonicity-condition} in the whole space $D = \R^d$.
%Furthermore, 
Assume both the drift coefficient $f \colon \R^d \rightarrow \R^d$ and the diffusion coefficients 
$g_j \colon \R^d \rightarrow \R^d, j \in \{ 1, 2,..., m \}$ of SDEs \eqref{eq:SODE} are 
%twice continuously differentiable, 
continuously differentiable in $\R^d$,  and
there exist some positive constants $\gamma \in [ 1, \infty )$ and $ p^* \in [ 6 \gamma - 4, \infty) $ such that,
\begin{align} 
\label{assum:eq-monotone-growth}
& \langle x , f ( x )  \rangle + \tfrac{ p^* - 1}{2} \lVert g ( x ) \rVert ^ 2   \leq C ( 1 +  \lVert x \rVert ^ 2 ),
\quad
\forall x \in \R^d,
\\
\label{eq:drift-polynomial-growth-condition}
\big \| \big ( \tfrac{ \partial f  } { \partial x}  ( x ) -  \tfrac{ \partial f  } { \partial x} ( \tilde{x} ) \big) y \big \| 
& 
\leq C ( 1 + \| x \| + \| \tilde{x} \| )^{ \gamma - 2} \| x - \tilde{x} \| \cdot \| y \|, 
%\quad i \in \{ 1, 2 \},
\quad
\forall x, \tilde{x}, y \in \R^d,
\\
\label{eq:diffusion-polynomial-growth-condition}
\big \| \big ( \tfrac{ \partial g_j  } { \partial x}  ( x ) -  \tfrac{ \partial g_j  } { \partial x} ( \tilde{x} ) \big) y \big \|^2 
& 
\leq C ( 1 + \| x \|+ \| \tilde{x} \| )^{ \gamma - 3}  \| x - \tilde{x} \|^2 \cdot \| y \|^2, 
%\quad i \in \{ 1, 2 \},
\quad
\forall x, \tilde{x}, y \in \R^d,
\,
 j \in \{ 1, 2,..., m \}.
\end{align}
%\begin{equation} 
%\end{equation}
%When the method parameter $\eta \neq 0$, 
Additionally we assume that the vector functions 
$ \eta \mathcal{L}^{j} g_j \colon \R^d \rightarrow \R^d, \eta \in [0,1], j \in \{ 1, 2,..., m \}$ are 
%three times continuously differentiable 
continuously differentiable and
\begin{equation}
\label{eq:ass-Ljgj}
\eta \big \| \big ( \tfrac{ \partial (\mathcal{L}^{j} g_j ) } { \partial x} ( x) - \tfrac{ \partial (\mathcal{L}^{j} g_j ) } { \partial x} ( \tilde{x} ) \big ) y\big\| 
\leq
C ( 1 + \| x \|+ \| \tilde{x} \| )^{ \gamma - 2 } 
\| x - \tilde{x} \| \cdot \| y \|,
\quad
\forall  x, \tilde{x}, y \in \R^d.
\end{equation}
Moreover, the initial data $X_0 $ is supposed to be $ \mathcal{F}_0$-adapted, satisfying
\begin{equation}
\label{eq:ass-initial-data}
\| X_0 \|_{ L^{ p^* } ( \Omega; \R^d ) } < \infty.
\end{equation}
\end{assumption}
Recall that we use $ \tfrac{ \partial \phi } { \partial x} $ to denote the Jacobian matrix of a vector function
$\phi \colon \mathbb{R}^d \rightarrow \mathbb{R}^{d }$.
{\color{black}We mention that the condition \eqref{assum:eq-monotone-growth} is usually called a coercivity condition,
which is a classical one in the literature  to guarantee  that the exact solution has finite $p^*$-th moments, i.e.,
$\sup_{ t \in [ 0, T] }
\lVert X_{ t }  \rVert_{L^{p^*} ( \Omega; \R^d )} < \infty$.
The remaining conditions \eqref{eq:drift-polynomial-growth-condition}-\eqref{eq:ass-Ljgj} are a kind of polynomial growth conditions,
which have been also used in \cite{beyn2017stochastic,kumar2019milstein} to carry out the error analysis of Milstein type methods.
In subsection \ref{subsect:SPDE-example}, we present a system of SDEs that fulfill the above conditions.}
Note that the condition \eqref{eq:drift-polynomial-growth-condition} immediately implies
\begin{equation} \label{eq:fderivative-polynomial-growth}
\| \tfrac{ \partial f }{ \partial x} ( x ) y \| \leq
C ( 1 + \| x \| ) ^{\gamma-1} \| y \|,
\quad
\forall x, y \in \R^d,
\end{equation}
which in turn implies
\begin{align} 
\label{eq:f-polynomial-growth0}
\| f ( x ) - f ( \tilde{x} ) \|
& \leq
C
( 1 + \| x \| + \| \tilde{x} \| )^{\gamma-1} \| x - \tilde{x} \|,
\quad
\forall x, \tilde{x} \in \R^d,
\\
\label{eq:f-polynomial-growth}
\| f ( x ) \| 
& \leq
C ( 1 + \| x \| )^\gamma,
\quad
\forall x \in \R^d.
\end{align}
Likewise, the assumption \eqref{eq:diffusion-polynomial-growth-condition} ensures
%combining \eqref{eq:drift-polynomial-growth-condition} with \eqref{eq:mono-condition1}  offers
\begin{equation}
\label{eq:partial-g-growth}
\| \tfrac{ \partial g_j }{ \partial x} ( x ) y \|^2 \leq
C ( 1 + \| x \| ) ^{\gamma-1} \| y \|^2,
\quad
\forall x, y \in \R^d,
\,
 j \in \{ 1, 2,..., m \},
\end{equation}
and therefore
\begin{equation} \label{eq:multiplicative-g-difference}
\| g_j ( x ) - g_j ( \tilde{x} ) \|^2
\leq
C
( 1 + \| x \| + \| \tilde{x} \| )^{\gamma-1} \| x - \tilde{x} \|^2,
\quad
\forall x, \tilde{x} \in \R^d,
 j \in \{ 1, 2,..., m \}.
\end{equation}
This in turn gives
\begin{equation}
\| g ( x) \|^2 \leq
C ( 1 + \| x \| )^{ \gamma + 1 },
\quad
\forall x \in \R^d.
\end{equation}
Similarly as above, the assumption \eqref{eq:ass-Ljgj} promises
\begin{equation}
\eta
\big \| \tfrac{ \partial \mathcal{L}^{j} g_j } { \partial x} ( x)  y\big\| 
\leq
C ( 1 + \| x \| )^{ \gamma - 1 } 
\| y \|,
\quad
\forall  \, x,  y \in \R^d,
\end{equation}
and hence
\begin{equation}\label{eq:Ljgj-difference-growth}
\begin{split}
\eta \big \| \mathcal{L}^{j} g_j  ( x ) - \mathcal{L}^{j} g_j ( \tilde{x} ) \big \|
& \leq
C
( 1 + \| x \| + \| \tilde{x} \| )^{\gamma-1} \| x - \tilde{x} \|,
\quad
\forall \, x, \tilde{x} \in \R^d,
\\
\eta \big \| \mathcal{L}^{j} g_j  ( x )  \big \|
& \leq
C
( 1 + \| x \| )^{ \gamma },
\quad
\forall \, x \in \R^d.
\end{split}
\end{equation}
Further, Assumption \ref{ass:f-polynomial-growth} together with Assumption \ref{ass:monotonicity-condition} in $D = \R^d$
suffices to guarantee Assumption \ref{ass:well-possedness} holds in $D = \R^d$.
More formally, under these assumptions, the SDE \eqref{eq:SODE} possesses
a unique adapted solution with continuous sample paths, $X \colon [0,T] \times \Omega \rightarrow \mathbb{R}^d $,
satisfying 
\begin{equation}\label{eq:moment-bound-Xt}
\sup_{ t \in [ 0, T] }
\lVert X_{ t }  \rVert_{L^p ( \Omega; \R^d )}
\leq
C \big( 1 + \lVert X_{ 0 }  \rVert_{L^p ( \Omega; \R^d )} \big)
< \infty ,
\quad
p \in [2,  p^*  ],
%\:
%p^* \in [ 4 \gamma - 2, \infty),
\end{equation}
and thus
$\sup_{ s \in [ 0, T]} \E [ \| X_{ s } \|^2 ] + \sup_{ s \in [ 0, T]} \E [ \| f ( X_{ s } ) \|^2 ] + 
\sum_{j_1, j_2 = 1}^{m} \E [ \| \mathcal{L}^{j_1} g_{j_2} ( X_{ 0 } ) \|]^2 < \infty$,
where $ p^* \in [ 6 \gamma - 4, \infty) $ comes from Assumption \ref{ass:f-polynomial-growth}.
%obeying, for any $p \leq p^*$,
%\begin{equation}
%\E [ \|  X_t \|^p ] \leq ...
%\end{equation}
Further, the condition \eqref{eq:mono-condition2} in $\R^d$ from Assumption \ref{ass:monotonicity-condition} ensures that
the implicit Milstein type methods \eqref{eq:general-scheme} are well-defined in $\R^d$.
%
%\begin{assumption}[Initial data]
%\label{ass:initial-data}
%Suppose that the initial data $X_0 $ is  $ \mathcal{F}_0$-adapted and
%\begin{equation}
%\label{eq:ass-initial-data}
%\| X_0 \|_{ L^{ p^* } ( \Omega; \R^d ) } < \infty.
%\end{equation}
%%
%\end{assumption}
%
Thanks to Assumption \ref{ass:f-polynomial-growth} as well as the above implications, one can straightforwardly show
%\begin{lem}
%\label{lemma-Xt1-Xt2}
%and for $ t_1 , t_2 \in [ 0, T ]$,
\begin{equation}\label{eq:Xt1-Xt2}
\left \lVert X_{ t_1 } - X_{ t_2 } \right \rVert_{L^\delta ( \Omega; \R^d )}
\leq
C  \Big( 1 + \sup_{ t \in [0, T] } \| X_t \|^\gamma_{L^{\gamma \delta } ( \Omega; \R^d ) } \Big) 
\lvert t_1 - t_2  \rvert ^{ \frac 12 },
\quad
\forall \,
\delta \in [1, \tfrac { p^* } { \gamma }  ],
\,
 t_1 , t_2 \in [ 0, T ].
\end{equation}
\subsection{Analysis of the mean-square convergence rate}
We are now ready to give the main result of this section that reveals the optimal mean-square convergence rate
of the considered schemes under Assumptions \ref{ass:monotonicity-condition}, \ref{ass:f-polynomial-growth}.
\begin{thm}[Mean-square convergence rates of the schemes]
\label{thm: convergence-rate-of-STM}
Let coefficients of SDEs \eqref{eq:SODE} and method parameters of the schemes \eqref{eq:general-scheme} obey
Assumption \ref{ass:monotonicity-condition} in the whole space $D = \R^d$.
Let Assumption \ref{ass:f-polynomial-growth} be fulfilled and let the step-size $ h = \tfrac{T}{N} \in (0, \tfrac{ 1 }{ 2 \theta L_2 } ) $ with
$\theta \in [\tfrac12, 1], N \in \N$. Then SDEs \eqref{eq:SODE} and the schemes \eqref{eq:general-scheme}  
admit unique adapted solutions in $\R^d$, denoted by 
$ \{ X_{t} \}_{ t \in [ 0, T] } $ and $ \{ Y_n \}_{ 0 \leq n \leq N}$, respectively. Furthermore, there exists a constant $C>0$, 
independent of $N \in \N$, such that, for any $N \in \N$,
	\begin{equation} \label{eq:thm-multiplicative-converg-rate}
	\begin{split}
	\sup_{ 0 \leq n \leq N}  \lVert X_{t_n} - Y_n  \rVert_{ L^2 ( \Omega; \R^d ) }
	\leq &
	C 
	\Big ( 
	1 +  \| X_0 \|^{\max \{ 2 \gamma, 3 \gamma - 2\} }_{L^{ \max\{ 4 \gamma, 6\gamma - 4 \} }(\Omega; \R^d)}
         \Big)
	h.
	\end{split}
	\end{equation}
\end{thm}
{\it Proof of Theorem \ref{thm: convergence-rate-of-STM}.}
The above discussion reminds us that all conditions in  Assumptions \ref{ass:monotonicity-condition}, 
\ref{ass:well-possedness} hold in $D = \R^d$. 
Therefore, Theorem \ref{thm:upper-error-bound} is applicable here
and we only need to properly estimate two error terms $\E [ \lVert R_{ i } \rVert ^2 ]$ and 
$\E [ \lVert \E ( R_{ i } | \mathcal{F}_{t_{ i -1} } ) \rVert ^2] $
before arriving at the expected mean-square convergence rate.
Recalling the definition of $\{R_i\}_{0 \leq i \leq N}$ given by \eqref{eq:Error-Remainder-Defn} 
and using a triangle inequality yield
\begin{align}\label{eq:R_i+1}
%\begin{split}
%\E \big[ 
\lVert R_{ i } \rVert_{ L^2 ( \Omega; \R^d ) } 
%\big]
\leq &
\theta
\bigg \lVert \int_{ t_{i-1} }^{ t_{ i } } f( X_s ) - f( X_{ t_{ i } } ) \, \dd s \bigg \rVert_{ L^2 ( \Omega; \R^d ) }
\nonumber
\\ &
+
( 1 - \theta ) \bigg \lVert \int_{ t_{i - 1 } }^{ t_{ i } } f( X_s ) - f( X_{ t_{i-1} } ) \, \dd s \bigg \rVert_{ L^2 ( \Omega; \R^d ) }
\nonumber
\\ & +
\bigg \lVert \int_{ t_{i - 1 } }^{ t_{ i } } g( X_s ) - g( X_{ t_{ i - 1 } } ) \, \dd W_s  
- \sum^m_{j_1,j_2=1}
            \mathcal{L}^{j_1} g_{j_2}( X_{ t_{i-1} } ) I_{j_1,j_2}^{t_{i-1}, t_{i}} \bigg \rVert_{ L^2 ( \Omega; \R^d ) }
 \nonumber \\ &
            +
            \bigg \lVert
            \tfrac {\eta} {2} \sum_{ j = 1}^{ m} \mathcal{L}^{j } g_{j }( X_{ t_{i} } ) h
            -
             \tfrac {\eta} {2} \sum_{ j = 1}^{ m} \mathcal{L}^{j } g_{j }( X_{ t_{i-1} } ) h
            \bigg \rVert_{ L^2 ( \Omega; \R^d ) }
\nonumber \\
=: & \
\mathbb{I}_1 + \mathbb{I}_2 + \mathbb{I}_3 + \mathbb{I}_4.
%\end{split}
\end{align}
Next we handle the first term in (\ref{eq:R_i+1}) and the second term can be treated similarly.
Using the H\"{o}lder inequality, \eqref{eq:f-polynomial-growth0},  \eqref{eq:moment-bound-Xt} and \eqref{eq:Xt1-Xt2} shows
\begin{equation}\label{eq:estimate_I1}
\begin{split}
\mathbb{I}_1
& \leq
\theta
\int_{ t_{i-1} }^{ t_{ i } } \lVert  f( X_s ) - f( X_{ t_{ i } } )  \rVert_{ L^2 ( \Omega; \R^d ) } \, \dd s
\\ &  \leq 
C \int_{ t_{i-1} }^{ t_{ i } } \left( 1 +  \lVert  X_s   \rVert^{ \gamma - 1 }_{ L^{ 4 \gamma - 2 } ( \Omega; \R^d ) } 
+  \lVert X_{ t_{ i  } }  \rVert^{ \gamma - 1 }_{ L^{ 4 \gamma - 2 } ( \Omega; \R^d ) } \right) 
 \lVert X_s - X_{ t_{ i  } }  \rVert_{ L^{ ( 4 \gamma - 2 )/ \gamma } ( \Omega; \R^d ) } \, \dd s 
\\ & \leq
C \big( 1 + \sup_{t \in [0, T]} \| X_t \|^{2 \gamma - 1}_{ L^{  4 \gamma - 2  } ( \Omega; \R^d ) }  \big) h^\frac32.
\end{split}
\end{equation}
In the same way, one can also obtain
\begin{equation}
\mathbb{I}_2
\leq
C \big( 1 + \sup_{t \in [0, T]} \| X_t \|^{2 \gamma - 1}_{ L^{  4 \gamma - 2  } ( \Omega; \R^d ) }  \big) h^\frac32.
\end{equation}
Before coming to the estimate of $\mathbb{I}_3$, we note that, for any differentiable functions $\phi \colon \R^d \rightarrow \R^d$,
\begin{equation} \label{eq:phi-expansion}
\begin{split}
\phi ( X_ t ) - \phi ( X_{ s })
& =
\tfrac{ \partial \phi }{ \partial x } ( X_{ s } ) ( X_t - X_{ s } ) 
+ \mathcal{R}_{\phi} ( X_{ s }, X_t  )
\\
&
=
\tfrac{ \partial \phi }{ \partial x } ( X_{ s } ) 
\Big(
\int_{ s }^ t f ( X_\xi ) \, \dd \xi 
+
\int_{ s }^t g ( X_\xi ) \, \dd W_\xi 
\Big)
+ \mathcal{R}_{\phi } ( X_{ s }, X_t  ),
\quad
s < t,
\end{split}
\end{equation}
where for short we denote
\begin{equation}
\mathcal{R}_{\phi} ( X_{ s }, X_t  ) 
:=
\int_0^1 
\big[
\tfrac{ \partial \phi }{ \partial x }
\big( X_{ s } + r ( X_t - X_{ s } )  \big)
-
\tfrac{ \partial \phi }{ \partial x } ( X_{ s } )
\big]
( X_t - X_{ s } )
\,
\dd r.
\end{equation}
As a direct consequence of \eqref{eq:derivative-operator} and \eqref{eq:phi-expansion},
one can show
\begin{equation}
\begin{split}
& 
\int_{ t_{i - 1 } }^{ t_{ i } } g( X_s ) - g( X_{ t_{ i - 1 } } ) \, \dd W_s  
- \sum^m_{j_1,j_2=1}
            \mathcal{L}^{j_1} g_{j_2}( X_{ t_{i-1} } ) I_{j_1,j_2}^{t_{i-1},t_{i}}
\\
& \quad =
\sum_{ j_2=1}^m
\int_{ t _{i-1} }^{ t_i }
\Big[
g_{j_2} ( X_s) - g_{j_2} ( X_{t_{i-1} }) - \sum_{ j_1 =1}^m \mathcal{L}^{j_1} g_{j_2} ( X_{t_{i-1} } ) ( W_s^{j_1} - W_{t_{i-1}}^{j_1} )
\Big]
\dd W_s^{j_2}
\\
& \quad =
\sum_{ j_2=1}^m
\int_{ t _{i-1} }^{ t_i }
\Big[
g_{j_2} ( X_s) - g_{j_2} ( X_{t_{i-1} }) - \sum_{ j_1 =1}^m \tfrac{ \partial g_{j_2} }{ \partial x } ( X_{t_{i-1} } ) g_{j_1} ( X_{t_{i-1} } ) 
( W_s^{j_1} - W_{t_{i-1}}^{j_1} )
\Big]
\dd W_s^{j_2}
\\
& \quad
=
\sum_{j=1}^m
\int_{ t _{i-1} }^{ t_i }
\Big[
\tfrac{\partial g_j}{\partial x} ( X_{t_{i-1}} ) 
\Big ( \! 
\int_{t_{i-1}}^s \! f ( X_\xi ) \, \dd \xi
+
%\tfrac{\partial g_j}{\partial x} ( X_{t_{i-1}} ) 
\int_{t_{i-1}}^s 
\!
\big [ g ( X_\xi ) - g ( X_{t_{i-1}} ) \big]  \dd W_ \xi
\Big )
\\ 
& \qquad \qquad \qquad
+
\mathcal{R}_{g_j} ( X_{t_{i-1} }, X_s  ) 
\Big]
\dd W_s^j.
\end{split}
\end{equation}
Bearing this in mind, one can utilize the It\^o isometry to obtain
\begin{equation}
\label{eq:estimate_I3}
\begin{split}
| \mathbb{I}_3 |^2
 & =
\sum_{j_2=1}^m
%\bigg (
\int_{ t _{i-1} }^{ t_i }
\E \Big[
\Big \|
g_{j_2} ( X_s) - g_{j_2} ( X_{t_{i-1} }) 
- 
\sum_{ j_1 =1}^m 
\mathcal{L}^{j_1}  g_{j_2} ( X_{t_{i-1} } ) ( W_s^{j_1} - W_{t_{i-1}}^{j_1} )
\Big \|^2
\Big]
\dd s
%\bigg )^\frac12
\\
& =
\sum_{j=1}^m
%\bigg (
\int_{ t _{i-1} }^{ t_i }
\E \Big[
\Big \|
\tfrac{\partial g_j}{\partial x} ( X_{t_{i-1}} ) 
\Big ( \! 
\int_{t_{i-1}}^s \! f ( X_\xi ) \, \dd \xi
+
%\tfrac{\partial g_k}{\partial x} ( X_{t_{i-1}} ) 
\int_{t_{i-1}}^s 
\!
\big [ g ( X_\xi ) - g ( X_{t_{i-1}} ) \big]  \dd W_ \xi
\Big )
+
\mathcal{R}_{g_j} ( X_{t_{i-1} }, X_s  ) 
\Big \|^2
\Big]
\dd s
%\bigg )^\frac12
\\
& \leq
3 \sum_{j=1}^m
%\bigg (
\int_{ t _{i-1} }^{ t_i }
\E \Big[
\Big \|
\tfrac{\partial g_j}{\partial x} ( X_{t_{i-1}} )  
\int_{t_{i-1}}^s \! f ( X_\xi ) \, \dd \xi
\Big \|^2
\Big]
\dd s
+
3
\sum_{j=1}^m
%\bigg (
\int_{ t _{i-1} }^{ t_i }
\E \big[
\big \|
\mathcal{R}_{g_j} ( X_{t_{i-1} }, X_s  )
\big \|^2
\big]
\dd s
\\
&
\quad +
3 \sum_{j=1}^m
%\bigg (
\int_{ t _{i-1} }^{ t_i }
\E \Big[
\Big \|
\tfrac{\partial g_j}{\partial x} ( X_{t_{i-1}} )  
\int_{t_{i-1}}^s \! g ( X_\xi ) - g ( X_{t_{i-1}} ) \, \dd W_\xi
\Big \|^2
\Big]
\dd s.
%\\
%& \leq
%C h^3
\end{split}
\end{equation}
In the following we cope with the above three items separately.
Thanks to \eqref{eq:f-polynomial-growth}, \eqref{eq:partial-g-growth} and  the H\"{o}lder inequality, we first get
\begin{equation} \label{eq:I3-estimate1}
\begin{split}
\Big \|
\tfrac{\partial g_j}{\partial x} ( X_{t_{i-1}} )  
\int_{t_{i-1}}^s \! f ( X_\xi ) \, \dd \xi
\Big \|_{L^2 ( \Omega; \R^d )}
& \leq
C
\int_{t_{i-1}}^s \!
\Big\|
( 1 + \| X_{t_{i-1}} \| )^{\frac{\gamma - 1}{2}}
\| f ( X_\xi ) \|
\Big \|_{L^2(\Omega; \R) }
\, \dd \xi
\\
& \leq
C
\int_{t_{i-1}}^s \!
\Big\|
( 1 + \| X_{t_{i-1}} \| )^{\frac{\gamma - 1}{2}}
( 1 +  \|X_\xi \| )^\gamma
\Big \|_{L^2(\Omega; \R) }
\, \dd \xi
\\
&
\leq
C h 
\Big ( 
1 + \sup_{ t \in [0, T]} \| X_t \|^{\frac{3 \gamma - 1}{2} }_{L^{3\gamma - 1}(\Omega; \R^d)}
\Big).
\end{split}
\end{equation}
Again, using the It\^o isometry, the H\"{o}lder inequality, \eqref{eq:partial-g-growth}, 
\eqref{eq:multiplicative-g-difference} and \eqref{eq:Xt1-Xt2} yields
\begin{equation} \label{eq:I3-estimate2}
\begin{split}
&
\E
\Big[
\Big \|
\tfrac{\partial g_j}{\partial x} ( X_{t_{i-1}} )  \!
\int_{t_{i-1}}^s \! g ( X_\xi ) - g ( X_{t_{i-1}} ) \, \dd W_\xi
\Big \|^2
%_{L^2 ( \Omega; \R^d )}
\Big]
 \\
&\quad 
=
 \sum_{ l = 1}^{m}
\int_{t_{i-1}}^s
\!
\big \|
\tfrac{\partial g_j}{\partial x} ( X_{t_{i-1}} )  
\big[ 
g_l ( X_\xi ) - g_l ( X_{t_{i-1}} )
\big]
\big \|^2_{L^2 ( \Omega; \R^d )}
\dd \xi
\\
& \quad \leq
C
 \sum_{ l = 1}^{m}
\int_{t_{i-1}}^s \!
\Big\|
( 1 + \| X_{t_{i-1}} \| )^{\frac{\gamma - 1}{2}}
\| g_l ( X_\xi ) -  g_l ( X_{t_{i-1}} ) \|
\Big \|^2_{L^2(\Omega; \R) }
\, \dd \xi
\\
& \quad \leq
C
\int_{t_{i-1}}^s \!
\Big\|
( 1 + \| X_{t_{i-1}} \| )^{\frac{\gamma - 1}{2}}
( 1 +  \|X_\xi \| + \|X_{t_{i-1}} \| )^{\frac{\gamma - 1}{2}} \| X_\xi - X_{t_{i-1}} \|
\Big \|^2_{L^2(\Omega; \R) }
\, \dd \xi
\\
& \quad \leq
C
\int_{t_{i-1}}^s \!
\Big\|
( 1 +  \|X_\xi \| + \|X_{t_{i-1}} \| )^{\gamma - 1} \| X_\xi - X_{t_{i-1}} \|
\Big \|^2_{L^2(\Omega; \R) }
\, \dd \xi
\\
&
\quad \leq
C h^2 
\Big ( 
1 + \sup_{ t \in [0, T]} \| X_t \|^{4 \gamma - 2}_{L^{4\gamma - 2}(\Omega; \R^d)}
\Big).
\end{split}
\end{equation}
In light of \eqref{eq:diffusion-polynomial-growth-condition}, \eqref{eq:moment-bound-Xt} and \eqref{eq:Xt1-Xt2}, 
one can further use the H\"{o}lder inequality to acquire
\begin{equation} \label{eq:I3-estimate3}
\begin{split}
&
\big \|
\mathcal{R}_{g_j} ( X_{t_{i-1} }, X_s  )
\big \|_{L^2 ( \Omega; \R^d )}
\\ & \quad
\leq
\int_0^1 
\big \|
\big[
\tfrac{ \partial g_j }{ \partial x }
\big( X_{t_{i-1} } + r ( X_s - X_{t_{i-1} } )  \big)
-
\tfrac{ \partial g_j }{ \partial x } ( X_{t_{i-1} } )
\big]
( X_s - X_{t_{i-1} } )
\big \|_{L^2 ( \Omega; \R^d )}
\dd r
\\
& \quad
\leq
C
\int_0^1 \!
\Big\|
( 1 + \| r X_s + ( 1 - r ) X_{t_{i-1} }  \| +  \| X_{t_{i-1}} \| )^{\frac{\gamma - 3}{2}}
 \| X_s - X_{t_{i-1} } \|^2
\Big \|_{L^2(\Omega; \R) }
\, \dd r
\\
& \quad
\leq
C h 
\Big ( 
1 + \sup_{ t \in [0, T]} \| X_t \|^{\frac{\max \{ 4 \gamma, 5 \gamma - 3 \} }{2} }_{L^{ \max\{ 4 \gamma, 5\gamma - 3 \} }(\Omega; \R^d)}
\Big).
\end{split}
\end{equation}
Plugging the above three estimates \eqref{eq:I3-estimate1}-\eqref{eq:I3-estimate3} into \eqref{eq:estimate_I3} gives
\begin{equation}
\mathbb{I}_3 
\leq
C h^\frac32 
\Big ( 
1 + \sup_{ t \in [0, T]} \| X_t \|^{\frac{\max \{ 4 \gamma, 5 \gamma - 3 \} }{2} }_{L^{ \max\{ 4 \gamma, 5\gamma - 3 \} }(\Omega; \R^d)}
\Big). 
\end{equation}
With regard to $\mathbb{I}_4$, we utilize \eqref{eq:Ljgj-difference-growth}-\eqref{eq:Xt1-Xt2} 
and  the H\"{o}lder inequality to obtain
%the term vanishes when $\eta = 0$. So we only consider $\eta \neq 0$:
\begin{equation}
\label{eq:estimate_I4}
\begin{split}
\mathbb{I}_4
& \leq
\tfrac {\eta h} {2} 
\sum_{ j = 1}^{ m}
\big \lVert
              \mathcal{L}^{j } g_{j }( X_{ t_{i} } )
            -
             \mathcal{L}^{j } g_{j }( X_{ t_{i-1} } )
            \big \rVert_{ L^2 ( \Omega; \R^d ) }
\\
& \leq
C h
\Big \lVert
( 1 + \| X_{ t_{i} } \| + \| X_{ t_{i-1} } \| )^{\gamma-1} \| X_{ t_{i} } - X_{ t_{i-1} } \|
\Big \rVert_{ L^2 ( \Omega; \R ) }
\\
& \leq
C h^{\frac32}
\Big ( 
1 + \sup_{ t \in [0, T]} \| X_t \|^{ 2 \gamma - 1 }_{L^{4\gamma - 2}(\Omega; \R^d)}
\Big).
\end{split}
\end{equation}
Putting all the above estimates together we derive from \eqref{eq:R_i+1} that
%$
\begin{equation}\label{eq:multiplicative-ms-local-error}
\begin{split}
%\E [ \lVert R_{ i  } \rVert ^2 ]
\lVert R_{ i } \rVert_{ L^2 ( \Omega; \R^d ) } 
\leq
C h^\frac32
\Big ( 
1 + \sup_{ t \in [0, T]} \| X_t \|^{\frac{\max \{ 4 \gamma, 5 \gamma - 3 \} }{2} }_{L^{ \max\{ 4 \gamma, 5\gamma - 3 \} }(\Omega; \R^d)}
\Big). 
\end{split}
\end{equation}
Noting that the stochastic integral vanishes under the conditional expectation, one can, similarly as in \eqref{eq:R_i+1},
infer that
\begin{equation}
\label{eq:E_E_R_i+1^2}
\begin{split}
\lVert \E ( R_{ i } \vert \mathcal{ F }_{ t_{ i - 1 } } ) \rVert_{ L^2 ( \Omega; \R^d ) } 
\leq &
 \theta \bigg \lVert \E \bigg(  \int_{ t_{i-1} }^{ t_{ i } }  f( X_s ) - f( X_{ t_{ i  } } ) \dd s \Big \vert \mathcal{ F }_{ t_{ i - 1} } \bigg) 
\bigg \rVert_{ L^2 ( \Omega; \R^d ) } 
\\ & + 
 (1 - \theta ) \bigg \lVert \E \bigg( \int_{ t_{i-1} }^{ t_{ i } } f( X_s ) - f( X_{ t_{ i - 1 } } ) \dd s \Big \vert \mathcal{ F }_{ t_{ i - 1 } } 
\bigg) \bigg \rVert_{ L^2 ( \Omega; \R^d ) } 
\\ & +
\tfrac {\eta h } {2} \sum_{ j = 1}^{ m}
\Big \lVert \E \Big(
\big[ 
\mathcal{L}^{j } g_{j }( X_{ t_{i} } )
            -
 \mathcal{L}^{j } g_{j }( X_{ t_{i-1} } )  
\big]
\Big \vert \mathcal{ F }_{ t_{ i - 1 } }
\Big) \Big \rVert_{ L^2 ( \Omega; \R^d ) } 
\\ =: &
\ \mathbb{I}_5 + \mathbb{I}_6 + \mathbb{I}_7.
\end{split}
\end{equation}
In order to estimate $\mathbb{I}_5$, we first note that
\begin{equation} \label{eq:conditional-expectation-zero}
\begin{split}
\E \bigg(  \int_{ t_{i-1} }^{ t_{ i } } 
\tfrac{ \partial f }{ \partial x } ( X_{s } ) 
\int_s^{t_{i}} \! g ( X_\xi ) \, \dd W_\xi  
 \, \dd s \Big \vert \mathcal{ F }_{ t_{ i - 1} } \bigg)
& =
\int_{ t_{i-1} }^{ t_{ i } } 
\E \bigg( 
\int_s^{t_{i}} \! 
\tfrac{ \partial f }{ \partial x } ( X_{s } )  g ( X_\xi ) \, \dd W_\xi  
 \Big \vert \mathcal{ F }_{ t_{ i - 1} } \bigg)
 \dd s
\\
& =
\int_{ t_{i-1} }^{ t_{ i } } 
\!
\E \bigg( 
\E \bigg( \!
\int_s^{t_{i}} \! 
\tfrac{ \partial f }{ \partial x } ( X_{s } )  g ( X_\xi ) \, \dd W_\xi  
\Big \vert \mathcal{ F }_{ s} \bigg)
 \bigg \vert \mathcal{ F }_{ t_{ i - 1} } \bigg)
 \dd s
\\
& =
0.
\end{split}
\end{equation}
Using this and \eqref{eq:phi-expansion} with $\phi = f$ ensures
\begin{equation}
\E \bigg(  \int_{ t_{i-1} }^{ t_{ i } } f( X_{ t_{ i  } }) - f( X_s )  \, \dd s \Big \vert \mathcal{ F }_{ t_{ i - 1} } \bigg)
=
\E \bigg( \!
\int_{ t_{i-1} }^{ t_{ i } }
\Big[
\tfrac{ \partial f }{ \partial x } ( X_{s } ) 
\int_s^{t_{i}} \! f ( X_\xi ) \, \dd \xi 
+ 
\mathcal{R}_{f} ( X_s, X_{t_{i} }  )
\Big]
\dd s
\Big \vert \mathcal{ F }_{ t_{ i - 1} } \bigg),
\end{equation}
and thus
\begin{equation} \label{eq:I5-estimate0}
\begin{split}
\mathbb{I}_5 
& =
 \theta \bigg \lVert
\E \bigg( \!
\int_{ t_{i-1} }^{ t_{ i } }
\Big[
\tfrac{ \partial f }{ \partial x } ( X_{s } ) 
\int_s^{t_{i}} \! f ( X_\xi ) \, \dd \xi 
+ 
\mathcal{R}_{f} ( X_s, X_{t_{i} }  )
\Big]
\dd s
\Big \vert \mathcal{ F }_{ t_{ i - 1} } \bigg)
\bigg \rVert_{ L^2 ( \Omega; \R^d ) } 
\\
&
\leq
 \theta \bigg \lVert
\int_{ t_{i-1} }^{ t_{ i } }
\Big[
\tfrac{ \partial f }{ \partial x } ( X_{s } ) 
\int_s^{t_{i}} \! f ( X_\xi ) \, \dd \xi 
+ 
\mathcal{R}_{f} ( X_s, X_{t_{i} }  )
\Big]
\dd s
\bigg \rVert_{ L^2 ( \Omega; \R^d ) } 
\\
&
\leq
 \theta
\int_{ t_{i-1} }^{ t_{ i } }
\int_s^{t_{i}} 
\big \lVert
\tfrac{ \partial f }{ \partial x } ( X_{s } ) f ( X_\xi ) 
\big \rVert_{ L^2 ( \Omega; \R^d ) } 
\, \dd \xi 
\,
\dd s
+
 \theta 
\int_{ t_{i-1} }^{ t_{ i } }
\big \lVert
\mathcal{R}_{f} ( X_s, X_{t_{i} }  )
\big \rVert_{ L^2 ( \Omega; \R^d ) } 
\dd s,
\end{split}
\end{equation}
where the Jensen inequality was used for the second step.
Here we employ \eqref{eq:fderivative-polynomial-growth}, \eqref{eq:f-polynomial-growth} and the H\"{o}lder inequality to show
\begin{equation} \label{eq:parial-f-f-estimate}
\begin{split}
\big \lVert
\tfrac{ \partial f }{ \partial x } ( X_{s } ) f ( X_\xi ) 
\big \rVert_{ L^2 ( \Omega; \R^d ) } 
& \leq
C
\big \lVert
\big(
1 + \| X_{s } \| 
\big)^{\gamma - 1 }
\big(
1 
+
\|
X_\xi 
\|
\big)^\gamma
\big \rVert_{ L^2 ( \Omega; \R ) }
\\
& \leq
C
\big ( 
1 + \sup_{ t \in [0, T]} \| X_t \|^{ 2 \gamma - 1 }_{L^{4\gamma - 2}(\Omega; \R^d)}
\big),
\end{split}
\end{equation}
and employ \eqref{eq:drift-polynomial-growth-condition}, \eqref{eq:Xt1-Xt2} and the H\"{o}lder inequality to arrive at
\begin{equation} \label{eq:Rf-estimate}
\begin{split}
\big \|
\mathcal{R}_{f} ( X_s, X_{t_{i} }  )
\big \|_{L^2 ( \Omega; \R^d )}
& \leq
\int_0^1 
\big \|
\big[
\tfrac{ \partial f }{ \partial x }
\big( X_s + r ( X_{t_{i} } - X_s )  \big)
-
\tfrac{ \partial f }{ \partial x } ( X_{ s } )
\big]
(  X_{t_{i} } - X_{ s })
\big \|_{L^2 ( \Omega; \R^d )}
\dd r
\\
& \leq
C
\int_0^1 \!
\Big\|
( 1 + \| r X_{t_i} + ( 1 - r ) X_s   \| +  \| X_s \| )^{ \gamma - 2}
 \| X_{t_{i} } - X_s \|^2
\Big \|_{L^2(\Omega; \R) }
\, \dd r
\\
&
\leq
C h 
\Big ( 
1 + \sup_{ t \in [0, T]} \| X_t \|^{\max \{ 2 \gamma, 3 \gamma - 2\} }_{L^{ \max\{ 4 \gamma, 6\gamma - 4 \} }(\Omega; \R^d)}
\Big).
\end{split}
\end{equation}
Inserting \eqref{eq:parial-f-f-estimate} and \eqref{eq:Rf-estimate} into \eqref{eq:I5-estimate0} implies
\begin{equation}
\mathbb{I}_5 
\leq
C h^2 
\Big ( 
1 + \sup_{ t \in [0, T]} \| X_t \|^{\max \{ 2 \gamma, 3 \gamma - 2\} }_{L^{ \max\{ 4 \gamma, 6\gamma - 4 \} }(\Omega; \R^d)}
\Big).
\end{equation}
The estimates of $\mathbb{I}_6 $ and $\mathbb{I}_7 $ are similar and one can also get
\begin{equation}
\mathbb{I}_6
+
\mathbb{I}_7 
\leq
C h^2 
\Big ( 
1 + \sup_{ t \in [0, T]} \| X_t \|^{\max \{ 2 \gamma, 3 \gamma - 2\} }_{L^{ \max\{ 4 \gamma, 6\gamma - 4 \} }(\Omega; \R^d)}
\Big).
\end{equation} 
Therefore, from \eqref{eq:E_E_R_i+1^2} it immediately follows that
\begin{equation}
\lVert \E ( R_{ i } \vert \mathcal{ F }_{ t_{ i - 1 } } ) \rVert_{ L^2 ( \Omega; \R^d ) } 
\leq
C h^2 
\Big ( 
1 + \sup_{ t \in [0, T]} \| X_t \|^{\max \{ 2 \gamma, 3 \gamma - 2\} }_{L^{ \max\{ 4 \gamma, 6\gamma - 4 \} }(\Omega; \R^d)}
\Big).
\end{equation}
In view of Theorem \ref{thm:upper-error-bound} and \eqref{eq:moment-bound-Xt}, 
we validate the desired assertion \eqref{eq:thm-multiplicative-converg-rate}.
\qed

As already mentioned at the end of Section \ref{sect:upper-error-bound}, 
Theorem \ref{thm:upper-error-bound} still holds when Assumption \ref{ass:monotonicity-condition}
is replaced by Assumption \ref{ass:cor-monotonicity-condition-theta-Milstein}  
or Assumption  \ref{ass:cor-monotonicity-condition-backward-Milstein}. 
Therefore, the following two corollaries follow directly from
Corollaries \ref{cor:error-bounds-ass1}, \ref{cor:backward-Milstein-upper-error-bound}.
\begin{cor} \label{cor:MS-rate-Milstein1}
Let Assumption \ref{ass:cor-monotonicity-condition-theta-Milstein} be fulfilled with $D = \R^d$ 
and let $ \theta L_3 h \leq \nu $ for some $\theta \in [\tfrac12, 1]$, $\nu \in (0, 1)$.
Let conditions in Assumption \ref{ass:f-polynomial-growth} be all satisfied. 
Then SDEs \eqref{eq:SODE} and the semi-implicit Milstein methods \eqref{eq:semi-implicit-scheme}
admit unique adapted solutions in $ \R^d$ and \eqref{eq:thm-multiplicative-converg-rate} holds, namely,
the schemes \eqref{eq:semi-implicit-scheme}  retain a mean-square convergence rate of order one.
\end{cor}
%%%%%%
\begin{cor} \label{cor:MS-rate-Milstein2}
Let Assumption \ref{ass:cor-monotonicity-condition-backward-Milstein} be fulfilled with $D = \R^d$
and let conditions in Assumption \ref{ass:f-polynomial-growth} be all satisfied. 
Let $ \theta L_6 h \leq \nu $ for some $\theta \in (\tfrac12, 1]$, $\eta \in [0,1]$ and for some $\nu \in (0, 1)$.
Then SDEs \eqref{eq:SODE} and the proposed schemes \eqref{eq:general-scheme}  
admit unique adapted solutions in $ \R^d$ and \eqref{eq:thm-multiplicative-converg-rate} holds, namely,
the schemes \eqref{eq:general-scheme}  retain a mean-square convergence rate of order one.
\end{cor}
%%%%%%%%%%%%%
%\begin{example}
%...
%\end{example}
%%%%%%%%%%%%%
%
%%%
{\color{black}
\subsection{An example with numerical simulations}
\label{subsect:SPDE-example}
In this subsection, we aim to give an example SDE that satisfies Assumptions \ref{ass:monotonicity-condition}, \ref{ass:f-polynomial-growth}.
To this end, let us first consider the following semi-linear stochastic partial differential equation (SPDE) \cite{majee2018optimal,liu2021strong}:
\begin{equation}
\label{eq:SPDE}
\begin{split}
\left\{
    \begin{array}{ll}
 \dd  u (t, x )  = [ \tfrac{\partial^2 }{ \partial x^2} u (t, x) + u (t, x) - u^3 (t, x) ] \, \dd t + g ( u(t, x) ) \, \dd W_t,
\quad
t\in (0,T],  \: x \in (0, 1),
\\
 u (t, 0) = u (t, 1) = 0,
\\
 u ( 0, x )  = u_0 (x),
 \end{array}\right.
 \end{split}
\end{equation}
where $g \colon \mathbb{R} \rightarrow \mathbb{R} $ and  $W_{\cdot} \colon [0,T] \times \Omega \rightarrow \mathbb{R} $ 
is the real-valued standard Brownian motions.
Such an SPDE is usually termed as the stochastic Allen-Cahn equation.
Next we want to spatially discretize the above SPDE to obtain an SDE system.
On the interval $[0, 1]$ we construct a uniform mesh with stepsize $\Delta x : = \tfrac{1}{K}$ 
and denote $x_i = i \Delta x$, $i = 1, 2,..., K- 1$.
Discretizing the SPDE \eqref{eq:SPDE} spatially by a finite difference method yields a system of SDEs:
\begin{equation}
\label{eq:SPDE-SODE-system}
\dd  X_t = [\mathbb{A} X_t + \mathbb{F} ( X_ t ) ] \, \dd t + \mathbb{G} ( X_t ) \, \dd W_t,
\quad
t \in (0, T],
\quad
X_0 = x_0,
\end{equation}
where
$X_t 
=( X_{1,t},  X_{2,t}, \cdots, X_{K-1,t} ) ^T
:= (u(t, x_1),  u(t, x_2), \cdots, u(t, x_{K-1}) )^T
$,
$\mathbb{A} \in R^{ (K-1) \times (K-1) }$,
$x_0 = (u_0 ( x_1 ), u_0 ( x_2 ), ..., u_0 ( x_{K-1} ) )^T $
and
\begin{equation*}
%R^{ (K-1) \times (K-1) }
%\ni
\mathbb{A} 
= 
K^2 \left[\begin{array}{cccccc}
-2 & 1 & 0 & \cdots & 0 & 0   \\
1 & -2 & 1 & \cdots & 0 & 0 \\
0 & 1 & -2 & \cdots & 0 & 0 \\
 &  \cdots &  & \cdots &  &  \\
0 & 0 & 0 &  \cdots & -2 & 1 \\
0 & 0 & 0 &  \cdots & 1 & -2
\end{array}\right],
%_{K \times K},
%\text{tridiag}[1, -2, 1]
%\in \R^{ K \times K},
\:
\mathbb{F} ( X ) = \left[\begin{array}{c} X_1 - (X_1)^3 \\ X_2 - (X_2)^3  \\ \vdots \\ X_{K-1} - (X_{K-1})^3 \end{array}\right],
\:
\mathbb{G} ( X ) = \left[\begin{array}{c}  g(X_1) \\  g( X_2)   \\ \vdots \\  g ( X_{K-1} ) \end{array}\right].
\end{equation*}
%
%The reason why we are interested in \eqref{eq:SPDE-SODE-system} is that such an SDE system is obtained
%when using the finite difference method to spatially discretize the following stochastic partial differential equation
%with the spatial stepsize $\Delta x : = \tfrac{1}{m+1}$:
%
%
We do not consider the error caused by the spatial discretization but focus on the temporal discretization of 
the SDE system \eqref{eq:SPDE-SODE-system}, done by the  semi-implicit Milstein method ($\theta = 1, \eta = 0$).
%To fulfill conditions in Assumption \ref{ass:f-polynomial-growth}, 
Moreover, we assume $g \in C_b^3 (\mathbb{R}, \mathbb{R}) $, i.e., $g$ is three times differentiable with derivatives bounded.  
%
%\begin{equation}
%\begin{split}
%\tfrac{ 13 }{2} |  g ( u ) | ^ 2   & \leq c ( 1 +  | u |^ 2 ) + | u |^4,
%\\
%%| g'' (x) | & \leq c,
%| g' ( u ) - g ' ( v ) | & \leq c | u - v | ,
%\\
%%|g''(x)g (x) + ( g'(x) )^2 -  g''(y)g (y) - ( g'(y) )^2 | & \leq c ( 1 + | x | + | y | ) | x - y|.
%| ( g'( u )g ( u ) )' -  ( g'( v ) g ( v ) )' | & \leq c ( 1 + | u | + | v | ) | u -  v |.
%\end{split}
%\end{equation}
%$u_0 (x) \equiv 1 $, $f (u) = u - u^3$ and $g ( u ) = 2 u$.
%
It is easy to check all conditions in Assumption \ref{ass:f-polynomial-growth} are fulfilled in $ D = \mathbb{R}^{K-1}$ 
with $\gamma = 3$ and for any $p^* \geq 14$.
By setting $\theta = 1, \eta = 0$, conditions \eqref{eq:mono-condition1}, \eqref{eq:mono-condition2} 
in Assumption \ref{ass:monotonicity-condition} are also both satisfied in $D = \mathbb{R}^{K-1}$. Therefore, 
Theorem \ref{thm: convergence-rate-of-STM} is applicable, with the first convergence rate obtained for the semi-implicit Milstein method.
Since the SDE system has commutative noise, the Milstein type methods do not involve the Levy area
\cite{kloeden1992numerical,milstein2013stochastic} and can be implemented as easily as the Euler type methods. 

In what follows we set $g (u) = \sin (u) + 1$ and $u_0 ( x ) \equiv 1$ and do some numerical experiments.
In Figure \ref{fig:convergence-rates-SPDE-SODE-system}, we plot mean-square errors of 
the  semi-implicit Milstein method ($\theta = 1, \eta = 0$) for the SDE system \eqref{eq:SPDE-SODE-system} with $K = 4$.
There one can observe a convergence rate of order one, as the step-sizes shrink.
Here and below numerical approximations are performed using six different stepsizes $h = 2^{-i}, i =2,3,...7$.
The ``exact'' solution is identified as the numerical one using a fine stepsize $h_{\text{exact}} = 2^{-12}$
and the expectations are approximated by computing averages over $10^4$ samples.
%
%is numerically solved by  and the tamed Milstein method \cite{kumar2019milstein,wang2013tamed},
For comparison, we also discretize  \eqref{eq:SPDE-SODE-system} by the tamed Milstein method 
for  non-Lipschitz SDEs \cite{kumar2019milstein,wang2013tamed}.
Tables 1-3 provide mean-square approximation errors of  these two methods for three cases $K = 4, 8, 16$.
Clearly,  the tamed Milstein method gives satisfactory results
in the low dimension case $K = 4$ when the time stepsize is small,
i.e., $h = 2^{-5}, 2^{-6}, 2^{-7}$. As the dimension $K$ increases ($K = 8, 16$),  the tamed Milstein method
gives large errors and the approximations become unreliable for even small stepsizes.
However,   the semi-implicit Milstein method  performs much better,
%with stepsizes $h = 2^{-4}, 2^{-5}, 2^{-6}$ always  produce reliable results with small errors, 
even in high dimension case $K = 16$.
This happens because the eigenvalues $\{ \lambda_i \}_{i = 1}^{K-1}$ of $A$ are $\lambda_i = - 4 K^2 \sin^2 (\tfrac{i \pi}{2 K}) < 0$
and  the problem \eqref{eq:SPDE-SODE-system} turns to be a very stiff system \cite{milstein2013stochastic} as $K$ increases. 
As a kind of explicit method, the tamed Milstein method applied to solve stiff system,
faces severe time step-size reduction due to the stability issue.
%Indeed, the eigenvalues $\{ \lambda_i \}_{i = 1}^{K-1}$ of $A$ are $\lambda_i = - 4 K^2 \sin^2 (\tfrac{i \pi}{2 K}) < 0$
%and  the problem \eqref{eq:SPDE-SODE-system} turns to be a very stiff system \cite{milstein2013stochastic} as $K$ increases. 
%
%As far as these two implicit methods are concerned,  the SSBM method has considerably better accuracy than the SSBE method.
On the contrary, the semi-implicit Milstein method  has excellent stability property and is well suited for such stiff system.
%This guarantees that the SSBM method is more efficient than the SSBE method for the current example.

%%%
\begin{figure}[htp]
\centering
      \includegraphics[width=4in,height=3.5in]  {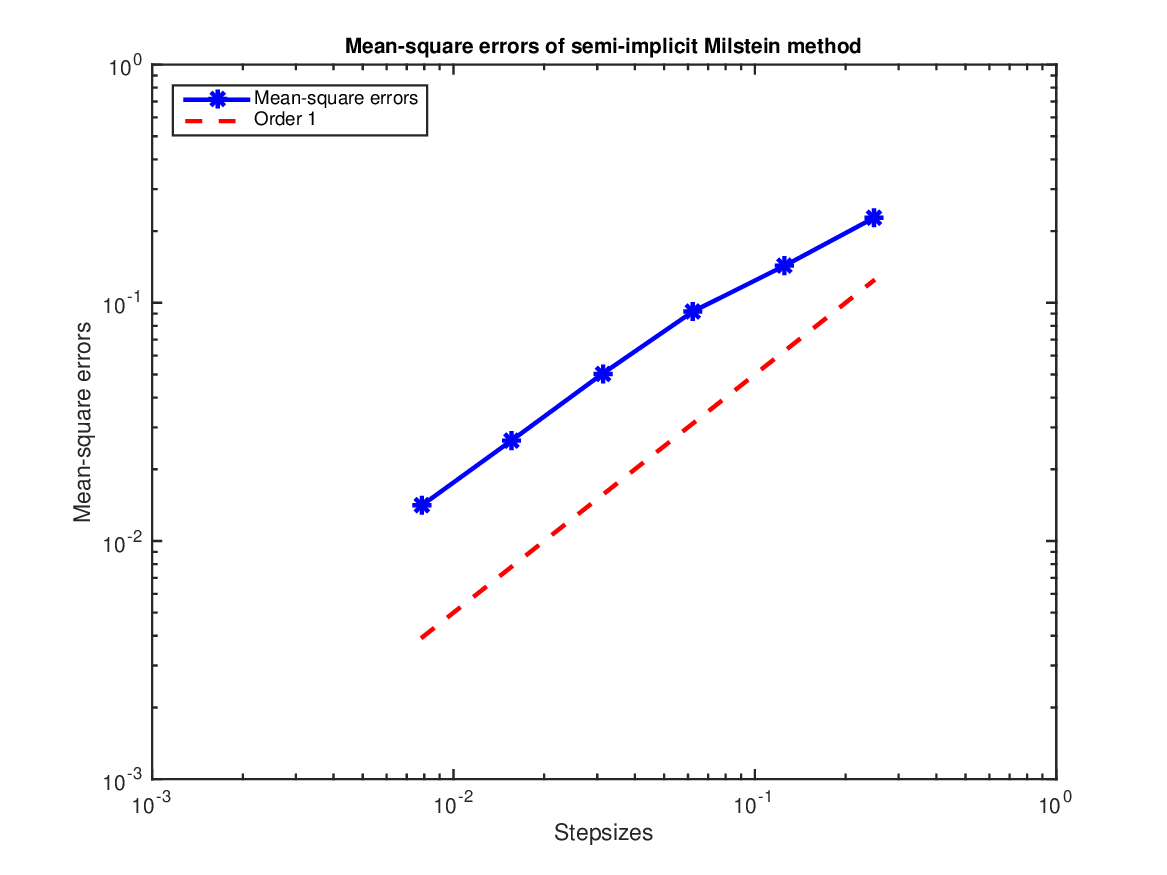}
 \caption{Mean-square convergence rate of the  semi-implicit Milstein method 
 for \eqref{eq:SPDE-SODE-system} ($K = 4$).}
\label{fig:convergence-rates-SPDE-SODE-system}
\end{figure}   
%%
%
%%%
\begin{table}[htp]
\label{table:m4SPDE}
\begin{center}% \footnotesize
\caption{Mean-square approximation errors for two schemes ($K = 4 $).}
\begin{tabular*}{12cm}{@{\extracolsep{\fill}}ccc}
\hline Stepsizes $h$ &  Semi-implicit Milstein   &  Tamed Milstein
\\ \hline
$h = 2^{-2}$ & 0.228228472003678   & 1.334521881473836
 \\ %\hline
$h = 2^{-3}$ &0.142671496841737   & 0.669681337348534
\\ %\hline
$h = 2^{-4}$ & 0.092138829109993   & 0.304845858687944   \\
 $h = 2^{-5}$ & 0.050402455908956   & 0.104293855220492
 \\ %\hline
$h = 2^{-6}$ & 0.026477850950294   & 0.044846913728710
 \\ %\hline
$h = 2^{-7}$ & 0.014040231850694   & 0.023917375308279
\\ %\hline
 \hline
\end{tabular*}
\end{center}
\end{table}
\begin{table}[htp]
\label{table:m8SPDE}
\begin{center}% \footnotesize
\caption{Mean-square approximation errors for two schemes ($K = 8 $).}
\begin{tabular*}{12cm}{@{\extracolsep{\fill}}ccc}
\hline Stepsizes $h$ &  Semi-implicit Milstein   &  Tamed Milstein
\\ \hline
$h = 2^{-2}$ & 0.337954132405219   & 3.127906338055271
 \\ %\hline
$h = 2^{-3}$ & 0.215927776446030   & 2.264234907688349
\\ %\hline
$h = 2^{-4}$ & 0.143858604122065   & 1.393606951102787   \\
 $h = 2^{-5}$ & 0.082829804151649   & 0.792573908322782
 \\ %\hline
$h = 2^{-6}$ & 0.045812280417151   & 0.375592176659368
 \\ %\hline
$h = 2^{-7}$ & 0.025766349691283   & 0.060932654697185
\\ %\hline
 \hline
\end{tabular*}
\end{center}
\end{table}
%%
%%%%%
\begin{table}[htp]
\label{table:m16SPDE}
\begin{center}% \footnotesize
\caption{Mean-square approximation errors for two schemes ($K = 16 $).}
\begin{tabular*}{12cm}{@{\extracolsep{\fill}}ccc}
\hline Stepsizes $h$ &  Semi-implicit Milstein   &  Tamed Milstein
\\ \hline
$h = 2^{-2}$ & 0.483493085665317   & 5.551900376818316
 \\ %\hline
$h = 2^{-3}$ & 0.310800712759207   & 4.923675828972162
\\ %\hline
$h = 2^{-4}$ & 0.209040389203629   & 4.088007760382119  \\
 $h = 2^{-5}$ & 0.122832739545349  & 3.113087309657318
 \\ %\hline
$h = 2^{-6}$ & 0.070290414827683   & 2.028019061710102
 \\ %\hline
$h = 2^{-7}$ & 0.041888961361398  & 1.128188804503420
\\ %\hline
 \hline
\end{tabular*}
\end{center}
\end{table}
%
%%%
%%%%%%
%\begin{table}\label{table:m32SPDE}
%\begin{center}% \footnotesize
%\caption{Mean-square errors for three schemes ($m = 32 $).}
%\begin{tabular*}{12cm}{@{\extracolsep{\fill}}cccc}
%\hline Stepsizes $h$ &  SSBM & SSBE  &  Tamed Milstein
%\\ \hline
%$h = 2^{-2}$ & 0.3048 & 0.3068  & 3.9822
% \\ %\hline
%$h = 2^{-3}$ & 0.1769 & 0.1874  & 2.4529
%\\ %\hline
%$h = 2^{-4}$ & 0.0801& 0.1108  & 1.6350   \\
% $h = 2^{-5}$ & 0.0417 & 0.0735  & 1.2847
% \\ %\hline
%$h = 2^{-6}$ & 0.0198 & 0.0469  & 0.7634
%\\ %\hline
% \hline
%\end{tabular*}
%\end{center}
%\end{table}
}
%%%
%%
%\section{Applications to stochastic volatility processes}
%\section{Applications to SDE models from finance}
%\section{Applications: convergence rates of positivity preserving schemes for financial SDEs with non-Lipschitz coefficients}
\section{Convergence rates of positivity preserving schemes for SDEs with non-globally Lipschitz coefficients}
\label{sect:applications}
In the present section, we turn our attention to the aforementioned scalar SDE models \eqref{eq:intro-32model} and 
\eqref{eq:intro-Ait-Sahalia-model-SDE} arising from mathematical finance.
Unlike general SDEs studied in the previous section, the considered financial models do not evolve in the whole space $\R$, 
but only in the positive domain $D = (0, \infty)$. This thus makes the convergence theory developed 
in the previous section not applicable in this situation.
Moreover, preservation of positivity is usually a desirable modeling property and
positivity of the approximation is,  in many cases, necessary in order for the numerical scheme
to be well defined (see, e.g., \eqref{eq:32model-scheme} and \eqref{eq:Milstein-Ait-Sahalia-model} below). 
However, numerical schemes are, in general, not able to preserve positivity. 
For example, the classical Euler-Maruyama method fails to preserve positivity for any scalar SDE \cite{kahl2008structure}.
In this section we choose two particular schemes from \eqref{eq:scheme-dm1} to 
approximate these two models, which are capable of preserving positivity of the continuous models.
By means of  the previously obtained error bound,  we carefully analyze the expected 
mean-square convergence rate of the resulting numerical approximations.
\subsection{The double implicit Milstein scheme for the Heston $ \frac{ 3 }{ 2 } $-volatility model}
\label{subsection:Milstein-32model}
As the first considered financial model, let us look at the Heston $\tfrac32$-volatility model \cite{heston1997simple,lewis2000option}:
\begin{equation}
\label{eq:32model}
d X_t = X_t ( \mu - \alpha X_t ) dt + \beta X_t ^ {3/2} \, \dd W(t), 
 \quad
 X_0 = x_0 > 0,
 \quad
 \mu , \alpha , \beta > 0,
 \,
  t > 0,
\end{equation}
which can be viewed as an inverse of a Cox-Ingersoll-Ross (CIR) process \cite{neuenkirch2014first}. Such an equation is
also used for modelling term structure dynamics \cite{ahn1999parametric}.
%can be viewed as a stochastic extension to the logistic equation \cite{Guard88}.
%@RELATED NUMERICAL RESULTS@
%%%%%%%%%%%%
Recently, some researchers \cite{chassagneux2016explicit,higham2013convergence,neuenkirch2014first} 
proposed and analyzed different positivity-preserving numerical schemes for strong approximations of the $ \frac{ 3 }{ 2 } $-process.
Similarly to \cite{higham2013convergence},  we choose a particular double implicit Milstein scheme 
\eqref{eq:scheme-dm1} with $ \theta = \eta = 1 $ to approximate the above $ \frac{ 3 }{ 2 } $-process. 
Furthermore, we attempt to prove the expected convergence rate for the scheme,
which is missing in \cite{higham2013convergence}.
%%%%%%%%%%%%
%here $ x ( t ) \in \R $ for each $  t > 0 $.
%@In general, numerical approximations do not preserve nonnegativity and hence 
%convergence theorems developed in [26, 27] cannot be applied in this situation. 
%Kahl et al. [17] proved that the EM scheme does
%not preserve positivity for any scalar SDE.@
%

Given $T \in (0, \infty)$ and $N \in \N$,  one can construct a uniform mesh on the interval $[0, T]$ with the uniform stepsize $ h = \tfrac{T}{N}$.
Based on the uniform mesh, we apply the drift-diffusion double implicit Milstein scheme \eqref{eq:scheme-dm1} with $ \theta = \eta = 1 $ 
to the model \eqref{eq:32model}, resulting in, for $n \in \{ 0, 1, 2, ..., N-1 \}$,
\begin{equation}
\label{eq:32model-scheme}
Y_{ n + 1 } =  Y_n + Y_{ n+ 1} ( \mu - \alpha Y_{ n+ 1} ) h 
+
\beta Y_n ^ {3/2} \Delta W_n + \tfrac{ 3 \beta^2 }{ 4 }  Y_{ n }^2 | \Delta W_n |^2 
-
\tfrac{ 3 \beta^2 }{ 4 }  Y_{ n+1 }^2 h 
,
\quad
Y_0 = X_0,
\end{equation}
which is a quadratic equation and has a unique positive solution explicitly given by $Y_0 = X_0$ and
\begin{equation}
\begin{split}
\label{eq:model32-Milstein-explicit}
Y_{n+1} & = 
\Big(
\sqrt{ 
( 1 - \mu h)^2 + 4  h ( \alpha + \tfrac34 \beta^2 )  ( Y_{n} + \beta |Y_n|^\frac32 \Delta W_n 
+ 
\tfrac34 \beta^2 | Y_n|^2 |\Delta W_n|^2 )
}
\\
& \qquad 
- 
( 1 - \mu h)   
\Big)
\Big/
(2 \alpha h + \tfrac32 \beta^2 h) > 0,
%\quad
%Y_0 = X_0,
\end{split}
\end{equation}
given that $Y_n >0$, $n \in \{ 0, 1, 2, ..., N-1 \}$.
We mention that no additional restriction is put on the stepsize $h > 0$ to ensure the positivity of the above approximations.
In order to carry out the error analysis for the scheme using Theorem \ref{thm:upper-error-bound}, 
we should first
%well pose the model and its approximations and to 
justify all conditions required 
in Assumptions \ref{ass:monotonicity-condition}, \ref{ass:well-possedness}, 
which are clarified in the forthcoming lemma.
\begin{lem}
\label{lem:32model-assumptions}
Let $\mu , \alpha , \beta > 0$, $X_0 >0$.
Then the Heston -$ \frac{ 3 }{ 2 } $ volatility model \eqref{eq:32model} has a unique global solution in $(0, \infty)$
and the scheme \eqref{eq:32model-scheme} produces unique positivity preserving approximations given by \eqref{eq:model32-Milstein-explicit}.
When  $\alpha > \tfrac32 \beta^2$,
%and $h \leq [ 2 \alpha - 3 (q - 1 ) \beta^2 ]/ [ \tfrac32 \beta^2 \mu +  2 \mu \alpha ] $ 
%for some $2 < q < 1 + \tfrac{ 2 \alpha }{ 3 \beta^2 } $, 
the SDE model \eqref{eq:32model} and the scheme \eqref{eq:32model-scheme} obey 
Assumptions \ref{ass:monotonicity-condition}, \ref{ass:well-possedness} in the domain $D =  (0, \infty)$
for some $2 < q < 1 + \tfrac{ 8 \alpha }{ 9 \beta^2 } $.
\end{lem}
{\it Proof of Lemma \ref{lem:32model-assumptions}.}
The well-posedness of the considered model \eqref{eq:32model} and the scheme \eqref{eq:32model-scheme} 
in the positive domain $(0, \infty)$ can be found in \cite{neuenkirch2014first,higham2013convergence}.
It remains to validate the other conditions in Assumptions \ref{ass:monotonicity-condition}, \ref{ass:well-possedness}.
For brevity,  we denote the drift and diffusion coefficients of SDE \eqref{eq:32model} by
\begin{equation}
f ( x ) :=  x ( \mu - \alpha x  ),
\quad
g ( x ) := \beta x^{ \frac32 },
\quad
x \in \R_+.
\end{equation}
As a result, $ g ' (x) g ( x ) = \tfrac32 \beta^2 x^2, \, x \in \R_+ $ and one can find a positive constant $\tilde{c} > 0$ such that
\begin{equation}
\begin{split}
\Xi (x, y, h)
& :=
 \varrho \tfrac{ h }{ 2 } \| g'g ( x ) - g'g ( y ) \| ^2
  + \eta h \langle  g'g ( x ) - g'g ( y ) , f ( x ) - f ( y ) \rangle
  - h \| f ( x ) - f ( y ) \|^2
\\
&  =
\tfrac98 \varrho \beta^4 h ( x^2 - y^2 )^2 + \tfrac32 \beta^2 \eta \mu h ( x^2 - y^2 ) ( x - y)
 - \tfrac32 \beta^2 \eta \alpha h ( x^2 - y^2 )^2
 \\
& \quad
 - h \big[ \mu^2 ( x - y)^2
 - 2 \mu \alpha ( x - y ) ( x^2 - y^2 ) + \alpha^2 ( x^2 - y^2 )^2 \big]
\\
&  = 
\big \{ [ \tfrac98 \varrho \beta^4 - \tfrac32 \beta^2 \alpha - \alpha^2 ] ( x + y )^2 
+ 
[ \tfrac32 \beta^2 \mu +  2 \mu \alpha ] ( x + y ) - \mu^2   \big \}  ( x - y)^2 h
\\
& \leq
%\big[ \tfrac32 \beta^2 \mu +  2 \mu \alpha \big] ( x + y )  ( x - y)^2 h,
\tilde{c} ( x - y)^2 h,
\qquad
\forall \, x , y \in \R_+,
\end{split}
\end{equation}
where we used the facts that $\eta = 1$ and that $ \tfrac98 \varrho \beta^4 - \tfrac32 \beta^2 \alpha - \alpha^2 < 0$ for some $\varrho > 1$ 
since  
$ 
\alpha > \tfrac32 \beta^2 
%> \tfrac{3(\sqrt{3} -1)}{4} \beta^2
$ by assumption.
Further, we take some $2 < q < 1 + \tfrac{ 8 \alpha }{ 9 \beta^2 } $ 
%and $h \leq [ 2 \alpha - 3 (q - 1 ) \beta^2 ]/ [ \tfrac32 \beta^2 \mu +  2 \mu \alpha ] $ 
%to promise $\tfrac32 \beta^2 \mu h +  2 \mu \alpha h - 2 \alpha + 3 (q - 1 ) \beta^2 \leq 0$. 
to promise $ \tfrac94 (q - 1 ) \beta^2 - 2 \alpha \leq 0$
and hence
\begin{equation}
\begin{split}
&
2 \langle x - y ,f ( x ) - f ( y ) \rangle
  +  ( q - 1 ) \|g ( x ) - g ( y )\| ^2
+
\Xi (x, y, h)
\\
& = 
2 \mu | x - y |^2 -  2\alpha (x^2 - y^2) ( x - y) 
+
(q - 1 ) \beta^2 ( x^\frac32 - y^\frac32 )^2
+
\Xi (x, y, h)
\\
& \leq
2\mu | x - y |^2 
+
% \big[ \tfrac32 \beta^2 \mu h +  2 \mu \alpha h 
% - 2 \alpha + 3 (q - 1 ) \beta^2 \big] ( x + y )  ( x - y)^2
\big[ 
\tfrac94 (q - 1 ) \beta^2  - 2 \alpha \big] ( x + y )  ( x - y)^2
+
\tilde{c} ( x - y)^2 h
\\
& \leq
(2 \mu + \tilde{c} T ) | x - y |^2,
\quad
 \forall \, x , y \in \R_+,
\end{split}
\end{equation}
which means the condition \eqref{eq:mono-condition1} in Assumption \ref{ass:monotonicity-condition} is fulfilled.
%
%The drift coefficient is $ f( x ) = x ( \mu - \alpha x  ), \quad f \in C^1 ( \R,\R ) $,
%and the diffusion coefficient is
%$ g( x ( t ) ) = \beta x ( t ) ^ {\tfrac{3}{2}},
%\quad g \in C^2 ( \R,\R ) $ .
%
Now we validate \eqref{eq:mono-condition2} as follows:
\begin{equation}
\begin{split}
 &
 \big \langle x - y , f ( x ) - f ( y ) 
 - \tfrac{ 1 }{ 2 } 
  [ g' g ( x )
  - g' g ( y )  ] \big \rangle
   \\
   &
   \quad =
   \mu | x - y |^2 -  \alpha (x^2 - y^2) ( x - y) 
   -
   \tfrac34 \beta^2 ( x^2 - y^2 ) ( x - y)
  \\
  & \quad
   \leq
   \mu | x - y |^2,
   \qquad
   \forall \, x , y \in \R_+.
\end{split}
\end{equation}
%The condition \eqref{eq:mono-condition2} is thus validated.
Next we note that for any $p^* \leq 4$
\begin{equation} \label{eq:moment-bounds-condition}
 \big \langle
 x, f ( x ) 
 \big \rangle
 +
 \tfrac{ p^* - 1 }{ 2 }
 \| g ( x ) \|^2
 =
 \mu x^2 - ( \alpha - \tfrac{ p^* - 1 }{ 2 } \beta^2 ) x^3
 \leq
 \mu x^2,
 \quad
 \forall x \in \R_+,
\end{equation}
where the assumption $ \alpha > \tfrac32 \beta^2 $ was again used.
This assures $ \sup_{ t \in [ 0, T] }
\lVert X_{ t }  \rVert_{L^4 ( \Omega; \R^d )} < \infty$ and thus
$\sup_{ s \in [ 0, T]} \E [ \|  X_{ s } \|^2 ] + \sup_{ s \in [ 0, T]} \E [ \| f ( X_{ s } ) \|^2 ] 
%+ \sup_{ s \in [ 0, T]} \E [ \| g' ( X_{ s } )  g ( X_{ s } ) \|]^2 
< \infty$, 
as required in Assumption \ref{ass:well-possedness}.
\qed

Now we are able to apply Theorem \ref{thm:upper-error-bound} to deduce the convergence rate of the numerical scheme.
\begin{thm}
\label{thm:Milstein-convergence-rate-32model}
Let $X_0 > 0$ and let $\mu , \alpha , \beta > 0$ satisfy $ \alpha \geq \tfrac52 \beta^2 $. 
Let $\{ X_{ t} \}_{ t \in [0, T]}$ and $ \{Y_n\}_{0\leq n\leq N} $ be uniquely given by
\eqref{eq:32model} and \eqref{eq:model32-Milstein-explicit}, respectively.
Let $ h \in (0, \tfrac{1}{2 \mu} ) $.
%and $h \leq [ 2 \alpha - 3 (q - 1 ) \beta^2 ]/ [ \tfrac32 \beta^2 \mu +  2 \mu \alpha ] $ 
%for some $2 < q < 1 + \tfrac{ 8 \alpha }{ 9 \beta^2 } $.
Then there exists a constant $C>0$, independent of $N \in \N$, such that
\begin{equation} \label{eq:thm-32model-converg-rate}
	\begin{split}
	\sup_{ 0 \leq n \leq N} \lVert Y_n - X_{t_n} \rVert_{ L^2 ( \Omega; \R ) }
	\leq &
	C h 
%	\big( 1 +  \| X_0 \|_{L^6 ( \Omega; \R) }^3 \big)
	.
	\end{split}
	\end{equation}
\end{thm}
%
%%%
{\it Proof of Theorem \ref{thm:Milstein-convergence-rate-32model}.}
As already clarified in the proof of Lemma \ref{lem:32model-assumptions},  
the considered model and the scheme obey Assumptions \ref{ass:monotonicity-condition}, \ref{ass:well-possedness} 
in the domain $D =  (0, \infty)$. 
Therefore, Theorem \ref{thm:upper-error-bound} is applicable here
and it remains to estimate two error terms $\E [ \lVert R_{ i } \rVert ^2 ]$ and 
$\E [ \lVert \E ( R_{ i } | \mathcal{F}_{t_{ i -1} } ) \rVert ^2] $, $ i \in \{ 1, 2, ..., N\}$.
before attaining the convergence rate.
First of all, we recall 
$
f ( x ) :=  x ( \mu - \alpha x  ),
\,
g ( x ) := \beta x^{ \frac32 },
\,
x \in \R_+.
$
Following the notation used in Theorem \ref{thm:upper-error-bound},  one can easily see
\begin{align}\label{eq:23model-Ri}
%\begin{split}
%\E \big[ 
\lVert R_{ i } \rVert_{ L^2 ( \Omega; \R ) } 
%\big]
\leq &
\bigg \lVert \int_{ t_{i-1} }^{ t_{ i } } f( X_s ) - f( X_{ t_{ i } } ) \, \dd s \bigg \rVert_{ L^2 ( \Omega; \R ) }
            +
             \tfrac {h} {2}
            \big \lVert
            g' g ( X_{ t_{i} } ) 
            -
            g'g ( X_{ t_{i-1} } ) 
            \big \rVert_{ L^2 ( \Omega; \R ) }
\nonumber
\\ & +
\bigg \lVert \int_{ t_{i - 1 } }^{ t_{ i } } 
\big[
g( X_s ) - g( X_{ t_{ i - 1 } } ) 
- 
g'g ( X_{ t_{i-1} } ) ( W_s - W_{ t_{i - 1 } } )
\big]
\, \dd W_s  
\bigg \rVert_{ L^2 ( \Omega; \R ) }
\nonumber \\
=: &
I_1 + I_2 + I_3.
%\end{split}
\end{align}
Applying the It\^o formula to the quadratic polynomial $f (x) = x ( \mu - \alpha x  ), x \in \R_+$ yields
\begin{equation}
\label{eq:f-Ito}
f( X_{ t_{ i } } ) - f( X_s )
=
\int_{s}^{t_i} 
\big[
f' ( X_r ) f ( X_r ) + \tfrac12 f '' ( X_r ) g^2 ( X_r )
\big] 
\, \dd r
+
\int_{s}^{t_i} 
f' ( X_r ) g ( X_r )
\, \dd W_r
\end{equation}
and thus
\begin{equation} \label{eq:32-model-I1-estimate}
\begin{split}
I_1 
&  
%=
%\Big \lVert \int_{ t_{i-1} }^{ t_{ i } } 
%\Big( 
%...
%\Big)
%\, \dd s  \Big \rVert_{ L^2 ( \Omega; \R ) }
%\\
%&
\leq
\int_{ t_{i-1} }^{ t_{ i } } 
\int_{s}^{t_i} 
\lVert
f' ( X_r ) f ( X_r ) + \tfrac12 f '' ( X_r ) g^2 ( X_r )
\rVert_{ L^2 ( \Omega; \R ) }
\, \dd r
\, \dd s
\\
&
\quad +
\int_{ t_{i-1} }^{ t_{ i } } 
\Big (
\int_{s}^{t_i}
\E
[
\|
f' ( X_r ) g ( X_r )
\|^2
]
\, \dd r
\Big )^{\frac12}
\dd s
\\
&
\leq
C h^\frac32 \big( 1 + \sup_{s \in [0, T]} \| X_s \|_{L^6 ( \Omega; \R) }^3 \big),
\end{split}
\end{equation}
where one used the It\^o isometry and computed that $f' (x) f (x) = x ( \mu - \alpha x ) ( \mu - 2 \alpha x )$, $ f '' (x) g^2 ( x ) = - 2 \alpha  \beta^2 x^3 $
and $ f ' ( x ) g ( x ) = \beta ( \mu - 2 \alpha x ) x^{\frac32} $, $x \in \R_+$.
%on the condition $ \alpha \geq \tfrac52 \beta^2 $.
Since $g' ( x )g (x) = \tfrac32 \beta^2 x^2, x \in \R_+$ is also a quadratic polynomial, one can repeat the same lines as above to arrive at
%, when $ \alpha \geq \tfrac52 \beta^2 $,
\begin{equation} \label{eq:I2-estimate}
\begin{split}
I_2 
\leq 
C h^\frac32 \Big( 1 + \sup_{s \in [0, T]} \| X_s \|_{L^6 ( \Omega; \R) }^3 \Big)
.
\end{split}
\end{equation}
Also, applying the It\^o formula applied to $g (x) = \beta x^{ \frac32 }$ and $g' ( x )g (x) = \tfrac32 \beta^2 x^2$,  $x \in \R_+$,
using the It\^o isometry and considering \eqref{eq:I2-estimate}, one can show
%, 
%when $ \alpha \geq \tfrac52 \beta^2 $,
%%
\begin{equation} \label{eq:32model-I3}
\begin{split}
| I_3 |^2
& =
%\bigg ( 
\int_{ t_{i - 1 } }^{ t_{ i } } 
\lVert
g( X_s ) - g( X_{ t_{ i - 1 } } ) 
- 
g'g ( X_{ t_{i-1} } ) ( W_s - W_{ t_{i - 1 } } ) \rVert_{ L^2 ( \Omega; \R ) }^2 \, \dd s  
%\bigg )^\frac12
\\
& =
\int_{ t_{i - 1 } }^{ t_{ i } } 
\Big \lVert
\int_{ t_{i - 1 } }^{ s }
[
g' ( X_r ) f ( X_r ) 
+ \tfrac12 g'' ( X_r ) g^2 ( X_r ) ] \, \dd r
\\ &
\quad
+
\int_{ t_{i - 1 } }^{ s }
[
g' g ( X_r ) - g'g ( X_{ t_{i-1} } )
] \, \dd W_r
\Big \rVert_{ L^2 ( \Omega; \R ) }^2
\dd s
\\
&
\leq
2 h
\int_{ t_{i - 1 } }^{ t_{ i } } 
\int_{ t_{i - 1 } }^{ s }
\E
[
\|
g' ( X_r ) f ( X_r ) 
+ \tfrac12 g'' ( X_r ) g^2 ( X_r ) 
\|^2
]
\, \dd r
\\ & \quad
+
2
\int_{ t_{i - 1 } }^{ t_{ i } } 
\int_{ t_{i - 1 } }^{ s }
\E
[
\|
g' g ( X_r ) - g'g ( X_{ t_{i-1} } )
\|^2
]
\, \dd r
\\ 
&
\leq
C h^3 \big( 1 + \sup_{s \in [0, T]} \| X_s \|_{L^6 ( \Omega; \R) }^6 \big),
%\\ &
%\leq
%C h^3.
\end{split}
\end{equation}
where we computed that $ g' (x) f (x) = \tfrac32 \beta x^\frac32 ( \mu - \alpha x  ) $, 
$ g'' (x) g^2 (x) = \tfrac34 \beta^3 x^\frac52 $, $x \in \R_+$. 
Gathering the above three estimates together, we derive from \eqref{eq:23model-Ri} that
\begin{equation}
\label{eq:Ri-mean-square}
\lVert R_{ i } \rVert_{ L^2 ( \Omega; \R ) } 
\leq
C h^\frac32 \big( 1 + \sup_{s \in [0, T]} \| X_s \|_{L^6 ( \Omega; \R) }^3 \big).
\end{equation}
At the moment it remains to bound $\lVert \E ( R_{ i } \vert \mathcal{ F }_{ t_{ i - 1 } } ) \rVert_{ L^2 ( \Omega; \R^d ) }$,
which, similarly to \eqref{eq:E_E_R_i+1^2}, can be decomposed into two terms by a triangle inequality:
\begin{equation} \label{eq:Ri-conditional-decomposition}
\begin{split}
\lVert \E ( R_{ i } \vert \mathcal{ F }_{ t_{ i - 1 } } ) \rVert_{ L^2 ( \Omega; \R ) } 
\leq &
\bigg \lVert 
%\E \bigg(  
\int_{ t_{i-1} }^{ t_{ i } } 
\E \big( 
[ f( X_s ) - f( X_{ t_{ i  } } ) ]
\vert \mathcal{ F }_{ t_{ i - 1} } 
\big)
\dd s 
%\bigg) 
\bigg \rVert_{ L^2 ( \Omega; \R ) } 
\\ & +
\tfrac { h } {2} 
\big \lVert \E \big(
[ 
g' g ( X_{ t_{i} } ) 
            -
g'g ( X_{ t_{i-1} } )   
]
\big \vert \mathcal{ F }_{ t_{ i - 1 } }
\big) \big \rVert_{ L^2 ( \Omega; \R ) }.
%\\ 
%\leq &
%C h^3.
\end{split}
\end{equation}
Keeping \eqref{eq:f-Ito} in mind, recalling that the It\^o integral vanishes under the conditional expectation 
(see  \eqref{eq:conditional-expectation-zero} for clarification) and utilizing the Jensen inequality,  we derive
\begin{equation} \label{eq:f-conditional-expect-32-model}
\begin{split}
&
\bigg \lVert 
%\E \bigg(  
\int_{ t_{i-1} }^{ t_{ i } } 
\E \big( 
[ f( X_s ) - f( X_{ t_{ i  } } ) ]
\vert \mathcal{ F }_{ t_{ i - 1} } 
\big)
\,
\dd s 
%\bigg) 
\bigg \rVert_{ L^2 ( \Omega; \R ) } 
\\ &
\quad =
\bigg \lVert 
%\E \bigg(  
\int_{ t_{i-1} }^{ t_{ i } } 
\int_{s}^{t_i} 
\E
\big[
\big(
f' ( X_r ) f ( X_r ) + \tfrac12 f '' ( X_r ) g^2 ( X_r )
\big)
\vert \mathcal{ F }_{ t_{ i - 1} } 
\big] 
\, \dd r
\, \dd s 
%\bigg) 
\bigg \rVert_{ L^2 ( \Omega; \R ) } 
\\ &
\quad \leq
\int_{ t_{i-1} }^{ t_{ i } } 
\int_{s}^{t_i} 
\lVert
f' ( X_r ) f ( X_r ) + \tfrac12 f '' ( X_r ) g^2 ( X_r )
\rVert_{ L^2 ( \Omega; \R ) }
\, \dd r
\, \dd s
\\ &
\quad \leq
C h^2 \big( 1 + \sup_{s \in [0, T]} \| X_s \|_{L^6 ( \Omega; \R) }^3 \big).
\end{split}
\end{equation}
Following the same arguments as before,  one can derive
\begin{equation}
\tfrac{h}{2}
\big \lVert \E \big(
[ 
g' g ( X_{ t_{i} } ) 
            -
g'g ( X_{ t_{i-1} } )   
]
\big \vert \mathcal{ F }_{ t_{ i - 1 } }
\big) \big \rVert_{ L^2 ( \Omega; \R ) }
\leq
C h^2 \big( 1 + \sup_{s \in [0, T]} \| X_s \|_{L^6 ( \Omega; \R) }^3 \big).
\end{equation}
Plugging these two estimates into \eqref{eq:Ri-conditional-decomposition} results in
\begin{equation}
\label{eq:Ri-conditional-32-model}
\lVert \E ( R_{ i } \vert \mathcal{ F }_{ t_{ i - 1 } } ) \rVert_{ L^2 ( \Omega; \R ) } 
\leq
C h^2 \big( 1 + \sup_{s \in [0, T]} \| X_s \|_{L^6 ( \Omega; \R) }^3 \big).
\end{equation}
%
%Since the estimates of $J_1, J_2$ are the same, we just show that of $J_2$.
%\begin{equation}
%J_2
%=
%...
%\end{equation}
%
%%%%
%Again, by use of the It\^o formula  it is easy to check...
Analogously to \eqref{eq:moment-bounds-condition}, 
the assumption $ \alpha \geq \tfrac52 \beta^2 $ ensures
\begin{equation}\label{eq:32model-L6-MB}
\sup_{ t \in [ 0, T] }
\lVert X_{ t }  \rVert_{L^6 ( \Omega; \R )} < \infty.
\end{equation}
Thanks to \eqref{eq:32model-L6-MB} and Theorem \ref{thm:upper-error-bound}, 
the assertion \eqref{eq:thm-32model-converg-rate} follows based on \eqref{eq:Ri-mean-square} 
and \eqref{eq:Ri-conditional-32-model}.
\qed
\begin{rem}
Recall that strong convergence of the implicit Milstein scheme \eqref{eq:32model-scheme} for the $\tfrac32$ process
was analyzed by Higham et al. \cite{higham2013convergence}, 
with no convergence rates recovered.
Later in \cite{neuenkirch2014first}, with the aid of the Lamperti transformation, 
Neuenkirch and Szpruch \cite{neuenkirch2014first}  proposed a Lamperti transformed backward Euler method 
for a class of scalar SDEs in a domain including  the $\tfrac32$ process as a special case. 
There a mean-square convergence rate of order $1$ was proved for 
the Lamperti-backward Euler method solving the $\tfrac32$ process when the model parameters obey 
$\tfrac{ \alpha} { \beta^2 } > 5 $  (see Propositions 3.2 from \cite{neuenkirch2014first}).
%for a kind of Lamperti-backward Euler method
%Neuenkirch and Szpruch proved a mean-square convergence rate of order $1$ 
%for a kind of Lamperti-backward Euler method solving the Heston-$\tfrac32$ process 
%when $\tfrac{ \alpha} { \beta^2 } > 5 $  (see Propositions 3.2 from \cite{neuenkirch2014first}). 
In this work we turn to the implicit Milstein scheme \eqref{eq:32model-scheme}, 
covered by \eqref{eq:general-scheme} and also studied in \cite{higham2013convergence}, 
and successfully prove a mean-square convergence rate of order $1$ for the scheme 
on the condition $\tfrac{ \alpha} { \beta^2 } > \tfrac52 $. This not only fills the gap left by \cite{higham2013convergence},
but also significantly relaxes the restriction put on the model parameters as required in \cite{neuenkirch2014first}. 
\end{rem}
%
%%%%%%
{\color{black}
\subsection{The double implicit Milstein scheme for the stochastic Lotka-Volterra competition model}
In this subsection, we consider the scalar stochastic Lotka-Volterra (LV) competitive model \cite{mao2021positivity}
\begin{equation}
\label{eq:LVmodel}
\dd X_t =  [ b X_t  - a X_t^2 ] \, \dd t + \sigma X_t  \, \dd W_t,
\quad
X_0 = x_0 > 0
\end{equation}
for a single species, where individuals within the species are competitive and  $b, a, \sigma$ are all positive numbers.
The well-posedness of the model in the positive domain $(0, \infty)$ is known in the paper \cite{mao2021positivity},
where a positivity-preserving scheme is proposed, but with no convergence rate revealed.
%\cite[Theorem 2.1 on p. 381]{mao2007stochastic}.
%
On the uniform mesh, we apply the double implicit Milstein scheme \eqref{eq:scheme-dm1} with $ \theta = \eta = 1 $ 
to numerically solve the model \eqref{eq:LVmodel} as follows:
\begin{equation}
\label{eq:LVmodel-scheme}
Y_{ n + 1 } =  Y_n +  ( b Y_{ n+ 1}- a Y_{ n+ 1}^2 ) h 
+
\sigma Y_n \Delta W_n 
+  
\tfrac12 \sigma^2 Y_n  | \Delta W_n |^2 
-
\tfrac12 \sigma^2 Y_{n+1} h  
,
\quad
Y_0 = X_0.
\end{equation}
Obviously, it is a quadratic equation and has a unique positive solution:
\begin{equation}
\label{eq:LVmodel-Milstein-explicit}
Y_{n+1} = \frac{ - ( 1 - b h + \frac{\sigma^2}{2} h )
+
\sqrt{ ( 1 - b h + \frac{\sigma^2}{2} h )^2 + 4 a h Y_n (  1 + \sigma \Delta W_n + \frac{\sigma^2}{2}  |\Delta W_n |^2 ) }  }{2 a h} > 0,
\end{equation}
given that $Y_n >0$, $n \in \{ 0, 1, 2, ..., N-1 \}$.
We highlight that no additional restriction is put on the stepsize $h > 0$ to ensure the positivity of the above approximations.
Also, one can easily verify that
Assumptions \ref{ass:monotonicity-condition}, \ref{ass:well-possedness} are both fulfilled in the domain $D =  (0, \infty)$.
\begin{lem} \label{lem:LVmodel-assumptions}
Let $ b, a, \sigma > 0$, $X_0 >0$.
Then the stochastic LV competitive model \eqref{eq:LVmodel} has a unique global solution in $(0, \infty)$ satisfying
\begin{equation}
\label{eq:lema-LVmodel-moment}
\E \Big[ \sup_{t \in [0, T]} |X_t|^p \Big]
< \infty,
\quad
\forall p \geq 2.
\end{equation}
Moreover, the scheme \eqref{eq:LVmodel-scheme} produces unique positivity preserving approximations given by \eqref{eq:LVmodel-Milstein-explicit}.
The SDE model \eqref{eq:LVmodel} and the scheme \eqref{eq:LVmodel-scheme} obey 
Assumptions \ref{ass:monotonicity-condition}, \ref{ass:well-possedness} in the domain $D =  (0, \infty)$.
\end{lem}
{\it Proof of Lemma \ref{lem:LVmodel-assumptions}.}
The well-posedness of the model in the positive domain $(0, \infty)$ is known in \cite{mao2021positivity}
and the moment bound \eqref{eq:lema-LVmodel-moment} comes from \cite[Lemma 2.2]{mao2021positivity}. 
As discussed above, the scheme \eqref{eq:LVmodel-scheme}  has a unique positive solution given by \eqref{eq:LVmodel-Milstein-explicit}.
Consequently, all conditions in Assumption \ref{ass:well-possedness} are satisfied with $D = (0, \infty)$.
Now it remains to validate conditions in Assumption \ref{ass:monotonicity-condition}.
For brevity,  we denote the drift and diffusion coefficients of SDE \eqref{eq:32model} by
\begin{equation}
f ( x ) :=  b x - a x^2,
\quad
g ( x ) := \sigma x,
\quad
x > 0.
\end{equation}
%As a result, $ g ' (x) g ( x ) = \sigma^2 x, \, x \in \R_+ $.
By setting $\theta = \eta = 1$, the conditions \eqref{eq:mono-condition1}, \eqref{eq:mono-condition2} 
in Assumption \ref{ass:monotonicity-condition} reduce to
 \begin{align}
 \label{eq:LV-mono-condition1}
% \begin{split}
 2 \langle x - y , & f ( x ) - f ( y ) \rangle
   +  ( q - 1 ) \|g ( x ) - g ( y )\| ^2
  + \tfrac{ \varrho }{ 2 } h
   \| 
   g' (x) g ( x ) - g' (y) g ( y ) \| ^2
 \\ \nonumber
 & \quad
  -
   h \| f ( x ) - f ( y ) \|^2
   \leq L _ 1 \| x - y \|^2 ,
\\
 \label{eq:LV-mono-condition2}
 &
 \big \langle x - y ,  [ f ( x ) - f ( y ) ] - \tfrac12 [ g' (x) g (x) -  g' (y) g (y) ]
 \big \rangle
   \leq L_2 \| x - y \|^2,
\:
\forall x, y \in (0, \infty)
.
% \end{split}
 \end{align} 
Note that the diffusion $g$ is a linear function and $ g ' (x) g ( x ) = \sigma^2 x, \, x \in \R_+ $ is also linear.
Further, note that $f' ( x ) = b - 2 a x \leq b, \, \forall x \in (0, \infty) $.
These facts ensure that conditions \eqref{eq:LV-mono-condition1}-\eqref{eq:LV-mono-condition2} are both satisfied,
which validates Assumption \ref{ass:monotonicity-condition}.
\qed
%%%%%%%%%%%
%

Thanks to Lemma \ref{lem:LVmodel-assumptions}
and similarly to the proof of Theorem \ref{thm:Milstein-convergence-rate-32model},
we are now able to apply Theorem \ref{thm:upper-error-bound} to deduce the convergence rate 
of the numerical scheme \eqref{eq:LVmodel-scheme}.
\begin{thm}
\label{thm:Milstein-convergence-rate-LVmodel}
%
%Let $X_0 \in (0, \infty)$ and let $b, a, \sigma$ be positive numbers. 
Let $ b, a, \sigma > 0$, $X_0 >0$.
Let $\{ X_{ t} \}_{ t \in [0, T]}$ and $ \{Y_n\}_{0\leq n\leq N} $ be uniquely given by
\eqref{eq:LVmodel} and \eqref{eq:LVmodel-scheme}, respectively.
For $ h > 0$ satisfying $ ( 2 b - \sigma^2 ) h < 1 $,
%and $h \leq [ 2 \alpha - 3 (q - 1 ) \beta^2 ]/ [ \tfrac32 \beta^2 \mu +  2 \mu \alpha ] $ 
%for some $2 < q < 1 + \tfrac{ 8 \alpha }{ 9 \beta^2 } $.
there exists a constant $C>0$, independent of $N \in \N$, such that
\begin{equation} \label{eq:thm-LVmodel-converg-rate}
	\begin{split}
	\sup_{ 0 \leq n \leq N} \lVert Y_n - X_{t_n} \rVert_{ L^2 ( \Omega; \R ) }
	\leq &
	C h 
%	\big( 1 +  \| X_0 \|_{L^6 ( \Omega; \R) }^3 \big)
	.
	\end{split}
	\end{equation}
\end{thm}
%
%\begin{rem}
%...
%\end{rem}
By estimating $\E [ \lVert R_{ i } \rVert ^2 ]$ and  $\E [ \lVert \E ( R_{ i } | \mathcal{F}_{t_{ i -1} } ) \rVert ^2] $, $ i \in \{ 1, 2, ..., N\}$,
the proof of Theorem \ref{thm:Milstein-convergence-rate-LVmodel} is  similar to that of 
Theorem \ref{thm:Milstein-convergence-rate-32model} and omitted here.
Different from Theorem \ref{thm:Milstein-convergence-rate-32model} for the Heston-$ \frac{ 3 }{ 2 } $ volatility model,  
no further restriction is put on the parameters of the model \eqref{eq:LVmodel} because the diffusion coefficient is linear
and all required conditions are satisfied for full parameters $ b, a, \sigma > 0$.
}
%%%%%%
\subsection{The semi-implicit Milstein scheme for the Ait-Sahalia-type interest rate model}

%In what follows we focus on mean-square approximations of 
The next SDE financial model that we aim to numerically investigate is
the generalized Ait-Sahalia-type interest rate model 
\cite{ait-sahalia1996testing}, described by
\begin{equation}
\label{eq:Ait-Sahalia-model-SDE}
\dd X_t = ( \alpha_{-1} X_t^{-1} - \alpha_0 + \alpha_1 X_t - \alpha_2 X_t^\kappa  ) \, \dd t + \sigma X_t^{\rho} \, \dd W_t,
\quad
t > 0,
\quad
X_0 = x_0 > 0,
\end{equation}
where $\alpha_{-1}, \alpha_0, \alpha_1, \alpha_2, \sigma > 0$ are positive constants and $\kappa > 1, \rho >1$.
%SDE \eqref{eq:Ait-Sahalia-model-SDE} is referred to as the generalized Ait-Sahalia model.
%Clearly, 
Compared with the previous financial model \eqref{eq:32model}, a complication in \eqref{eq:Ait-Sahalia-model-SDE} is due to
the drift containing a term $ \alpha_{-1} X_t^{-1} $ that does not behave well near the origin.
The well-posedness of the model \eqref{eq:Ait-Sahalia-model-SDE} has been already shown in 
\cite[Theorem 2.1]{szpruch2011numerical} and we repeat it as follows.
\begin{prop}\label{thm:wellposed-Ait-Sahalia}
Let $\alpha_{-1}, \alpha_0, \alpha_1, \alpha_2, \sigma > 0 $ be positive constants and $\kappa > 1, \rho >1$. 
Given any initial data $X_0 = x_0 > 0$, there exists a unique, positive global solution $ \{ X_t \}_{t \geq 0}$ 
to \eqref{eq:Ait-Sahalia-model-SDE}.
\end{prop}
Recently, such a model has been numerically studied by many authors
\cite{chassagneux2016explicit,neuenkirch2014first,szpruch2011numerical,wang2018mean}, 
with an emphasis on introducing and analyzing various positivity preserving strong approximation schemes
(see Remark \ref{rem:Ait-Sahalia-model} for more details).
Different from numerical schemes introduced in \cite{chassagneux2016explicit,neuenkirch2014first,szpruch2011numerical,wang2018mean},
we apply the newly proposed Milstein scheme to the model \eqref{eq:Ait-Sahalia-model-SDE} with $ \kappa + 1 \geq 2 \rho $, 
covering both the  standard regime $ \kappa + 1 > 2 \rho $ and the critical regime $ \kappa + 1 = 2 \rho $, 
and successfully recover the expected mean-square convergence rate, by use of the previously obtained error bounds.
Given a uniform mesh on the interval $[0, T]$ with the uniform stepsize $ h = \tfrac{T}{N}, N \in \N, T \in (0, \infty)$, 
we apply the proposed Milstein type scheme \eqref{eq:general-scheme} with $\theta = 1, \eta = 0$ 
(called the semi-implicit Milstein method) to the above model \eqref{eq:Ait-Sahalia-model-SDE}
%, after a slight modification, 
and obtain numerical approximations, given by $Y_0 = X_0$ and
\begin{equation}\label{eq:Milstein-Ait-Sahalia-model}
\begin{split}
Y_{n+1} & = Y_{ n  } + h [ \alpha_{-1} Y_{n+1}^{-1} - \alpha_0 + \alpha_1 Y_{n+1} - \alpha_2 Y_{n + 1}^\kappa  ] 
\\
& \quad 
+ 
\sigma  Y_{ n }^{ \rho }  \Delta W_{ n } 
+ 
\tfrac12 \rho \sigma^2 Y_{ n }^{ 2 \rho - 1 } ( | \Delta W_{ n } |^2 - h  )
,
%\quad
%Y_0 = X_0,
\quad
n \in \{ 0, 1, 2, ..., N-1 \}.
\end{split}
\end{equation}
The next lemma concerns the well-posedness of the scheme \eqref{eq:Milstein-Ait-Sahalia-model}, 
which can be easily checked based on the observation that the drift coefficient function satisfies a monotonicity condition
(consult \cite[Lemma 3.1]{szpruch2011numerical} and \eqref{eq:monocity-Ait-Sahalia-0}).
\begin{lem}
\label{lem:well-posed-Milstein-ASmodel}
Let conditions in Proposition \ref{thm:wellposed-Ait-Sahalia} are all satisfied. 
For $ h \in (0, \tfrac{1}{ \alpha_1}] $,  the semi-implicit Milstein scheme \eqref{eq:Milstein-Ait-Sahalia-model} is well-defined 
in the sense that it admits a unique solution, preserving positivity of the underlying model \eqref{eq:Ait-Sahalia-model-SDE}.
\end{lem}
%
%{\it Proof of Lemma \ref{lem:well-posed-Milstein-ASmodel}.}
%For $\theta \in (0, 1]$, we define the function $F \colon \R_+ \rightarrow \R$, by
%$ F ( x ) : =  x - \theta h [ \alpha_{-1} x^{-1} - \alpha_0 + \alpha_1 x - \alpha_2 x^\kappa  ]   $, $x \in \R_+$.
%Evidently, the function $F$ is continuous and coercive on $\R_+$, 
%{\color{red}{
%i.e., $\lim_{ x \rightarrow \infty } F ( x ) = + \infty$,
%$ \lim_{ x \rightarrow 0^+ } F ( x ) = - \infty $. 
%}}
%Moreover, the function $F$ is strictly monotone as $ F' ( x ) =  1 + \theta h [ \alpha_{-1} x^{-2} - \alpha_1 
%+  \kappa \alpha_2 x^{\kappa -1} ] > 1 - \theta h \alpha_1 \geq 0 $ for $x \in \R_+$, whenever $ h \leq \tfrac{1}{ \theta \alpha_1} $.
%This means for any $y \in \R$   there exists a unique  $x \in \R_+$, solving $ F ( x ) = y$.
%\qed
%

For simplicity of notation in the following analysis, we update the definitions of functions $f, g$ in subsection \ref{subsection:Milstein-32model}
and denote the coefficients of SDE \eqref{eq:Ait-Sahalia-model-SDE} by
\begin{equation}
\label{eq:Ait-Sahalia-f-g-defn}
f ( x ) :=  \alpha_{-1} x^{-1} - \alpha_0 + \alpha_1 x - \alpha_2 x^\kappa,
\quad
g ( x ) := \sigma x^{ \rho },
\quad
x \in \R_+.
\end{equation}
It is easy to check that 
%$ g' g ( x ) =  \rho \sigma^2 x^{ 2 \rho - 1 }$ and 
\begin{equation}
g' (x) g ( x ) =  \rho \sigma^2 x^{ 2 \rho - 1 },
\qquad
f ( x ) - f ( y ) 
=
\big( -  \tfrac{ \alpha_{-1} }{x y} + \alpha_1 - \alpha_2 \tfrac{ x^\kappa -  y^\kappa } {x -y} \big) ( x - y).
%\quad
%g' g ( x ) - g' g ( y ) = 
\end{equation}
%We just focus on the particular case $ \theta = 1, \eta = 0$,
%\begin{equation}
%\begin{split}
%& \varrho \tfrac{ h }{ 2 } \| g'g ( x ) - g'g ( y ) \| ^2
%  - h \| f ( x ) - f ( y ) \|^2
%\\
%& \quad =
%\Big(
%\varrho \tfrac{ h }{ 2 } \rho^2 \sigma^4  \big( \tfrac{ x^{2 \rho - 1} - y^{2 \rho - 1} }{ x - y}  \big)^2
%-
%h
%\big(
% \tfrac{ \alpha_{-1} }{x y} - \alpha_1 + \alpha_2 \tfrac{ x^\kappa -  y^\kappa } {x -y}
%\big)^2
%\Big)
%(x - y)^2
%\end{split}
%\end{equation}
%
%Before proceeding further, 
In addition, for $\iota \in [1, \infty) $ we introduce a function $ z_{\iota} \colon \R_+ \times \R_+  \rightarrow \R_+$ defined by
\begin{equation}
\label{eq:z-iota-defn}
z_{\iota} (x, y) := \tfrac{ x^\iota -  y^\iota } {x -y}, \quad x, y \in \R_+.
\end{equation}

In the following error analysis for the numerical approximations, we cope with the standard 
case $ \kappa + 1 > 2 \rho $ and the critical case $ \kappa + 1 = 2 \rho $ separately, 
since different cases own different model properties.
%the corresponding error analysis makes a difference for different cases owning different model properties.

\vspace{0.3cm}
%\subsection{The non-critical case $ \kappa + 1 > 2 \rho $}
{\bf 5.2.1 The standard case $ \kappa + 1 > 2 \rho $}
\vspace{0.3cm}

At first, we focus on the standard case $ \kappa + 1 > 2 \rho $ and recall a lemma concerning (inverse) moment bounds 
of the solution to \eqref{eq:Ait-Sahalia-model-SDE}, quoted from \cite[Lemma 2.1]{szpruch2011numerical}.
\begin{lem} \label{Lem:moment-bounds-general}
Let conditions in Proposition \ref{thm:wellposed-Ait-Sahalia} be fulfilled with  $ \kappa + 1 > 2 \rho $
and let $ \{ X_t \}_{t \geq 0}$  be the unique solution to \eqref{eq:Ait-Sahalia-model-SDE}. 
Then for any $ p \geq 2$ it holds that
\begin{equation}
\sup_{ t \in [ 0, \infty) } \E [ | X_t |^p ] < \infty, 
\quad
\sup_{ t \in [ 0, \infty) } \E [ | X_t |^{ - p } ] < \infty.
\end{equation}
\end{lem}
In order to achieve the mean-square convergence rate of the scheme by means of Theorem \ref{thm:upper-error-bound}, 
we need to check all conditions required in Assumptions \ref{ass:monotonicity-condition}, \ref{ass:well-possedness}, 
which are clarified in the forthcoming lemma.
\begin{lem}
\label{lem:AS-model-noncritical-assumptions}
Let conditions in Proposition \ref{thm:wellposed-Ait-Sahalia} be all fulfilled with  $ \kappa + 1 > 2 \rho $ and let $ h \in (0, \tfrac{1}{ \alpha_1}] $. 
Then the Ait-Sahalia model \eqref{eq:Ait-Sahalia-model-SDE} and the scheme \eqref{eq:Milstein-Ait-Sahalia-model} obey 
Assumptions \ref{ass:monotonicity-condition}, \ref{ass:well-possedness} in the domain $D =  (0, \infty)$.
\end{lem}
{\it Proof of Lemma \ref{lem:AS-model-noncritical-assumptions}.}
%
%In the non-critical case $ \kappa + 1 > 2 \rho $, things are much complicated.
%We write $z_1 : =  \tfrac{ x^{2 \rho - 1} - y^{2 \rho - 1} }{ x - y} $ and $ z_2 : = \tfrac{ x^\kappa -  y^\kappa } {x -y}$. 
Note first that the well-posedness of the model and the scheme in $D = (0, \infty)$ 
has been proven in Proposition \ref{thm:wellposed-Ait-Sahalia} 
and Lemma \ref{lem:well-posed-Milstein-ASmodel}. It remains to confirm the other conditions.
We first claim that, for any $c >0$ there exists $a_0 \in [0, \infty)$ such that $ z_{2\rho -1} \leq c z_{\kappa} +  a_0$,
where we recall that $z_\iota$ is defined by \eqref{eq:z-iota-defn}.
Clearly, $z_{2\rho -1} (x, y) > 0$ and $z_{\kappa} (x, y) >0$ for all $x, y >0$.
Without loss of generality, we assume $x >y > 0$. Since $  \kappa - 1 >  2 \rho - 2$, for any $c >0$
one can find $a_0 \in [0, \infty)$ such that $ \sup_{v > 0} \big( ( 2 \rho - 1 ) v ^{ 2 \rho - 2 } - c \kappa v ^{ \kappa - 1 } \big) \leq a_0$.
As a consequence,
\begin{equation}
\begin{split}
( x - y) ( z_{2\rho -1} - c z_{\kappa} ) 
& =
x^{2 \rho - 1} - y^{2 \rho - 1} - c (x^\kappa -  y^\kappa)
\\
& =
(x - y)
\int_0^1 
\big [
(2 \rho - 1) 
\big (  y + \xi ( x - y) \big)^{ 2 \rho - 2 } - c \kappa \big (  y + \xi ( x - y) \big)^{ \kappa - 1 }  
\big]
\,
\dd \xi
\\
& \leq
a_0 (x - y),
\quad
\forall \, x >y > 0.
\end{split}
\end{equation}
The claim is thus validated. So one can choose $c < \tfrac{\sqrt{2} \alpha_2 }{ \sqrt{\varrho} \rho \sigma^2 } $
for some $\varrho > 1$ such that 
$
\tfrac{ \varrho }{ 2 } \rho^2 \sigma^4 c^2 
-
\alpha_2^2  < 0
$
and thus
\begin{equation}  \label{eq:g'g-f-condition-Ait-Sah}
\begin{split}
& \varrho \tfrac{ h }{ 2 } \| g' ( x ) g ( x ) - g' ( y ) g ( y ) \| ^2
  - h \| f ( x ) - f ( y ) \|^2
\\
&
\quad
=
\Big(
h \tfrac{ \varrho }{ 2 } \rho^2 \sigma^4  \big( \tfrac{ x^{2 \rho - 1} - y^{2 \rho - 1} }{ x - y}  \big)^2
-
h
\big(
 \tfrac{ \alpha_{-1} }{x y} - \alpha_1 + \alpha_2 \tfrac{ x^\kappa -  y^\kappa } {x -y}
\big)^2
\Big)
(x - y)^2
\\
& \quad
=
h \big[
\tfrac{ \varrho }{ 2 } \rho^2 \sigma^4  
%\big( \tfrac{ x^{ \kappa } - y^{ \kappa } }{ x - y}  \big)^2
z_{2\rho -1}^2
-
\big(
 \tfrac{ \alpha_{-1} }{x y} - \alpha_1 + \alpha_2 z_{\kappa}
\big)^2
\big]
(x - y)^2
\\ 
& \quad
=
h
\big[
\tfrac{ \varrho }{ 2 } \rho^2 \sigma^4  z_{2\rho -1}^2 
- 
( \tfrac{ \alpha_{-1} }{x y} - \alpha_1 )^2
-
\alpha_2^2 z_{\kappa}^2
-
2 ( \tfrac{ \alpha_{-1} }{x y} - \alpha_1 ) \alpha_2 z_{\kappa}
\big]
( x - y)^2
\\
&
\quad
\leq
h \big[
\tfrac{ \varrho }{ 2 } \rho^2 \sigma^4  z_{2\rho -1}^2 
-
\alpha_2^2 z_{\kappa}^2
+
2 \alpha_1 \alpha_2  z_{\kappa}
\big]
( x - y)^2
\\
&
\quad
\leq
h
\big[
\tfrac{ \varrho }{ 2 } \rho^2 \sigma^4  ( c z_{\kappa} +  a_0 )^2 
-
\alpha_2^2 z_{\kappa}^2
+
2 \alpha_1 \alpha_2  z_{\kappa}
\big]
( x - y)^2
\\
&
\quad
=
h
\big[
( \tfrac{ \varrho }{ 2 } \rho^2 \sigma^4 c^2 -  \alpha_2^2 ) z_{\kappa}^2
+
(
\varrho \rho^2 \sigma^4  c a_0
+
2 \alpha_1 \alpha_2  
)
z_{\kappa}
+
\tfrac{ \varrho }{ 2 } \rho^2 \sigma^4  a_0^2
\big]
( x - y)^2
\\
&
\quad
\leq C h ( x - y )^2,
\quad
\forall \, x, y \in \R_+.
\end{split}
\end{equation}
Furthermore, one can readily compute that, for any $ \kappa + 1 > 2 \rho $ and for some $q >2$,
\begin{equation}
\label{eq:one-sided-non-critical}
\sup_{ x > 0}
\big (
f'(x) + \tfrac{q-1}{2} | g' (x) |^2 
\big)
= 
\sup_{ x > 0}
\Big (
- \alpha_{-1} x^{ - 2} + \alpha_1 -  \alpha_2 \kappa x^{ \kappa - 1 } 
+ \tfrac{ (q-1) \sigma^2 \rho^2 }{ 2 } x^{ 2 \rho - 2 } 
\Big)
=: L
<
\infty.
%,
%\quad
%\forall \, x \in \R_+.
\end{equation}
This implies that
\begin{equation}\label{eq:monocity-Ait-Sahalia-0}
\begin{split}
& \langle x - y, f ( x ) - f ( y ) \rangle + \tfrac {q-1} {2} \lVert g ( x ) - g ( y ) \rVert ^ 2 
\\
& \quad =
\int_0^1 f' ( y + \xi (x-y) ) \dd \xi \cdot ( x - y)^2 
+
\tfrac {q-1} {2}
\Big | \int_0^1 g' ( y + \xi (x-y) ) \dd \xi  \Big | ^ 2
\cdot ( x - y)^2
\\
& \quad \leq
\int_0^1 
\big[ f' ( y + \xi (x-y) ) + \tfrac {q-1} {2} |g' ( y + \xi (x-y) )|^2 \big] \, \dd \xi \cdot ( x - y)^2
\\
& \quad \leq
L ( x - y)^2,
\quad
\forall \, x, y \in \R_+.
\end{split}
\end{equation}
%where $L := \sup_{ x >0} \big( - \alpha_{-1} x^{ - 2} + \alpha_1 -  \alpha_2 \kappa x^{ \kappa - 1 } 
%+ \tfrac{ (q-1) \sigma^2 \rho^2 }{ 2 } x^{ 2 \rho - 2 } \big)$.
%
Gathering \eqref{eq:g'g-f-condition-Ait-Sah} and \eqref{eq:monocity-Ait-Sahalia-0} together, 
the condition \eqref{eq:mono-condition1} is hence justified in the domain $ D = ( 0, \infty)$ with $\theta = 1, \eta = 0$.
From \eqref{eq:monocity-Ait-Sahalia-0}, one can assert that \eqref{eq:mono-condition2} is satisfied 
in $ D = ( 0, \infty)$ with $\theta = 1, \eta = 0, L_2 = \alpha_1$. Thus all conditions in Assumption \ref{ass:monotonicity-condition} 
are fulfilled in the domain $D =  (0, \infty)$. Assumption \ref{ass:well-possedness} follows
by taking Proposition \ref{thm:wellposed-Ait-Sahalia}, Lemmas \ref{lem:well-posed-Milstein-ASmodel}, 
\ref{Lem:moment-bounds-general} into consideration.
\qed
%
%Moreover, the condition \eqref{eq:mono-condition2} is evidently satisfied and Assumption \ref{ass:fg-greater-c0-c1} is fulfilled 
%with $ c_0 = c_1 = 0 $, for the particular case $ \theta = 1, \eta = 0$.

At the moment, we are well prepared to carry out the error analysis for the numerical approximations 
with the help of Theorem \ref{thm:upper-error-bound}.
\begin{thm}\label{thm:Mistein-Ait-Sahalia-convergence-non-critcal}
Let $\{ X_{ t} \}_{ t \in [0, T]}$ and $ \{Y_n\}_{0\leq n\leq N} $ be solutions to 
\eqref{eq:Ait-Sahalia-model-SDE} and \eqref{eq:Milstein-Ait-Sahalia-model}, respectively.
Let $q \in (2, \infty)$, $\varrho \in (1, \infty)$, let $\alpha_{-1}, \alpha_0, \alpha_1, \alpha_2, \sigma > 0$,
let $\kappa > 1, \rho >1$ obey  $ \kappa + 1 > 2 \rho $,
%
%Let $ h \in (0,  \min \{\tfrac{1}{ \alpha_1}, \tfrac{1}{2 L} \} ) $,
%where $L := \sup_{ x >0} \big( - \alpha_{-1} x^{ - 2} + \alpha_1 -  \alpha_2 \kappa x^{ \kappa - 1 } 
%+ \tfrac{ (q-1) \sigma^2 \rho^2 }{ 2 } x^{ 2 \rho - 2 } \big)$. 
%
and let $ h \in (0,  \tfrac{1}{ 2  \alpha_1} ) $.
Then there exists a constant $C>0$, independent of $N \in \N$, such that
\begin{equation}
%\sup_{ N \in \N }
\sup_{ 1\leq n \leq N}
\| X_{ t_n } - Y_n \|_{ L^2 ( \Omega; \R ) } 
\leq
C 
%(
%@\| X_0 \|@
%)
h.
\end{equation}
\end{thm}
{\it Proof of Theorem \ref{thm:Mistein-Ait-Sahalia-convergence-non-critcal}.}
As implied by Lemma \ref{lem:AS-model-noncritical-assumptions}, 
all conditions in Assumptions \ref{ass:monotonicity-condition}, \ref{ass:well-possedness} are fulfilled in $D =  (0, \infty)$.
Based on Theorem \ref{thm:upper-error-bound},  one just needs to properly estimate $\lVert R_{ i } \rVert_{ L^2 ( \Omega; \R ) } $ 
and $\lVert \E ( R_{ i } \vert \mathcal{ F }_{ t_{ i - 1 } } ) \rVert_{ L^2 ( \Omega; \R^d ) }$, $ i \in \{1, 2, ..., N \}$.
Following the notation used in \eqref{eq:Error-Remainder-Defn} and 
\eqref{eq:Ait-Sahalia-f-g-defn},  we first split the estimate of $\lVert R_{ i } \rVert_{ L^2 ( \Omega; \R ) } $ as follows:
\begin{align}\label{eq:Ait-Sahal-model-Ri}
%\begin{split}
%\E \big[ 
\lVert R_{ i } \rVert_{ L^2 ( \Omega; \R ) } 
%\big]
\leq &
\bigg \lVert \int_{ t_{i-1} }^{ t_{ i } } f( X_s ) - f( X_{ t_{ i } } ) \, \dd s \bigg \rVert_{ L^2 ( \Omega; \R ) }
\nonumber
\\ & +
\bigg \lVert \int_{ t_{i - 1 } }^{ t_{ i } } g( X_s ) - g( X_{ t_{ i - 1 } } ) 
- 
g'g ( X_{ t_{i-1} } ) ( W_s - W_{ t_{i - 1 } } ) \, \dd W_s  
\bigg \rVert_{ L^2 ( \Omega; \R ) }
\nonumber \\
=: &
I_4 + I_5.
%\end{split}
\end{align}
Repeating the same arguments as used in \eqref{eq:32-model-I1-estimate},
we apply the It\^o formula to $f (x) = \alpha_{-1} x^{-1} - \alpha_0 + \alpha_1 x - \alpha_2 x^\kappa, x \in \R_+$
and use Lemma \ref{Lem:moment-bounds-general} to derive
\begin{equation} \label{eq:I4-estimate}
\begin{split}
I_4 &  \leq
\int_{ t_{i-1} }^{ t_{ i } } \lVert f( X_s ) - f( X_{ t_{ i } } ) \rVert_{ L^2 ( \Omega; \R ) } \, \dd s 
%\\
%&
%\leq
%C h^\frac32 \Big( 1 + \sup_{s \in [0, T]} \| X_s \|_{L^6 ( \Omega; \R) }^3 \Big)
\\
&
\leq
C h^\frac32 \Big( 1 + \sup_{s \in [0, T]} \| X_s \|_{L^{4 \kappa - 2} ( \Omega; \R) }^{ 2 \kappa - 1 } 
+ 
\sup_{s \in [0, T]} \| X_s^{-1} \|_{L^{6} ( \Omega; \R) }^{ 3 }
\Big)
\\
&
\leq
C h^\frac32.
\end{split}
\end{equation}
Similarly to \eqref{eq:32model-I3}, by means of the It\^o isometry and the It\^o formula 
applied to $g ( x ) = \sigma x^{ \rho }$ and $g' (x) g ( x ) =  \rho \sigma^2 x^{ 2 \rho - 1 } , x \in \R_+$ one can show
\begin{equation} \label{eq:I5-estimate}
\begin{split}
| I_5 |^2
&\leq
2 h
\int_{ t_{i - 1 } }^{ t_{ i } } 
\int_{ t_{i - 1 } }^{ s }
\E
[
\|
g' ( X_r ) f ( X_r ) 
+ \tfrac12 g'' ( X_r ) g^2 ( X_r ) 
\|^2
]
\, \dd r \, \dd s
\\ & \quad
+
2
\int_{ t_{i - 1 } }^{ t_{ i } } 
\int_{ t_{i - 1 } }^{ s }
\E
[
\|
g' g ( X_r ) - g'g ( X_{ t_{i-1} } )
\|^2
]
\, \dd r \, \dd s
\\
& \leq
C h^3,
%\\ 
%&
%\leq
%C h^3 \big( 1 + \sup_{s \in [0, T]} \| X_s \|_{L^? ( \Omega; \R) }^? \big)
%\\ &
%\leq
%C h^3
%\big( 1 + \sup_{s \in [0, T]} \| X_0 \|_{L^? ( \Omega; \R) }^? \big).
\end{split}
\end{equation}
where the (inverse) moment bounds in Lemma \ref{Lem:moment-bounds-general} were also used for the last step.
Inserting \eqref{eq:I4-estimate} and \eqref{eq:I5-estimate} into \eqref{eq:Ait-Sahal-model-Ri} implies
\begin{equation}
\lVert R_{ i } \rVert_{ L^2 ( \Omega; \R ) } 
\leq
C 
%(
%@\| X_0 \|@
%)
h^\frac32.
\end{equation}
In the same sprit of \eqref{eq:f-conditional-expect-32-model},
we rely on the use of It\^o formula applied to $f (x) = \alpha_{-1} x^{-1} - \alpha_0 + \alpha_1 x - \alpha_2 x^\kappa$
to show
\begin{equation}
\label{eq:Ri-conditional-expectation-non-critical}
\begin{split}
\lVert \E ( R_{ i } \vert \mathcal{ F }_{ t_{ i - 1 } } ) \rVert_{ L^2 ( \Omega; \R^d ) } 
& \leq 
\bigg \lVert 
%\E \bigg(  
\int_{ t_{i-1} }^{ t_{ i } } 
\E \big( 
[ f( X_s ) - f( X_{ t_{ i  } } ) ]
\vert \mathcal{ F }_{ t_{ i - 1} } 
\big)
\dd s 
%\bigg) 
\bigg \rVert_{ L^2 ( \Omega; \R^d ) } 
\\
& \leq
\int_{ t_{i-1} }^{ t_{ i } } 
\int_{s}^{t_i} 
\lVert
f' ( X_r ) f ( X_r ) + \tfrac12 f '' ( X_r ) g^2 ( X_r )
\rVert_{ L^2 ( \Omega; \R ) }
\, \dd r
\, \dd s
\\ &
 \leq 
C h^2,
\end{split}
\end{equation}
where we recalled that the It\^o integral vanishes under the conditional expectation 
and also used  the Jensen inequality and Lemma \ref{Lem:moment-bounds-general}.
Armed with these two estimates, one can apply Theorem \ref{thm:upper-error-bound} to arrive at the desired assertion.
\qed
%%%
%
%\begin{rem}
%In \cite{szpruch2011numerical}, the authors proved the strong convergence of the backward Euler method for the Ait-Sahalia model
%on the condition $  $
%\cite{chassagneux2016explicit,neuenkirch2014first,wang2018mean}
%\end{rem}

\vspace{0.3cm}
%\subsection{The critical case $ \kappa + 1 = 2 \rho $}
{\bf 5.2.2 The critical case $ \kappa + 1 = 2 \rho $}
\vspace{0.3cm}

In what follows we turn to the general critical case $ \kappa + 1 = 2 \rho $ and 
present first a lemma concerning (inverse) moment bounds 
of the solution process, which can be proved by following the same lines in the proof of 
Lemma \ref{Lem:moment-bounds-general} (cf. \cite[Lemma 2.1]{szpruch2011numerical}).
\begin{lem}\label{Lem:moment-bounds-critical}
Let conditions in Proposition \ref{thm:wellposed-Ait-Sahalia} be all fulfilled with $ \kappa + 1 = 2 \rho $
and let $ \{ X_t \}_{t \geq 0}$  be the unique solution to \eqref{eq:Ait-Sahalia-model-SDE}. 
Then we have, for any $ 2 \leq p_1 \leq \frac{ \sigma^2 + 2 \alpha_2 } { \sigma^2 }$ 
and for any $ p_2 \geq 2 $,
\begin{equation} \label{eq:moment-inverse-mom-critical}
\sup_{ t \in [ 0, \infty) } \E [ | X_t |^{p_1} ] < \infty, 
\quad
\sup_{ t \in [ 0, \infty) } \E [ | X_t |^{- p_2} ] < \infty.
\end{equation}
\end{lem}
For the purpose of analyzing the convergence rate of the numerical approximations, 
we validate all conditions of Assumptions \ref{ass:monotonicity-condition}, \ref{ass:well-possedness} in the next lemma, 
which is required by Theorem \ref{thm:upper-error-bound}.
\begin{lem}
\label{lem:AS-model-critical-assumptions}
Let conditions in Proposition \ref{thm:wellposed-Ait-Sahalia} be all fulfilled with $ \kappa + 1 = 2 \rho $
and let $ h \in (0, \tfrac{1}{ \alpha_1}] $. Let the model parameters obey 
$\tfrac{ \alpha_2 }{ \sigma^2 } \geq 2 \kappa - \tfrac32$ 
and $\tfrac{ \alpha_2 }{ \sigma^2 } > \tfrac { \kappa + 1 } { 2 \sqrt{2} }$. 
Then the SDE model \eqref{eq:Ait-Sahalia-model-SDE} and the scheme \eqref{eq:Milstein-Ait-Sahalia-model} satisfy 
Assumptions \ref{ass:monotonicity-condition}, \ref{ass:well-possedness} in the domain $D =  (0, \infty)$.
\end{lem}
{\it Proof of Lemma \ref{lem:AS-model-critical-assumptions}.}
Recall that the well-posedness of the model and the scheme in $ D = ( 0, \infty)$ has been proven in Proposition \ref{thm:wellposed-Ait-Sahalia} 
and Lemma \ref{lem:well-posed-Milstein-ASmodel}. It remains to verify the other conditions.
Thanks to the assumptions $ \kappa + 1 = 2 \rho $ and $\tfrac{ \alpha_2 }{ \sigma^2 } > \tfrac { \kappa + 1 } { 2 \sqrt{2} }$, 
one can find $\varrho > 1$ such that $ \alpha_2^2 > \tfrac{ \varrho }{ 2 } \rho^2 \sigma^4 
=  \tfrac{ \varrho }{ 8 } ( \kappa + 1 )^2 \sigma^4$ and thus
\begin{equation} \label{eq:mono-condition-critical}
\begin{split}
& h \tfrac{ \varrho }{ 2 } \| g'g ( x ) - g'g ( y ) \| ^2
  - h \| f ( x ) - f ( y ) \|^2
\\
& \quad
=
h
\big(
 \tfrac{ \varrho }{ 2 } \rho^2 \sigma^4  
%\big( \tfrac{ x^{ \kappa } - y^{ \kappa } }{ x - y}  \big)^2
z_{\kappa}^2
-
\big(
 \tfrac{ \alpha_{-1} }{x y} - \alpha_1 + \alpha_2 z_{\kappa}
\big)^2
\big)
(x - y)^2
\\ 
& \quad
=
h
\big[
\tfrac{ \varrho }{ 2 } \rho^2 \sigma^4  z_{\kappa}^2 
- 
( \tfrac{ \alpha_{-1} }{x y} - \alpha_1 )^2
-
\alpha_2^2 z_{\kappa}^2
-
2 ( \tfrac{ \alpha_{-1} }{x y} - \alpha_1 ) \alpha_2 z_{\kappa}
\big]
( x - y)^2
\\
&
\quad
\leq
h
\big[
(
\tfrac{ \varrho }{ 2 } \rho^2 \sigma^4 
-
\alpha_2^2 
)
z_{\kappa}^2
+
2 \alpha_1 \alpha_2  z_{\kappa}
\big]
( x - y)^2
\\
&
\quad
\leq C h ( x - y )^2,
\quad
\forall \, x, y \in \R_+.
\end{split}
\end{equation}
Noting $ \kappa + 1 = 2 \rho $ again, one can deduce from \eqref{eq:one-sided-non-critical} that
\begin{equation}
\label{eq:one-sided-critical}
\sup_{ x > 0}
\Big(
f'(x) + \tfrac{q-1}{2} | g' (x) |^2 
\Big)
%= - \alpha_{-1} x^{ - 2} + \alpha_1 -  \alpha_2 \kappa x^{ \kappa - 1 } 
%+ \tfrac{ (q-1) \sigma^2 \rho^2 }{ 2 } x^{ 2 \rho - 2 } 
\leq
\sup_{ x > 0}
\Big(
\alpha_1 -  \big( \alpha_2 \kappa 
- \tfrac{ (q-1) \sigma^2 \rho^2 }{ 2 } 
\big) x^{ \kappa - 1 }
\Big),
\quad
\forall \,
x \in \R_+.
\end{equation}
%Thanks to the assumption 
Since $ \tfrac{ \alpha_2 }{ \sigma^2 } \geq 2 \kappa - \tfrac32 > \tfrac{ \kappa + 3 }{ 8 } > \tfrac18 ( \kappa + 2 + \tfrac {1}{\kappa} )$
for $\kappa > 1$, one can find $q > 2$ such that $ \tfrac{ \alpha_2 }{ \sigma^2 } \geq  \tfrac{q-1}{8} ( \kappa + 2 + \tfrac {1}{\kappa} ) $, i.e.,
$ \alpha_2 \kappa - \tfrac{ (q-1) \sigma^2 \rho^2 }{ 2 } 
=  \alpha_2 \kappa - \tfrac{ (q-1) }{ 8 } \sigma^2 ( \kappa + 1 )^2 \geq 0 $ in \eqref{eq:one-sided-critical},
and thus, similarly to \eqref{eq:monocity-Ait-Sahalia-0}, 
\begin{equation}\label{eq:monocity-Ait-Sahalia-1}
\begin{split}
 \langle x - y, f ( x ) - f ( y ) \rangle + \tfrac {q-1} {2} \lVert g ( x ) - g ( y ) \rVert ^ 2 
 \leq
\alpha_1 ( x - y)^2,
\quad
\forall x, y \in \R_+.
\end{split}
\end{equation}
Combining this with \eqref{eq:mono-condition-critical} ensures that the condition 
\eqref{eq:mono-condition1} is fulfilled in 
$ D = ( 0, \infty)$ with $ \theta = 1, \eta = 0$.
The condition \eqref{eq:mono-condition2} follows from \eqref{eq:monocity-Ait-Sahalia-1} directly.
Finally,  since $ \frac{ \sigma^2 + 2 \alpha_2 } { \sigma^2 } \geq 4 \kappa - 2  > 2 \kappa$
by assumption $\tfrac{ \alpha_2 }{ \sigma^2 } \geq 2 \kappa - \tfrac32$, $ \kappa >1 $,
in view of Lemma \ref{Lem:moment-bounds-critical} one can infer $\sup_{ t \in [ 0, T] } \| X_t \|_{ L^2 ( \Omega; \R ) }  < \infty$ and
\begin{equation}
\sup_{ t \in [ 0, T] } \| f ( X_t ) \|_{ L^2 ( \Omega; \R ) }  
\leq
\sup_{ t \in [ 0, T] }
\big(
\alpha_{-1} \| X_t^{-1} \|_{ L^2 ( \Omega; \R ) }  
+
\alpha_0 
+ 
\alpha_1 \| X_t \|_{ L^2 ( \Omega; \R ) } 
+
\alpha_2 \| X_t \|^\kappa_{ L^{2\kappa} ( \Omega; \R ) } 
\big)
< \infty.
\end{equation}
%@MOMENT BOUNDS@
Therefore, all conditions in Assumptions \ref{ass:monotonicity-condition}, \ref{ass:well-possedness} are 
confirmed in the domain $D =  (0, \infty)$.
\qed

Now we are in a position to derive the convergence order with the aid of Theorem \ref{thm:upper-error-bound}.
\begin{thm}\label{thm:Milstein-Ait-Sahalia-convergence-critcal}
%
%{\color{red}{
Let $\{ X_{ t} \}_{ t \in [0, T]}$ and $ \{Y_n\}_{0\leq n\leq N} $ be solutions to 
\eqref{eq:Ait-Sahalia-model-SDE} and \eqref{eq:Milstein-Ait-Sahalia-model}, respectively.
%Let $\alpha_{-1}, \alpha_0$, $\alpha_1, \alpha_2, \sigma > 0$, $\kappa, \rho >1$ obey
%$ \kappa + 1 = 2 \rho $ and $ \tfrac{ \alpha_2 }{ \sigma^2 } \geq 2 \kappa - \tfrac32$ and 
Let conditions in Proposition \ref{thm:wellposed-Ait-Sahalia} be all fulfilled with $ \kappa + 1 = 2 \rho $
and let $ h \in (0, \tfrac{1}{2  \alpha_1} ) $.
Let the model parameters $\alpha_{-1}, \alpha_0, \alpha_1, \alpha_2, \sigma > 0$, $\kappa > 1, \rho >1$ obey 
$\tfrac{ \alpha_2 }{ \sigma^2 } \geq 2 \kappa - \tfrac32$ 
and $\tfrac{ \alpha_2 }{ \sigma^2 } > \tfrac { \kappa + 1 } { 2 \sqrt{2} }$. 
Then there exists a constant $C>0$, independent of $N \in \N$, such that
\begin{equation} \label{eq:thm-convergence-rate-critical}
%\sup_{ N \in \N }
\sup_{ 1\leq n \leq N}
\| X_{ t_n } - Y_n \|_{ L^2 ( \Omega; \R ) } 
\leq
C h.
\end{equation}
%}}
\end{thm}
{\it Proof of Theorem \ref{thm:Milstein-Ait-Sahalia-convergence-critcal}.}
As already verified in Lemma \ref{lem:AS-model-critical-assumptions}, 
all conditions in Assumptions \ref{ass:monotonicity-condition}, \ref{ass:well-possedness} are fulfilled in $D =  (0, \infty)$.
Based on Theorem \ref{thm:upper-error-bound}, one only needs to properly estimate $\lVert R_{ i } \rVert_{ L^2 ( \Omega; \R ) } $ 
and $\lVert \E ( R_{ i } \vert \mathcal{ F }_{ t_{ i - 1 } } ) \rVert_{ L^2 ( \Omega; \R^d ) }$.
Similarly as above, we split the the error term $\lVert R_{ i } \rVert_{ L^2 ( \Omega; \R ) } $ into two parts:
\begin{align}\label{eq:Ait-Sahal-model-Ri-criticle}
%\begin{split}
%\E \big[ 
\lVert R_{ i } \rVert_{ L^2 ( \Omega; \R ) } 
%\big]
\leq &
\bigg \lVert \int_{ t_{i-1} }^{ t_{ i } } f( X_s ) - f( X_{ t_{ i } } ) \, \dd s \bigg \rVert_{ L^2 ( \Omega; \R ) }
\nonumber
\\ & +
\bigg \lVert \int_{ t_{i - 1 } }^{ t_{ i } } g( X_s ) - g( X_{ t_{ i - 1 } } ) 
- 
g'g ( X_{ t_{i-1} } ) ( W_s - W_{ t_{i - 1 } } ) \, \dd W_s  
\bigg \rVert_{ L^2 ( \Omega; \R ) }
\nonumber \\
=: &
I_6 + I_7,
%\end{split}
\end{align}
where the coefficients $f,g$ are defined by \eqref{eq:Ait-Sahalia-f-g-defn}.
The It\^o formula applied to $f (x) = \alpha_{-1} x^{-1} - \alpha_0 + \alpha_1 x - \alpha_2 x^\kappa,
x \in \R_+$ gives
%and Lemma \ref{Lem:moment-bounds-critical} promise
\begin{equation}
\begin{split}
I_6 &  \leq
\int_{ t_{i-1} }^{ t_{ i } } \lVert f( X_s ) - f( X_{ t_{ i } } ) \rVert_{ L^2 ( \Omega; \R ) } \, \dd s 
\\
&
\leq
C h^\frac32 
\Big( 1 + \sup_{s \in [0, T]} \| X_s \|_{L^{4 \kappa - 2} ( \Omega; \R) }^{ 2 \kappa - 1 } 
+ 
\sup_{s \in [0, T]} \| X_s^{-1} \|_{L^{6} ( \Omega; \R) }^{ 3 }
\Big).
%\leq
%C h^\frac32.
\end{split}
\end{equation}
Following the same lines as in \eqref{eq:I5-estimate}, one can similarly show
\begin{equation}
| I_7 |^2
\leq
C h^3
\Big( 1 + \sup_{s \in [0, T]} \E [ | X_s |^{ 4 \kappa - 2 } ]
+ 
1_{\{ \rho < 2 \} }
\sup_{s \in [0, T]} \E[  | X_s |^{ - ( 4 - 2 \rho) }]
\Big),
%\big( 1 + \sup_{s \in [0, T]} \| X_0 \|_{L^? ( \Omega; \R) }^? \big).
\end{equation}
where we set $ 1_{\{ \rho < 2 \} } = 1 $ for $ \rho < 2 $ and $ 1_{\{ \rho < 2 \} } = 0 $ for $ \rho \geq 2 $.
Since $ \frac{ \sigma^2 + 2 \alpha_2 } { \sigma^2 } \geq 4 \kappa - 2 $ 
by the assumption $\tfrac{ \alpha_2 }{ \sigma^2 } \geq 2 \kappa - \tfrac32$,
we can plug these two estimates into \eqref{eq:Ait-Sahal-model-Ri-criticle}
and use Lemma \ref{Lem:moment-bounds-critical} 
to get
\begin{equation}
\lVert R_{ i } \rVert_{ L^2 ( \Omega; \R ) } 
\leq
C h^\frac32.
\end{equation}
Moreover, similarly to \eqref{eq:Ri-conditional-expectation-non-critical}, 
applying the It\^o formula to $f (x) = \alpha_{-1} x^{-1} - \alpha_0 + \alpha_1 x - \alpha_2 x^\kappa$
and noting the It\^o integral vanishes under the conditional expectation we deduce
\begin{equation}
%\label{eq:...}
\begin{split}
\lVert \E ( R_{ i } \vert \mathcal{ F }_{ t_{ i - 1 } } ) \rVert_{ L^2 ( \Omega; \R^d ) } 
& \leq
\bigg \lVert 
%\E \bigg(  
\int_{ t_{i-1} }^{ t_{ i } } 
\E \big( 
[ f( X_s ) - f( X_{ t_{ i  } } ) ]
\vert \mathcal{ F }_{ t_{ i - 1} } 
\big)
\dd s 
%\bigg) 
\bigg \rVert_{ L^2 ( \Omega; \R^d ) } 
\\
&
\leq
\int_{ t_{i-1} }^{ t_{ i } } 
\int_{s}^{t_i} 
\lVert
f' ( X_r ) f ( X_r ) + \tfrac12 f '' ( X_r ) g^2 ( X_r )
\rVert_{ L^2 ( \Omega; \R ) }
\, \dd r
\, \dd s
\\ 
& \leq  
C h^3 
\Big( 1 + \sup_{s \in [0, T]} \| X_s \|_{L^{4 \kappa - 2} ( \Omega; \R) }^{ 2 \kappa - 1 } 
+ 
\sup_{s \in [0, T]} \| X_s^{-1} \|_{L^{6} ( \Omega; \R) }^{ 3 }
\Big).
\end{split}
\end{equation}
In light of Lemma \ref{Lem:moment-bounds-critical} and with the help of Theorem \ref{thm:upper-error-bound}, 
one can obtain the assertion \eqref{eq:thm-convergence-rate-critical}.
\qed
%
%%%
\begin{rem} \label{rem:Ait-Sahalia-model}
Recall that Szpruch et al. \cite{szpruch2011numerical} examined the backward Euler method for
the Ait-Sahalia model \eqref{eq:Ait-Sahalia-model-SDE} and proved its strong convergence
only when $ \kappa + 1 > 2 \rho $,  but without revealing a rate of convergence. 
Very recently,  the authors of \cite{wang2018mean} fill the gap by
identifying the expected mean-square convergence rate of order $\tfrac12$ for stochastic theta methods
applied to the Ait-Sahalia model under conditions $ \kappa + 1 \geq 2 \rho $.
%
%
%The authors of \cite{neuenkirch2014first,wang2018mean}...
%
In 2014, a kind of Lamperti-backward Euler method was introduced in \cite{neuenkirch2014first} for the  
Ait-Sahalia model, with a mean-square convergence rate of order $1$ identified for 
the full parameter range in the general standard case $ \kappa + 1 > 2 \rho $ and
for a particular critical case $\kappa = 2, \rho = 1.5$  when $\tfrac{ \alpha_2 } { \sigma^2 } > 5 $  
(see Propositions 3.5, 3.6 from \cite{neuenkirch2014first}). As shown above, we apply the semi-implicit Milstein method 
\eqref{eq:Milstein-Ait-Sahalia-model} to the Ait-Sahalia model, which is able to treat both the general standard case
and a more general critical case $\kappa + 1 = 2 \rho$ for any $\kappa,  \rho> 1$. 
Moreover, we prove a mean-square convergence rate of order $1$ for the full parameter range in the general standard case and
for parameters satisfying $\tfrac{ \alpha_2 } { \sigma^2 } > 2 \kappa - \tfrac32 $ and 
$\tfrac{ \alpha_2 }{ \sigma^2 } > \tfrac { \kappa + 1 } { 2 \sqrt{2} }$ in  the general critical case.
For the special critical case $\kappa = 2, \rho = 1.5$, the restriction on parameters reduces into $\tfrac{ \alpha_2 } { \sigma^2 } > \tfrac52 $, 
which is moderately more relaxed than $\tfrac{ \alpha_2 } { \sigma^2 } > 5 $ as required in \cite{neuenkirch2014first}.
\end{rem}

%\subsection{The CIR and CEV model}

%\subsection{Applications with multi-level Monte Carlo approximations}

\subsection{Numerical tests}
\label{subsect:numer-results}

The aim of this subsection is to illustrate the above theoretical findings by providing several numerical examples. 
Two different schemes covered by \eqref{eq:general-scheme} are utilized to 
simulate the two previously studied financial models.
The resulting mean-square approximation errors are computed at the endpoint $T = 1$ and
the desired expectations are approximated by averages over $10000$ samples. Moreover,
the ``exact'' solutions are identified as numerical ones using a fine stepsize $h_{\text{exact}} = 2^{-12}$.

\begin{figure}[htp]
\centering
      \includegraphics[width=3.2in,height=3.5in]  {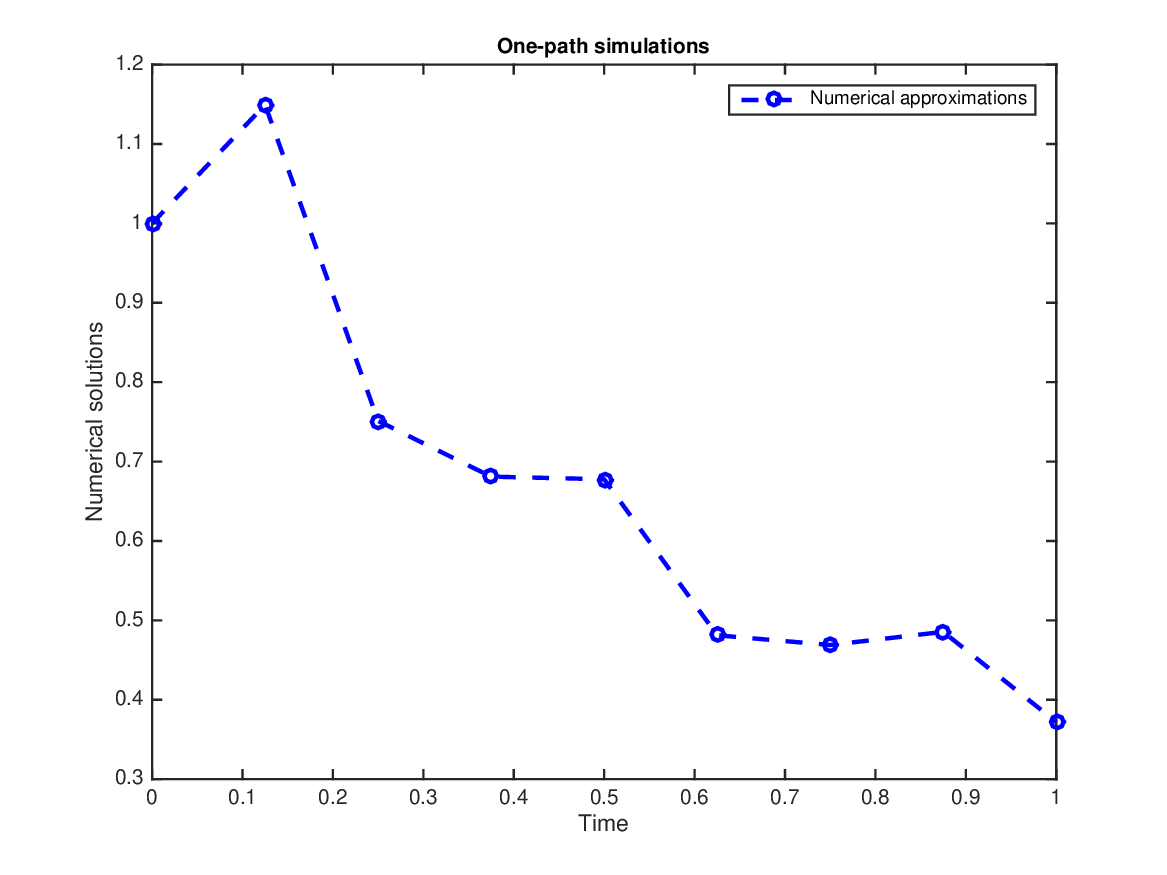}
      \includegraphics[width=3.2in,height=3.5in]  {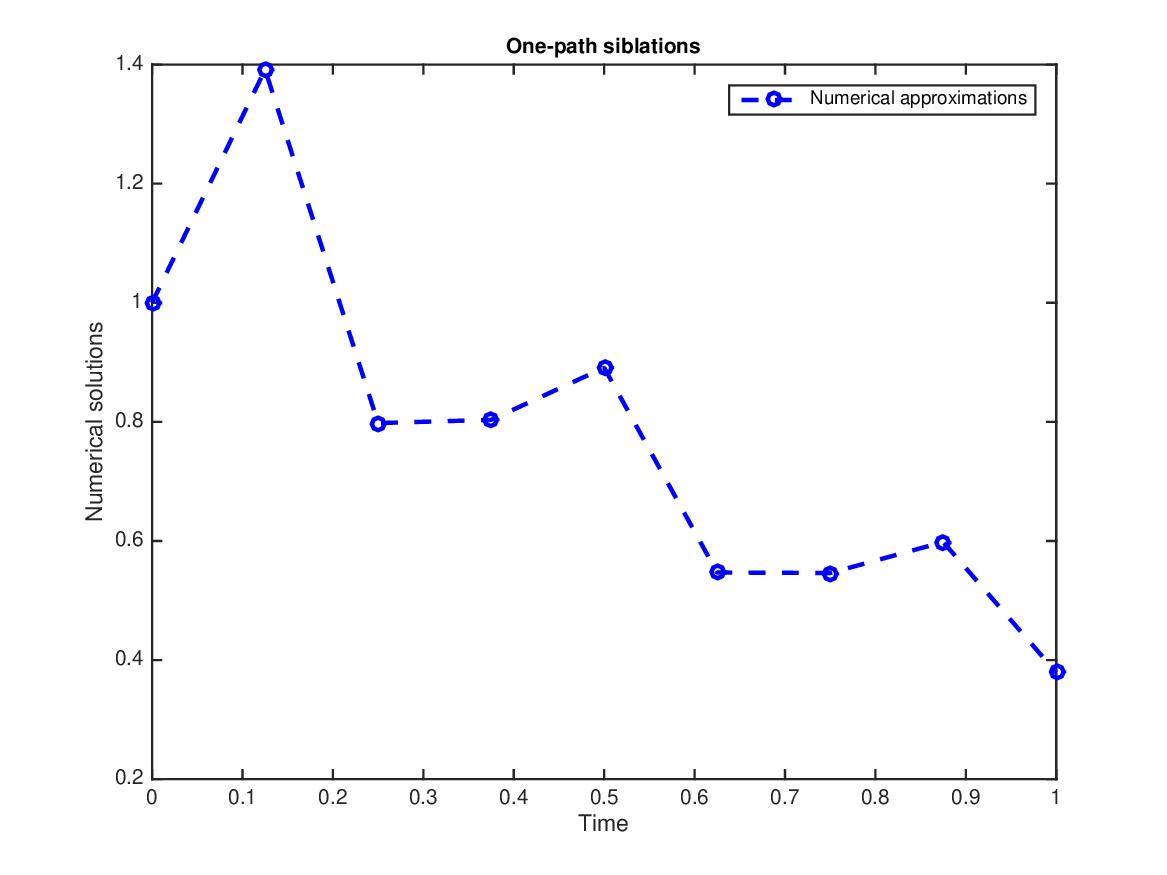}
 \caption{One-path simulations of the drift-diffusion double implicit Milstein method 
 for the Heston $\tfrac32$-volatility model (Left) and the stochastic LV model (Right).}
\label{fig:one-path-simulations}
\end{figure}

\begin{figure}[htp]
\centering
      \includegraphics[width=3.2in,height=3.5in]  {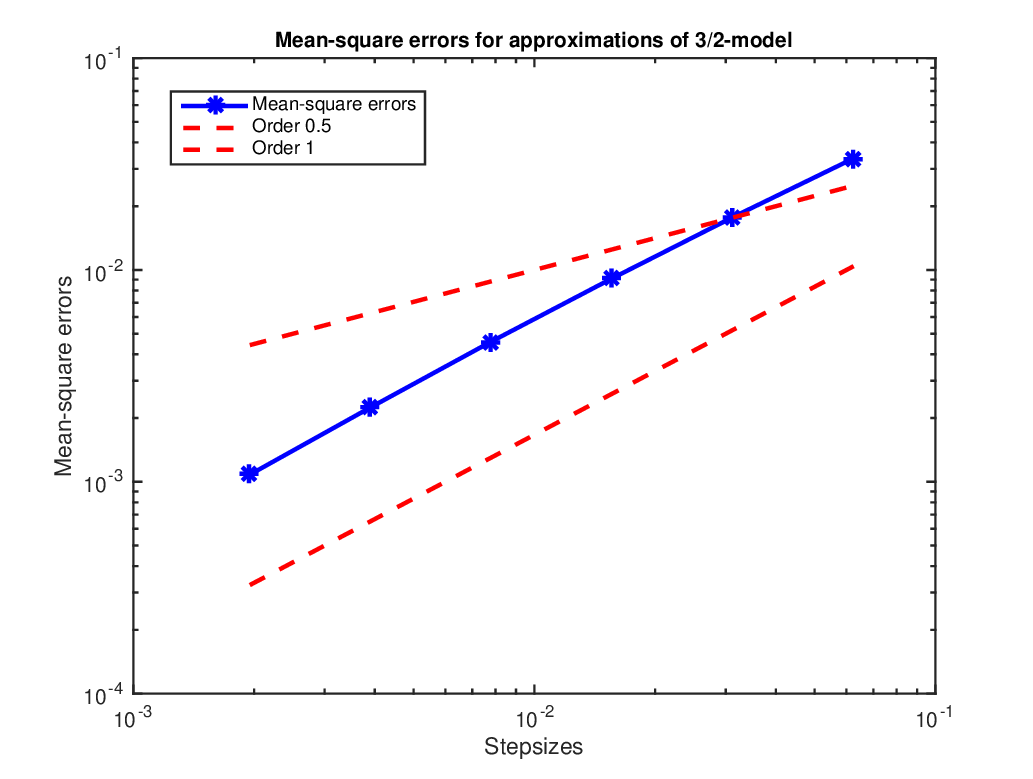}
      \includegraphics[width=3.2in,height=3.5in]  {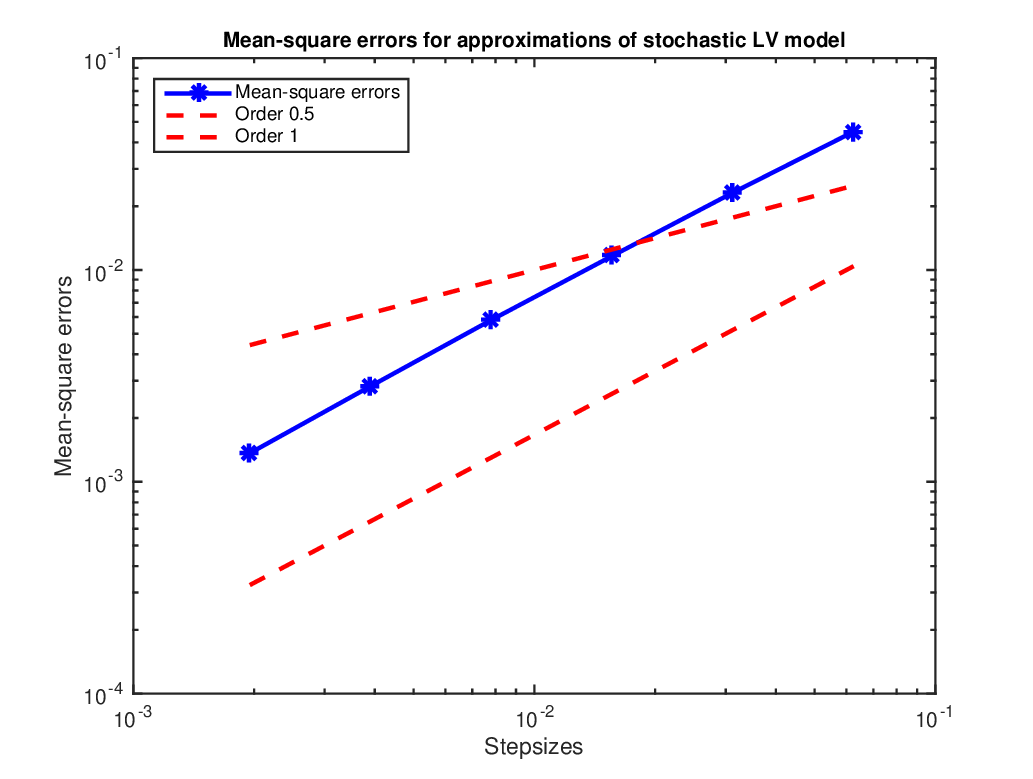}
 \caption{Mean-square convergence rates of the drift-diffusion double implicit Milstein method 
 for the Heston $\tfrac32$-volatility model (Left) and the stochastic LV model (Right).}
\label{fig:convergence-rates-32model}
\end{figure}

As the first example, let us first  look at the following SDE,
%Heston $\tfrac32$-volatility model,
{\color{black}
\begin{equation}
\label{eq:LV32model-numer-test}
d X_t = X_t ( \mu - \alpha X_t ) \, \dd t + ( \beta X_t ^ {3/2} + \sigma X_t ) \, \dd W(t),
 \quad 
 X_0 = 1,
 \quad
 t \in (0, 1].
\end{equation}
When $\sigma = 0 $ and $\beta = 0 $, the considered SDE \eqref{eq:LV32model-numer-test} reduces to 
the Heston $\tfrac32$-volatility model \eqref{eq:32model} and the stochastic LV competitive model \eqref{eq:LVmodel}, respectively.
We choose the parameters $ ( \mu, \alpha,  \beta, \sigma ) = ( 2,\tfrac52, 1, 0 )$ such that $ \alpha \geq \tfrac52 \beta^2 $ for 
the $\tfrac32$-model \eqref{eq:32model} and $ ( \mu, \alpha,  \beta, \sigma ) = ( 2,1, 0, 1)$ 
for the stochastic LV competitive model \eqref{eq:LVmodel}.   
By taking $ \theta = \eta = 1 $, we discrete these two models by the drift-diffusion double implicit Milstein method 
\eqref{eq:general-scheme}, which is explicitly solvable here (see \eqref{eq:model32-Milstein-explicit} and \eqref{eq:LVmodel-Milstein-explicit}).
%i.e., \eqref{eq:32model-scheme} or \eqref{eq:model32-Milstein-explicit}. 
In the following simulations, the expectations are approximated by computing averages over $10^4$ samples 
and the ``exact'' solutions are identified as approximations using a fine stepsize $h_{\text{exact}} = 2^{-12}$.
It turns out that the resulting numerical approximations always remain positive for all $10^4$ paths.
In Figure \ref{fig:one-path-simulations}, we present one-path simulations of the drift-diffusion double implicit Milstein method 
 for the Heston $\tfrac32$-volatility model (Left) and the stochastic LV model (Right), which are shown to be positive. 
}
To test the mean-square convergence rates, 
we depict in Figure \ref{fig:convergence-rates-32model} mean-square approximation errors $e_{h}$ 
%$e_{h} : = \sqrt{\E [ \| X_T - Y_N \|^2 ]}, T = 1$, at the endpoint $T = 1$
against six different stepsizes $h = 2^{-i}, i = 4, 5,...,9$ on a log-log scale.  Also, two reference lines of
slope $1$ and $\tfrac12$ are given there. From Figure \ref{fig:convergence-rates-32model}  
one can easily detect that the approximation errors decrease at a slope close to $1$ when stepsizes shrink, coinciding 
with the predicted convergence order obtained in Theorem \ref{thm:Milstein-convergence-rate-32model} 
and Theorem \ref{thm:Milstein-convergence-rate-LVmodel}.
Suppose that  the approximation errors $e_{h}$ obey a power law relation 
$e_{h} = C h^\delta $ for $C, \delta >0$, so that $\log e_{h} = \log C +  \delta \log h$.  Then we do a least squares power law 
fit for $\delta$ and get the value $0.9923$ for the rate $\delta$ with residual of $0.0719$. 
Again, this confirms the expected convergence rate in Theorem \ref{thm:Milstein-convergence-rate-32model}
and Theorem \ref{thm:Milstein-convergence-rate-LVmodel}.

%{\color{black}Assigning $\beta = 0 $, the considered SDE \eqref{eq:LV32model-numer-test}  
%reduces to the stochastic LV competitive model \eqref{eq:LVmodel}.
%Further, we choose the parameters $ ( \mu, \alpha,  \sigma ) = ( 2,1, 1)$
%and  discrete the model by the drift-diffusion double implicit Milstein method 
%\eqref{eq:LVmodel-scheme} or \eqref{eq:LVmodel-Milstein-explicit}. }
%

As the second example model, we look at the Ait-Sahalia interest rate model, given by
\begin{equation}
\label{eq:Ait-Sahalia-model-numer-test}
\dd X_t = ( \alpha_{-1} X_t^{-1} - \alpha_0 + \alpha_1 X_t - \alpha_2 X_t ^\kappa  ) \, \dd t + \sigma  X_t ^{\rho} \, \dd W_t,
\quad
X_0 = 1,
 \quad
 t \in (0, 1].
\end{equation}
%where $X_t, t \geq 0$ always remain positive 
%and we equivalently write  $|X_t|^\kappa$ and $| X_t |^{\rho}$ to 
%make the explicit Euler method  (not positivity preserving) work well. 
Let us consider both  the standard case $ \kappa + 1 > 2 \rho $ and the critical case $ \kappa + 1 = 2 \rho $,
by taking two sets of model parameters:
\newline

\noindent $\bullet$ Case I:  $\kappa = 4, \, \rho = 2, \,  \alpha_{-1} = \tfrac32, \alpha_0 = 2, \alpha_1 = 1, \alpha_2 = 1, \sigma = 1; $
\vspace{0.5cm}

\noindent $\bullet$ Case II:   $\kappa = 3, \, \rho = 2, \, \alpha_{-1} = \tfrac32, \alpha_0 = 2, \alpha_1 = 1,  \alpha_2 = \tfrac92, \sigma = 1 $.
\newline

%\begin{table}[!ht]
%\label{table:least-squares}
%\begin{center}% \footnotesize
%\caption{A least squares fit for the convergence rate $\delta$.} 
%\begin{tabular*}{12cm}{@{\extracolsep{\fill}}ccc}
%\hline  & Case I  &  Case II
% \\ \hline
%$\theta=\frac12$ & $\delta$ = 0.5588, resid = 0.0554  & $\delta$ = 0.5571, resid = 0.0499
%\\ \hline
%$\theta=1$ & $\delta$ = 0.5554, resid = 0.0546 & $\delta$ = 0.5633, resid = 0.0396   \\
% \hline
%\end{tabular*}
%\end{center}
%\end{table}

\begin{figure}
\centering
      \includegraphics[width=3.2in,height=3.5in]{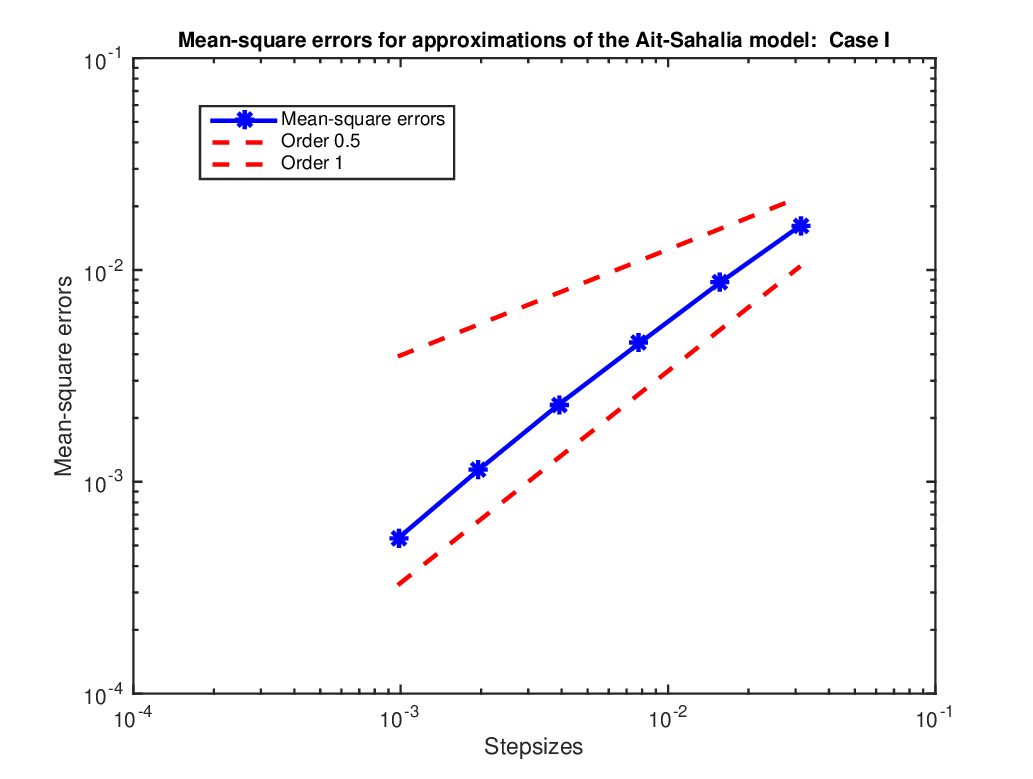}
      \includegraphics[width=3.2in,height=3.5in]{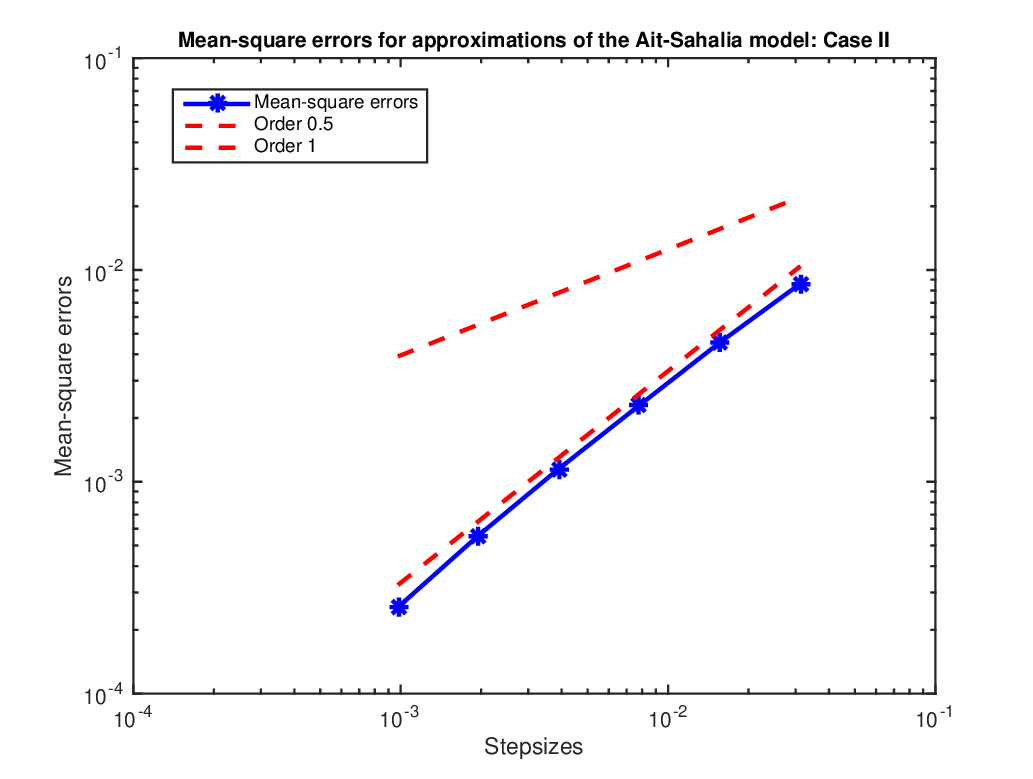}
 \caption{Mean-square convergence rates of the semi-implicit Milstein method 
 \eqref{eq:Milstein-Ait-Sahalia-model} for the Ait-Sahalia interest rate model (Left for Case I and right for Case II).}
\label{fig:convergence-rates-Ait-Sahalia-model}
\end{figure}

It is easy to check that Case I corresponds to the standard case and Case II corresponds to 
the critical case $ \kappa + 1 = 2 \rho $ satisfying $ \tfrac{ \alpha_2 }{ \sigma^2 } \geq 2 \kappa - \tfrac32$ 
and $\tfrac{ \alpha_2 }{ \sigma^2 } > \tfrac { \kappa + 1 } { 2 \sqrt{2} }$.
The semi-implicit Milstein scheme \eqref{eq:Milstein-Ait-Sahalia-model} is used to simulate the model 
\eqref{eq:Ait-Sahalia-model-numer-test} for these two cases.
As shown in Figure \ref{fig:convergence-rates-Ait-Sahalia-model}, the mean-square approximation error lines
have slopes close to $1$ for both cases.  A least squares fit produces a rate $ 0.9798$ with residual
of $0.0929$ for Case I  and a rate $1.0129$  with residual of $0.0968$ for Case II.
Hence, numerical results are consistent with strong order of convergence equal to one, 
as already revealed in Theorem \ref{thm:Mistein-Ait-Sahalia-convergence-non-critcal} 
and Theorem \ref{thm:Milstein-Ait-Sahalia-convergence-critcal}.

\section{Conclusion}
\label{sect:conclusion}
The present work introduces a family of implicit Milstein type methods for strong approximations of stochastic differential equations (SDEs)
with non-globally Lipschitz drift and diffusion coefficients. An easy and direct approach of the error analysis is developed 
to recover the expected mean-square convergence rate of order one for the proposed schemes. 
In particular, the optimal convergence rate of the positivity preserving schemes applied to 
three models in practice is obtained for the first time and more relaxed conditions are required, compared 
with existing results for first order schemes in the literature. 
In the future, we attempt to identify the general  $L^p$ rate of convergence with $p \geq 2$ for the schemes, 
which is highly non-trivial. 
%as the techniques of this paper are not enough.
\\
\newline
{\bf Funding }  This work was supported by Natural Science Foundation of China (12071488, 11971488)
                and Natural Science Foundation of Hunan Province for Distinguished Young Scholars (2020JJ2040). 

%\bibliographystyle{abbrv}

%\bibliography{../bib/bibfile}
%\bibliography{bibfile}
%
%%%%%

\section*{Declarations}

{\bf Conflict of interest}  The author declares no competing interests.

%\bibliographystyle{abbrv}

%\bibliography{bibfile}

\end{document}